\documentclass{ws-ijbc}     
\usepackage{graphicx}

\def\Z{\mathbb{Z}}

\def\E{\mathcal E}
\def\C{\mathcal C}
\def\p{\mathcal P}
\def\H{\mathcal H}
\def\M{\mathcal M}
\def\X{\mathcal X}
\let\e=\epsilon
\let\d=\delta
\let\D=\Delta

\newcommand{\be}{\begin{equation}}    
\newcommand{\ee}{\end{equation}}
\newcommand{\ba}{\begin{eqnarray}}
\newcommand{\ea}{\end{eqnarray}}
\newtheorem{oss}{Remark}
\begin{document}

\catchline{}{}{}{}{} 

\markboth{A. Marchesiello \& G. Pucacco}{Bifurcation sequences in the 1:1 Hamiltonian resonance}

\title{Bifurcation sequences in the symmetric 1:1 Hamiltonian resonance}

\author{Antonella Marchesiello}
\address{Faculty of Nuclear Sciences and Physical Engineering, \\
Czech Technical University in Prague, \\
D\v{e}\v{c}\'{\i}n Branch,
Pohranicn\'{\i} 1, 40501  D\v{e}\v{c}\'{\i}n \\
e-mail: anto.marchesiello@gmail.com}

\author{Giuseppe Pucacco}
\address{Dipartimento di Fisica and INFN -- Sezione di Roma II, \\
Universit\`a di Roma ``Tor Vergata", \\
Via della Ricerca Scientifica, 1 - 00133 Roma\\
e-mail: pucacco@roma2.infn.it}

\maketitle

%

\begin{abstract}
We present a general review of the bifurcation sequences of periodic orbits in general position of a family of resonant Hamiltonian normal forms with nearly equal unperturbed frequencies, invariant under $\Z_2\times\Z_2$ symmetry. The rich structure of these classical systems is investigated with geometric methods and the relation with the singularity theory approach is also highlighted. The geometric approach is the most straightforward way to obtain a general picture of the phase-space dynamics of the family as is defined by a complete subset in the space of control parameters complying with the symmetry constraint. It is shown how to find an energy-momentum map describing the phase space structure of each member of the family, a catastrophe map that captures its global features and formal expressions for action-angle variables. Several examples, mainly taken from astrodynamics, are used as applications.

\end{abstract}


\section{Introduction}
\label{intro}

Resonant Hamiltonian normal forms is a field of utmost interest for applications in a wide range of fields. Just to mention a few of them: celestial mechanics and astrodynamics, galactic dynamics, molecular and semi-classical quantum physics, structural engineering, etc. At the same time it is still the subject of active studies in analytical mechanics since it gives many clues for understanding non-integrable dynamics by the combination of different techniques such as group theory, algebraic and differential geometry, singularity theory, etc. For the readers involved in applications is not easy to appreciate the ability of these methods in characterizing the regular dynamics of non-integrable systems. Aim of this review is to provide an account of how to obtain a deep and `useful' framework to analyze a specific class of systems by exploiting a good deal of advanced tools, but avoiding a too formal apparatus and relying only on basic knowledges of analysis and geometry. Where strictly necessary, more sophisticated tools are introduced and described in a self-consistent way.

We investigate a class of two degrees of freedom Hamiltonian dynamical systems which at the same time have interesting applications and limited complexity. They are given by a planar harmonic oscillator with unperturbed frequencies close to each other and a non-linear coupling acting as a `perturbation' respecting some symmetries. These are chosen in order to simplify the description of the system without giving up the capability of representing good models for applications.

The nonlinearly perturbed harmonic oscillator is a reference system for everybody working in applied mathematics with an eye to physics, astronomy, chemistry, engineering, etc. \cite{LL,c1,fsgj,tv}, not to mention the investigation of effects due to even the tiniest perturbing forces \cite{torsion}. 

Among low-order Hamiltonian resonances (see e.g. Cushman {\it et al.}, 2007), the 1:1 resonance plays a prominent role. A huge amount of work has been devoted to this study leading to  advances that almost covered the subject. We recall the works of Kummer (1976), Deprit and coworkers \cite{D1,DE,M1}, Cushman and coworkers \cite{CR}, Broer and coworkers \cite{Br1:1} and van der Meer (2009). The general treatment of the non-symmetric 1:1 resonance has been done by Cotter \cite{Cotter} in his PhD thesis. With motivations mainly coming from galactic dynamics \cite{P09,MP11,MP13a}, the present review covers the most general case of a detuned 1:1--resonant Hamiltonian normal form invariant under $\Z_2 \times \Z_2$ symmetry by considering its generic expression (`universal unfolding') with three parameters plus detuning \cite{Henrard,Sch,Vf}. Although the treatment in Kummer's work \cite{Ku} is general enough to accommodate for detuning-like terms, their analysis is not explicit neither in his work nor in the others cited above. Moreover, bifurcation sequences in terms of the distinguished parameter (the `energy'), which are useful when comparing with numerical or laboratory experiments, are not explicitly given in the available references. Therefore we provide a general description of the phase-space structure as it is determined both by the `internal' parameters (energy and detuning) and the external (or `control') ones \cite{VU11,PHD}. As case study we will detail a natural system with potential associated to cubic and quartic non-linear coupling.

We exploit threshold values for bifurcations of periodic orbits as a latch to unlock the general structure of phase-space. This approach is based on the use of a regular reduction \cite{CB,KE} dividing out the $\mathbb{S}^1$ symmetry of the normal form. The reduced phase-space is invariant with respect to a second $\Z_2$ symmetry: we use a second reduction introduced by Han{\ss}mann and Sommer [2001] which allows us to divide out one of the $\Z_2$ symmetries. This trick provides an effective geometric strategy to understand how the phase-space structure is shaped by all possible combinations of the parameters. 

A general exploration of the space of external control parameter is possible by first examining a reference set and then analyzing its complement by exploiting the symmetries. In this way a complete description of both generic and degenerate cases is made possible, shedding light on the structural stability of the former ones and how to break the degeneracy of the latter. As a coronation of the geometric approach, a two-parameter combination (the `catastrophe' map, \cite{STK}) allows us to represent the general setting in a suitable 2-plane and all possible bifurcation sequences are clearly represented in the plane of the values of integrals of motion, the energy-momentum map, that can be easily plotted and used to get information on relative fractions of phase-space volume pertaining to each stable family. Quadrature formulas for the actions, periods and rotation number can also be obtained. 

In order to show the use of the theoretical approach, a set of examples is illustrated in a final section devoted to the applications. These are mainly coming from astrodynamics with reference to galactic dynamics and celestial mechanics; an example of a very well known degenerate case is also provided.


\section{The normal form}
\label{NOFO}

Let us consider a two-degrees of freedom system whose Hamiltonian is an analytic function in a neighborhood of an elliptic equilibrium. Its series expansion about the equilibrium point can be written as
\begin{equation}\label{hamiltonian}
H(\bm{p},\bm{q})=\sum_{j=0}^{\infty}H_{j}(\bm{p},\bm{q}),
\end{equation}
where $H_{j}(\bm{p},\bm{q})$, $j\geq0$, are homogeneous polynomials of degree $j+2$ in the canonical variables. Let us take the zero-order term ${H}_0$ to be in the positive definite form
\be\label{HZO}
H_0(\bm{p},\bm{q})=\frac12 \omega_1 (p_1^2 + q_1^2) + \frac12 \omega_2 (p_2^2 + q_2^2)
\ee
and assume that the `unperturbed' frequencies satisfy a quasi-resonance condition of the form
\be\label{det11}
\omega_1=\left(1+\d \right)\omega_2,
\ee
where  $\d \in \mathbb R$ is a `detuning' parameter \cite{Henrard,Sch,Vf}. With the rescaling of time
\be\label{scaling_t}
t\rightarrow \omega_2t
\ee
we can split the quadratic part into the isotropic oscillator Hamiltonian
\be\label{HZ}
H_0(\bm{p},\bm{q})=\frac12 (p_1^2 + q_1^2 + p_2^2 + q_2^2)
\ee
and the small quadratic term
\be\label{dett}
\frac12 \d (p_1^2+q_1^2);\ee
this is considered of higher-order and treated as a perturbation. We will investigate only time-reversible systems. Moreover, by means of translations, rotations and rescaling, it is always possible to suppress a set of higher-order terms. Therefore we assume that \eqref{hamiltonian} is invariant with respect to
\begin{eqnarray}
\mathcal S&:&(p_1,p_2,q_1,q_2)\rightarrow(p_1,-p_2,q_1,-q_2)\label{spatial_symmetry}\\
\mathcal T&:&(p_1,p_2,q_1,q_2)\rightarrow(-p_1,-p_2,q_1,q_2).\label{time_symmetry}
\end{eqnarray}
A standard example is the natural case in which
\begin{eqnarray}
H_{1}&=&b_{30} q_1^3 + b_{12} q_1 q_2^2 , \label{cubic_natural}\\
H_{2}&=&b_{40} q_1^4 + b_{22} q_1^2 q_2^2 + b_{04} q_2^4 \label{quartic_natural}
\end{eqnarray}
and so forth: the $b_{jk}$ are a set of real parameters specifying the original model problem. With the above mentioned transformations we see that, e.g. in $H_1$, the terms $b_{21} q_1^2 q_2$ and $b_{03} q_2^3$ have been suppressed.

The first step to simplify our system is the procedure of `normalization'. Resonant normal forms for the Hamiltonian system corresponding to \eqref{hamiltonian} are constructed with standard methods. The usual approach in Hamiltonian perturbation theory is to find a generating function in such a way to construct a `simpler' Hamiltonian: the first successful algorithm devoted to a resonant system was precisely introduced to treat the problem of the `third integral' in galactic dynamics \cite{gu} and was applied to construct the normal form of the {\it H\'enon-Heiles} model \cite{hh}. In subsections \ref{compl} and \ref{degenz} we will come back on the peculiarities of this fundamental example. The method based on the Lie transform offers several conceptual and technical advantages \cite{bp2,gior} with respect to the original approach and is nowadays the first choice \cite{mho,CHR,PY}. 

The idea behind the method is that of seeing the canonical transformation as a `flow' along the Hamiltonian vector field associated to the generating function. In spite of this apparently exotic statement, it turns out that this approach is that best suited in implementing recursive algorithms working with series of functions. Let us consider a sequence of polynomial functions in the phase-space 
\be\label{gene}
\{ G \} =G_{1},G_{2},...\ee 
To each term of the sequence (\ref{gene}) is naturally associated the linear differential operator
\be\label{diff}
{\rm e}^{{\cal L}_G} = \sum_k \frac1{k!} {\cal L}_G^k,\ee
whose action on a generic function $F$ is given by the Poisson bracket:
\be
{\cal L}_G F \doteq \{F,G\}.\ee
The original Hamiltonian system \eqref{hamiltonian} undergoes a canonical transformation to new variables $
(\bm{P},\bm{Q}),$ such that 
\be\label{TBS}
\bm{p} = {\rm e}^{{\cal L}_G} \bm{P}, \quad
\bm{q} = {\rm e}^{{\cal L}_G} \bm{Q}
\ee
and the new Hamiltonian is
\begin{equation}\label{HK}
     K(\bm{P},\bm{Q})=\left( {\rm e}^{{\cal L}_G} H \right) (\bm{P},\bm{Q})= H({\rm e}^{{\cal L}_G} \bm{P},{\rm e}^{{\cal L}_G} \bm{Q}),\,\ee
    where every function is assumed to be in the form of power series.

To construct $K$ starting from $H$ is a recursive procedure exploiting an algorithm based on the Lie transform \cite{gior,CHR}. To proceed we have to make some decision about the structure the new Hamiltonian must have, that is we have to chose a {\it normal} form for it. We therefore select the new Hamiltonian in such a way that it admit a new integral of motion, that is we take a certain function, say $F$, and impose that
\be\label{NFD}
\{K,F\}=0.
\ee
The usual choice (but not necessarily the only possible) is that of taking
\be\label{UC}
F=H_{0}=K_0
\ee
so that the function (\ref{HZ}) plays the double role of determining the specific form of the transformation and assuming the status of the second integral of motion. After $N$ steps of normalization, we get a truncated series of the form
\begin{equation}\label{kamiltonian}
K(\bm{P},\bm{Q};\alpha_i)=\sum_{j=0}^{N}K_{j}(\bm{P},\bm{Q}),
\end{equation}
each term of which satisfies condition \eqref{NFD}.
$K$ is characterized by a set of `external' control parameters (to be distinguished from the `internal' parameters fixed by the dynamics) that we collectively denote with $\alpha_i$. They are certain non-linear combinations of the parameters of the original physical model; for natural Hamiltonians there is a one-to-one correspondence between the original parameters and those in the normal form. Formally, a more direct way of applying this method is by using `action-angle'--{\it like} variables $\bm{J}, \bm{\phi}$ defined through the transformation
\be\label{Vaa}
Q_{\ell}=\sqrt{2J_{\ell}}\cos\phi_{\ell},\;\;\;\;\;P_{\ell}=\sqrt{2J_{\ell}}\sin\phi_{\ell}, \quad \ell=1,2,
\ee 
In this way we have
\be
H_0 = J_1 + J_2, \ee
so that
\be\label{LHAA}
 {\cal L}_{H_0} = \frac{\partial}{\partial \phi_1}  + \frac{\partial}{\partial \phi_2} . \ee
A generic polynomial series turns out to be a Poisson series in the angles, that is a Fourier series with coefficients depending on the actions only in such a way to comply with the polynomial structure of the terms of each degree \cite{mho}. It is clear that in the non-resonant case, the set of terms in normal form with $H_0$, namely those vanishing under the action of the operator \eqref{LHAA}, is composed only of powers of the actions: in this case we speak of a Birkhoff (or non-resonant) normal form. But, if the frequencies are almost in the ratio 1:1, we choose to retain also `resonant' terms, that is trigonometric terms in $\phi_1-\phi_2$. 
The result is \cite{MRS1,MP11}
\be\label{GNF}
K = J_1 + J_2 + \d J_1 + \alpha_1 J_1^2+ \alpha_2 J_2^2+ J_1 J_2 \left[\alpha_3 + 2 \alpha_4 \cos 2 (\phi_1-\phi_2)\right].
\ee 
As an example, the coefficients of the normal form of the standard natural system (\ref{cubic_natural},\ref{quartic_natural}) are given by
\ba
\alpha_1 & = & \frac32 \left( b_{40} - \frac52 b_{30}^2 \right), \label{11acq}\\
\alpha_2 & = & \frac12 \left( 3b_{04} - \frac56 b_{12}^2 \right), \label{11bcq}\\ 
\alpha_3 & = & b_{22} - b_{12} \left( 3 b_{30} + \frac23 b_{12} \right) ,\label{11ccq}\\
\alpha_4 & = & \frac14 \left[b_{22} + b_{12} \left( b_{30} - 2 b_{12} \right) \right]. \label{11dcq}\ea
It is easy to check that \eqref{GNF} is the most general form of a time reversible phase-space function of degree 2 in the actions which commutes with $H_0$: actually, we see that in this case $K_1$ vanishes (all terms coming from $H_1$ are conveyed in $G_1$) and $K$ appears to be invariant under the  the time reversion symmetry \eqref{time_symmetry} and the $\Z_2\times\Z_2$ group 
$$\{\rm{Id},\mathcal S_1,\mathcal S_2,\mathcal S_1 \circ \mathcal S_2\}$$ 
where $\mathcal S_2 \doteq \mathcal S$,
\be
\mathcal S_1:(P_1,P_2,Q_1,Q_2)\rightarrow(-P_1,P_2,-Q_1,Q_2).\label{spatial_symmetry1}
\ee
In this case the structure of the normal form imposes an extra symmetry which is apparent from the expression
\ba
K&=&\frac12 (1+\d)(P_1^2+Q_1^2) + \frac12 (P_2^2+Q_2^2) + 
\frac14 \alpha_1 (P_1^2+Q_1^2)^2 + \frac14 \alpha_2 (P_2^2+Q_2^2)^2 \\ 
&+&  \frac14 \alpha_3 (P_1^2+Q_1^2)(P_2^2+Q_2^2) +
\frac12 \alpha_4 [4 P_1 Q_1 P_2 Q_2 + (P_1^2-Q_1^2)(P_2^2-Q_2^2)].
\ea
To take into account the presence of reflection symmetries,
we can also speak of a `2:2 resonance'. The information given by the parameters of the original system is not simply scrambled in the $\alpha_i$ (see e.g. \ref{11acq}--\ref{11dcq}): the generating function $G$, which is an outcome of the normalization process together with $K$, can be used to invert the transformation and to recover apparently lost information.
In view of its peculiar role we include $\d$ in the category of internal (or `distinguished' parameters) \cite{Br1:1} and consider $\d J_1$ as a higher-order term with respect to $K_0$. We observe that the $\alpha_i$'s may in turn depend on $\d$ (as it happens, for example, in the family of natural systems with elliptical equipotentials \cite{MP13a}). The higher-order terms $K_{2j}(J_{\ell},\phi_{\ell}),j>1,$ are homogeneous polynomials of degree $2j$ in $J_{\ell}$ depending on the angles only through the combination $2 (\phi_1-\phi_2)$ and integer multiples of it. 

\section{Further simplification of the normal form}
\label{FUNOFO}

The canonical variables $J_{\ell},\phi_{\ell}$ are the most natural to investigate the dynamics in a perturbative framework. However, several other coordinate systems can be used to unveil the aspects of this class of systems. 
We list those that will be useful in the following.

First of all we use  coordinates `adapted to the resonance' \cite{SV}. There are various ways to do this: a standard recipe is to use the canonical transformation \cite{Br1:2}
\be \label{suber}
\left\{
  \begin{array}{l}
    J_1=J \\
    J_2=\E-J \\
     \psi=\phi_1-\phi_2\\
    \chi=\phi_2.
  \end{array}
\right.
\ee
This is used to perform a first reduction of the normal form, since $\chi$ is cyclic and its conjugate action $\E$ is the additional integral of motion. To first order, the reduced Hamiltonian is
\be\label{RK}
{\mathcal K} = \E + \alpha_2\E^2+\left(\d - (2 \alpha_2+\alpha_3)\E\right)J+(\alpha_1 +\alpha_2-\alpha_3)J^2 + \alpha_4 J(\E-J) \cos 2 \psi.\ee
A further reduction into a planar system, viewing $\E$
as a \emph{distinguished} parameter \cite{Br1:2} can be obtained via the canonical transformation \cite{Ku}
\be\label{ccoord}
\left\{
  \begin{array}{l}
    x=\sqrt{2J}\cos\psi \\
    y=\sqrt{2J}\sin\psi.
  \end{array}
\right.
\ee
 These coordinates are usually exploited to put ${\mathcal K}$ in the most `universal' form (see the Appendix). 
Following \cite{CB}, a related path to reduce the symmetry of the normal form passes through the introduction of the invariants of the isotropic harmonic oscillator:
\begin{equation}\label{invariants}
\left\{\begin{array}{ll}
         I_0= & \frac12(P_1^2+P_2^2+Q_1^2+Q_2^2)= K_0 =\E\\
         I_1= & P_1P_2+Q_1Q_2 \\
         I_2= & Q_1P_2-Q_2P_1 \\
         I_3= & \frac12(P_1^2-P_2^2+Q_1^2-Q_2^2).
       \end{array}
\right.
\end{equation}
The set $\{ I_0,I_1,I_2,I_3\}$ form a Hilbert basis of the ring of invariant polynomials and can be used as coordinates system for the reduced phase space. Their Poisson brackets are given by $\{I_a,I_b\}=2\e_{abc} I_c, \; a,b,c=1,2,3$. Notice that $I_0$ coincides with the linear part of the normal form $K_0=\E$, a Casimir of the Poisson structure. There is one relation between the new coordinates, namely
\begin{equation} \label{rel_invariants}
I_1^2+I_2^2+I_3^2=I_0^2=\E^2,
\end{equation}
hence the sphere
\begin{equation}\label{phase_sphere}
\mathcal S=\left\{(I_1,I_2,I_3)\in\mathbb R^3\;:\;I_1^2+I_2^2+I_3^2=\E^2\right\}
\end{equation}
is invariant under the flow defined by \eqref{GNF}. This provides a (geometric) second reduction to a one degree
of freedom system. The links between the two sets are given by the `Lissajous' relations
\cite{D1,DE}
\ba
I_1 &=& 2 \sqrt{J_1 J_2} \cos\psi = 2\sqrt{ J (\E-J)} \cos\psi,\label{inv1}\\
I_2 &=& 2 \sqrt{J_1 J_2} \sin\psi = 2 \sqrt{J (\E-J)} \sin\psi\label{inv2}
\ea
and
\be\label{stereo}
x=\frac{I_1}{\sqrt{\E-I_3}},\quad y=\frac{I_2}{\sqrt{\E-I_3}}.\ee
We remark that the coordinates $x,y$ are shown by Kummer [1976] to be associated with a variant of the stereographic projection of $\mathcal S$ on the $(I_1,I_2)$-plane.

The `normal modes' of the system are expressed in the following forms:
\ba
{\rm NM1}:&&J=0,\quad I_1 = I_2 = 0, \quad I_3 = -\E, \label{NM1}\\
{\rm NM2}:&&J=\E,\quad I_1 = I_2 = 0, \quad I_3 = \E. \label{NM2}
\ea
The periodic orbits `in general position' are most simply derived from the fixed points of the Hamiltonian vector field associated with \eqref{RK}. The family of `inclined' periodic orbits corresponds to the oscillations with phase difference $\psi=0,\pi$,
\ba
{\rm Ia}:&\psi= 0, &I_2 = 0, \quad I_3 = I_{3U}, \quad I_1 = + \sqrt{\E^2 - I_{3U}^2}, \label{Ia}\\
{\rm Ib}:&\psi= \pi, &I_2 = 0, \quad I_3 = I_{3U}, \quad I_1 = - \sqrt{\E^2 - I_{3U}^2}, \label{Ib}
\ea
whereas the family of `loop' periodic orbits corresponds to the oscillations  with phase difference $\psi=\pm\pi/2$,
\ba
{\rm La}:&\psi= \frac{\pi}2, &I_1 = 0, \quad I_3 = I_{3L}, \quad I_2 = + \sqrt{\E^2 - I_{3L}^2}, \label{La}\\
{\rm Lb}:&\psi= -\frac{\pi}2, &I_1 = 0, \quad I_3 = I_{3L}, \quad I_2 = - \sqrt{\E^2 - I_{3L}^2}. \label{Lb}
\ea
The expressions of $I_{3U}$ and $I_{3L}$ can be found by solving the conditions for the fixed points of the flow and will be recovered in Section \ref{GR} by relying on geometric arguments. 

We recall that in the unperturbed isotropic oscillator {\it all} orbits are periodic. A perturbation of the kind considered here always admits non-linear normal modes and `select' the families (\ref{Ia}--\ref{Lb}) among the set of all periodic orbits.
We also remark that the reduction illustrated above on the basis of the invariants of the isotropic oscillator can be in full analogy performed on the perturbations of the other paradigmatic superintegrable system, the Kepler problem, via Delaunay normalization \cite{D82}.

\section{Geometric reduction}\label{GR}

Exploiting the transformation \eqref{stereo} in order to use the invariant polynomials as phase-space variables, we can therefore adopt the function
\begin{equation}\label{hamiltonian_inv}
{\mathcal K}_I(I_1,I_2,I_3; \E)=\left(1+\Delta\right)\E + A_1 \E^2 +(B\E+\Delta)I_3 + A I_3^2 + C (I_1^2-I_2^2)
\end{equation}
on the reduced phase space given by the sphere \eqref{phase_sphere} to study the general behavior of the family. For simplicity we put
\be\label{ABC}
A\doteq\frac{\alpha_1+\alpha_2 - \alpha_3}{4},\;\;\;
     B\doteq\frac{\alpha_1-\alpha_2}{2},\;\;\;
C\doteq\frac {\alpha_4}{2},\;\;\; \Delta\doteq\frac{\d}2
\ee
and
$$
A_1 \doteq \frac{\alpha_1+\alpha_2 + \alpha_3}{4}.$$
There is a certain degree of redundancy in the external parameters, however as we see below there is  no strict reason not to keep them all, so we perform a general analysis of \eqref{hamiltonian_inv} for arbitrary values of the external parameters $A,B,C$ and the internal parameters $\Delta$ and $\E$. 

\subsection{Reduced phase space}
\label{RPS}

The two reflection symmetries now turn into the reversing symmetries $I_1\rightarrow-I_1$
and $I_2\rightarrow-I_2$. Their composition $(I_1,I_2,I_3)\rightarrow(-I_1,-I_2,I_3)$ gives a (non-reversing)
discrete symmetry of \eqref{hamiltonian_inv}. We perform a further reduction introduced by Han{\ss}mann and Sommer \cite{HS} to explicitly divide out this symmetry. This is given by the transformation
\begin{equation}\label{tr_lem}
\left\{
  \begin{array}{ll}
    X= I_1^2-I_2^2 \\
    Y= 2I_1I_2 \\
    Z= I_3
  \end{array}
\right.
\end{equation}
which turns the sphere \eqref{phase_sphere} into the `lemon' space
\begin{equation}\label{lemon}
\mathcal L=\left\{ (X,Y,Z)\in\mathbb R^3\;:\;X^2+Y^2=\left(\E+Z\right)^2\left(\E-Z\right)^2\right\}
\end{equation}
with Poisson bracket
$$\{f,g\}\doteq\left(\nabla f\times\nabla g,\nabla L\right)$$
where $(.,.)$ denotes the inner product and 
$$L\doteq X^2+Y^2-\left(\E+Z\right)^2\left(\E-Z\right)^2.$$ 
The Poisson structure of the $X,Y,Z$ variables is:
\begin{equation}\label{poistru}
\left\{
  \begin{array}{ll}
    \{X,Y\} = 8 \sqrt{X^2+Y^2} Z \\
    \{Y,Z\} = 4X \\
    \{Z,X\} = 4Y.
  \end{array}
\right.
\end{equation}
The Hamiltonian becomes
\begin{equation}\label{hamiltonian_lem}
{\mathcal K}_I(X,Z)=\left(1+\Delta\right)\E + A_1 \E^2 + C X+(B\E+\Delta)Z + A Z^2.
\end{equation}
We note that both reversing symmetries of \eqref{hamiltonian_inv} gives the invariance of \eqref{hamiltonian_lem}
with respect to the reflection $Y\rightarrow-Y$. To simplify the following formulae we omit the constant term from \eqref{hamiltonian_lem} by introducing 
\begin{equation}\label{eneconv}
\H \doteq {\mathcal K}_I - \left(1+\Delta\right)\E - A_1 \E^2 .
\end{equation} 
In this way we finally obtain
\begin{equation}\label{ham_lem_simply}
\H(X,Z)= C X+(B\E+\Delta)Z + A Z^2.
\end{equation}
Each integral curve of the reduced system defined by \eqref{ham_lem_simply} is given by the intersection between $\mathcal L$ and the surface
\begin{equation} \label{Ham_surface}
\{(Z,X)\in\mathbb R^2\;:\;\mathcal \H=h\}.
\end{equation}
Tangency points give equilibrium solutions. All information about bifurcations of periodic orbits in generic position and stability/transition of normal modes of the original system can be obtained by the study of the mutual positions of the surfaces $\H$ and $\mathcal L$  \cite{HS}.

\subsection{Reduced dynamics}
\label{RDY}

We can further simplify the approach by exploiting the fact that, since $Y$ does not enter in \eqref{ham_lem_simply}, the level sets $\{\H=h\}$ are parabolic cylinders.
For $C=0$ they degenerate into a pair of planes both perpendicular to the $Z$-axis. Notice that, at the quartic level,  $C=0$ corresponds to the integrable case of two uncoupled non-linear oscillators. For $A=0$ another degeneration occurs, since the level sets \eqref{Ham_surface}  reduce to the hyperplanes
\begin{equation}
CX +(B\E+\Delta)Z-h=0.
\end{equation}
We will refer to these particular cases as \emph{degenerate cases} and we will examine them in subsection \ref{dege}.

For $A\neq0,C\neq0$, if a tangency point occurs between  $\mathcal L$ and the surface \eqref{Ham_surface}, we have an (isolated) equilibrium for the reduced system.
Moreover,  two (degenerate) equilibria are represented by the singular points $\mathcal Q_1\equiv\left(0,0,-\E\right)$
and $\mathcal Q_2\equiv\left(0,0,\E\right)$. A tangent plane to $\mathcal L$ may coincide with a tangent plane to the parabolic cylinder $\{\mathcal K=h\}$
only at points where $Y$ vanishes. Hence all equilibria on $\mathcal L$  occur at $Y=0$: in order to study the existence and nature of the equilibria configuration of the system, it is then enough to  restrict the analysis to the phase-space section $\{Y=0\}$.

The contour $\mathcal C\equiv\mathcal L\cap\{Y=0\}$ in the $(Z,X)$-plane
is given by $\mathcal C_-\cup\mathcal C_+$, where
\be\label{lemon_arcs}
\C_\pm\equiv\left\{(Z,X)\in\mathbb R^2\;:\;|Z|\le\E,\;\,X=\pm\left(\E^2-Z^2\right)\right\}
\ee
and the set $\p\equiv\{\H=h\}\cap\{Y=0\}$ corresponds to the parabola
\begin{equation}\label{parabola}
X=\frac1{C} \left(h -(B\E+\Delta)Z-A Z^2 \right)\doteq \X(Z).
\end{equation}
All significant informations on the phase-space structure are therefore extracted from the analysis of the mutual positions of a parabola and two parabolic arcs on a plane. Given a certain set of `physical' parameters ($A,B,C,\Delta$) and some $\E>0$, a family of parabolas vertically `slides' on the $(Z,X)$-plane by varying $h$. Meaningful dynamics do exist in the interval of $h$ going from the first to the last contacts between the family of parabolas and the `projected lemon' \eqref{lemon_arcs}: in a sequence of pictures displayed in figg. \ref{RC1A} and \ref{RC1B} we will see all possibilities.  To enlighten the phase-space structure, each geometric frame of the reduced sections is matched with Poincar\`e sections $Q_2=0,P_2>0$ in the `physical' coordinates.

The reduced phase space $\mathcal L$ is invariant under reflection symmetries with respect to every coordinate-axes.
In particular, the reduced phase section $\C$ is invariant under both reflection transformations 
\begin{eqnarray}
R_1:\; Z\rightarrow -Z \label{R_z} \\
R_2:\; X\rightarrow -X \label{R_x}
\end{eqnarray}
and their compositions $R_1\circ R_2$ and $ R_2 \circ R_1$.  However the dynamics of the reduced system are not invariant under these actions.  Anyway it is easy to understand how they operate on the parabola $\eqref{parabola}$. When acting on $\X$, $R_1$ turns it into its symmetric with respect to the $X$-axis. Under the action of $  R_2$, $\X$ is reflected with respect to the $Z$-axis, that is, it reverses its concavity. Finally, the composition $ R_2  \circ R_1$ inverts the concavity of the parabola and then reflects it with respect to the $X$-axis. The application  of $ R_1\circ R_2$ on $\X$ gives the same result. Thus, we can restrict our analysis  to the case in which the parabola \eqref{parabola} is upward concave and for $\E=0$ achieves its minimum point on the negative $Z-$axis. If we choose a negative detuning, this corresponds to consider $A<0$ and $C>0$. Here and in the following we refer to this case as the \emph{reference case}. Then, by a simple application of \eqref{R_z}, $\eqref{R_x}$ and/or their compositions we obtain the bifurcation sequences in the remaining cases (cfr. the left panel in table 1).

\begin{table}

\tbl{Starting from the reference case \ref{caso_particolare}, we obtain all the complementary cases using reflection symmetries $ R_1$, $ R_2$ and $ R_1\circ R_2$. $\mathcal I$ stands for the identity transformation. The right panel shows how the fixed points of the system change under the action of the reflection symmetries of the twice reduced phase space.}
{\begin{tabular}{c|ll}

                                  & $A<0$& $A>0$   \\ \hline
                                  &       &           \\
  $C>0$                           & $\mathcal I$ &$R_2\circ R_1$\\
  $C<0$                           & $ R_2$ &$ R_1$         \\
  \hline
\end{tabular}
  \qquad\qquad
\begin{tabular}{c|llll}

                                  & $\mathcal Q_1$ & $\mathcal Q_2$ &$ \mathcal Q_L$       & $ {\mathcal Q}_U$       \\ \hline
                                  &       &       &             &              \\
  $ R_1$                  & $\mathcal Q_2$ & $\mathcal Q_1$ &$\widetilde {\mathcal Q}_L$ & $\widetilde {\mathcal Q}_U$ \\
  $ R_2$                  & $\mathcal Q_1$ & $\mathcal Q_2$ &${\mathcal Q}_U$        & $ {\mathcal Q}_L$ \\
  $ R_1\circ R_2$    & $\mathcal Q_2$ & $\mathcal Q_1$ &$\widetilde {\mathcal Q}_U$ & $\widetilde {\mathcal Q}_L$ \\
  \hline
\end{tabular}}

\label{T1}
\end{table}

On the section $\C$, the two degenerate equilibria are $\mathcal Q_1\equiv(-\E,0)$ and $\mathcal Q_2\equiv(\E,0)$.
It is always possible to fix $h$ such that \eqref{parabola} intersects $\C$ in one of these points, so that
\ba
h&=&h_1\doteq \E \left((A - B) \E - \Delta \right),\label{yax_energy}\\
h&=&h_2\doteq \E \left((A + B) \E + \Delta \right). \label{xax_energy}
\ea
 Thus, for $h=h_1$ the system stays in the point $\mathcal Q_1$
 and similarly for $h=h_2$. Comparing with \eqref{NM1}-\eqref{NM2} we see that they correspond to the two \emph{normal mode} solutions NM1 and NM2. 

\begin{figure}[!h]
\begin{center}
\includegraphics[width=7cm]{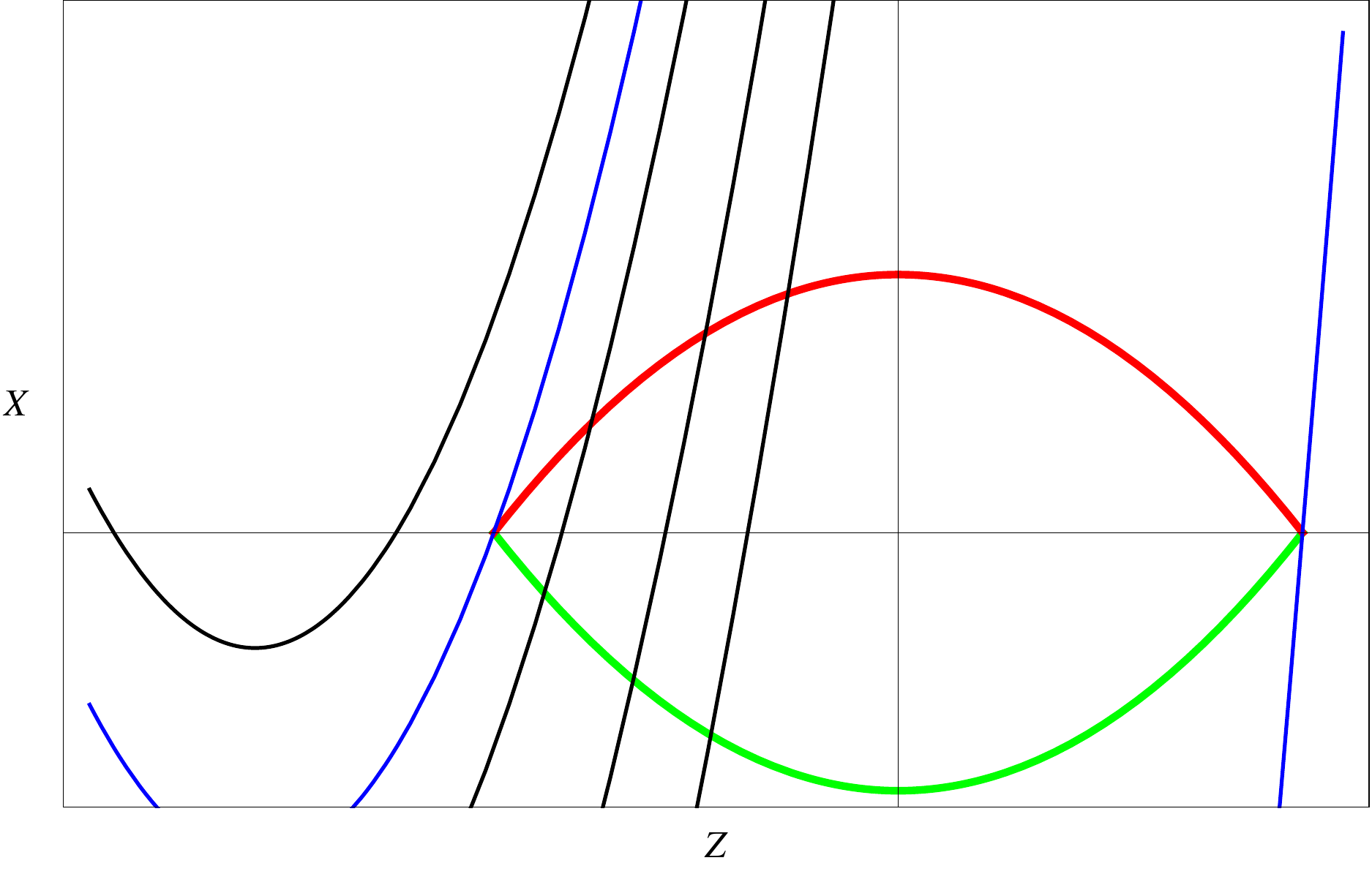}
\includegraphics[width=5.25cm]{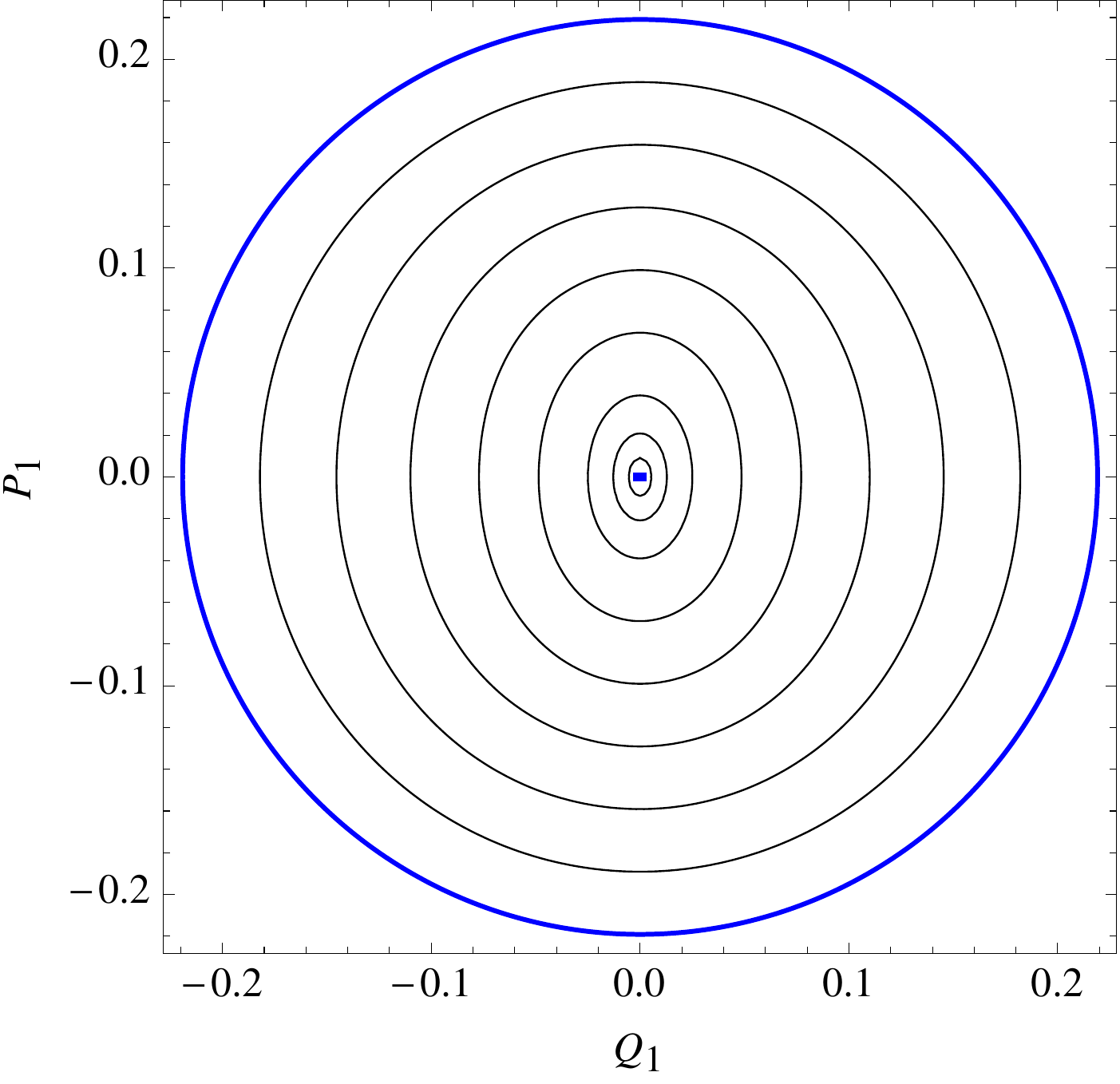}
\includegraphics[width=7cm]{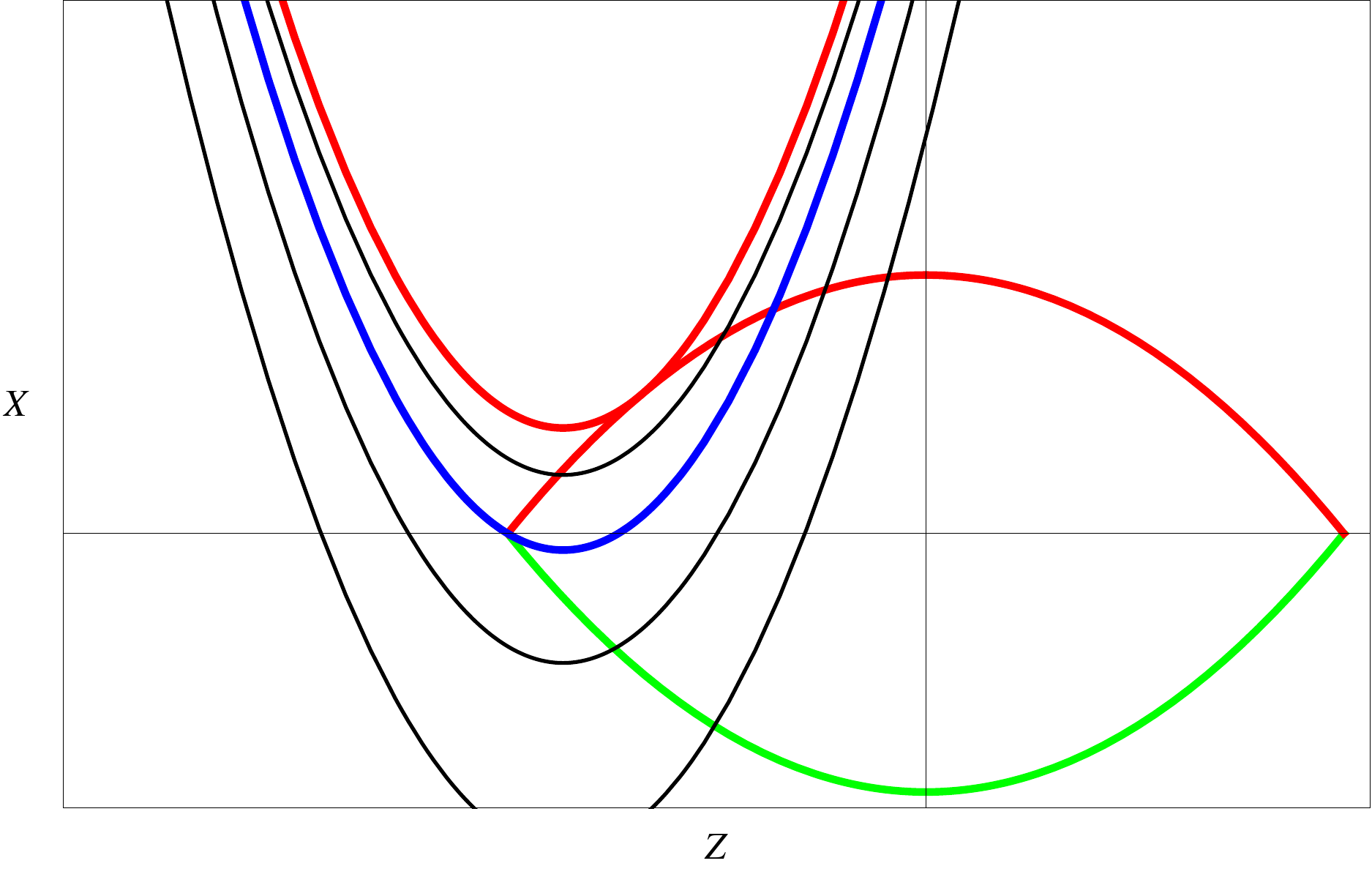}
\includegraphics[width=5.25cm]{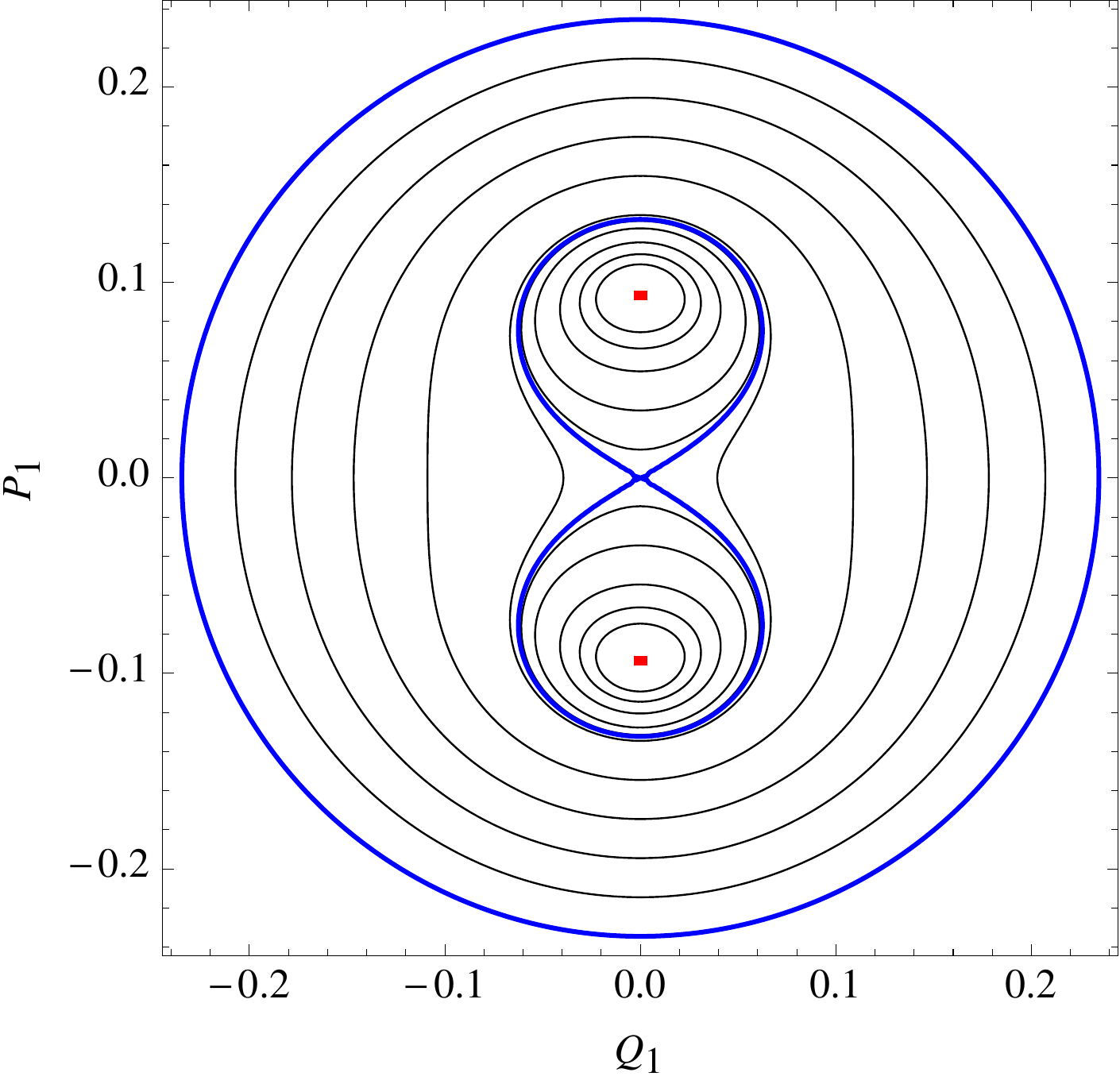}
\includegraphics[width=7cm]{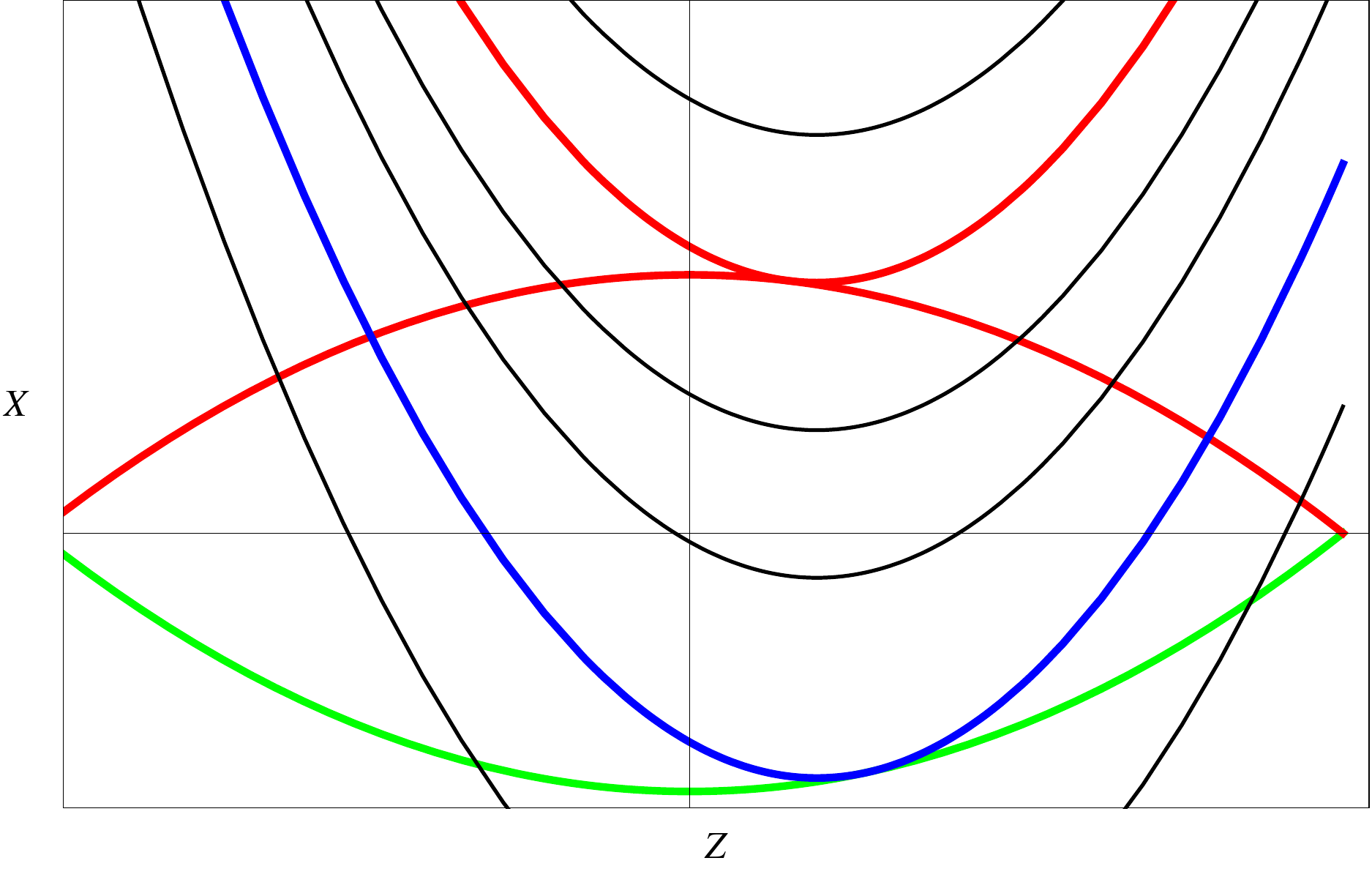}
\includegraphics[width=5.25cm]{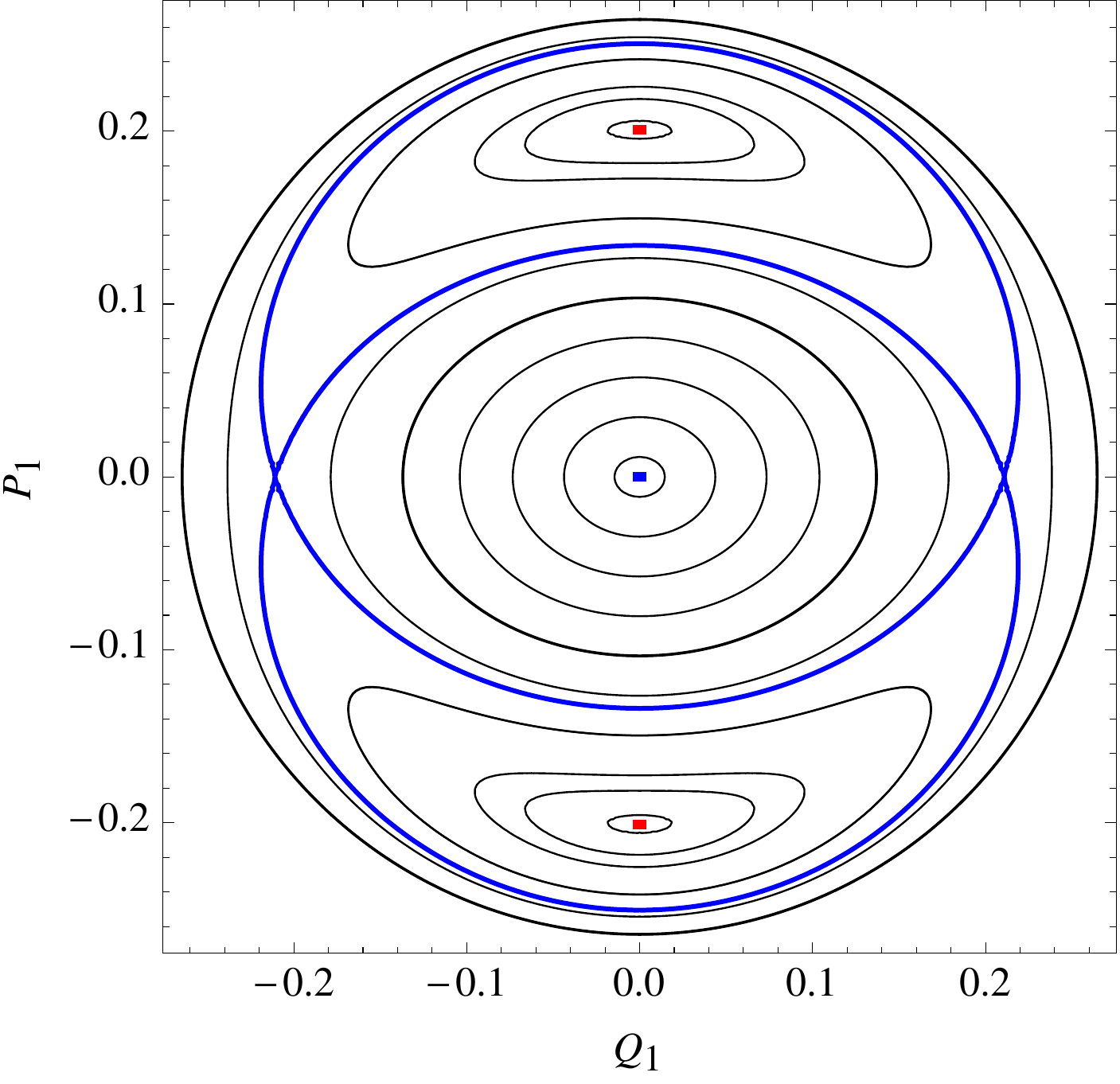}
\includegraphics[width=7cm]{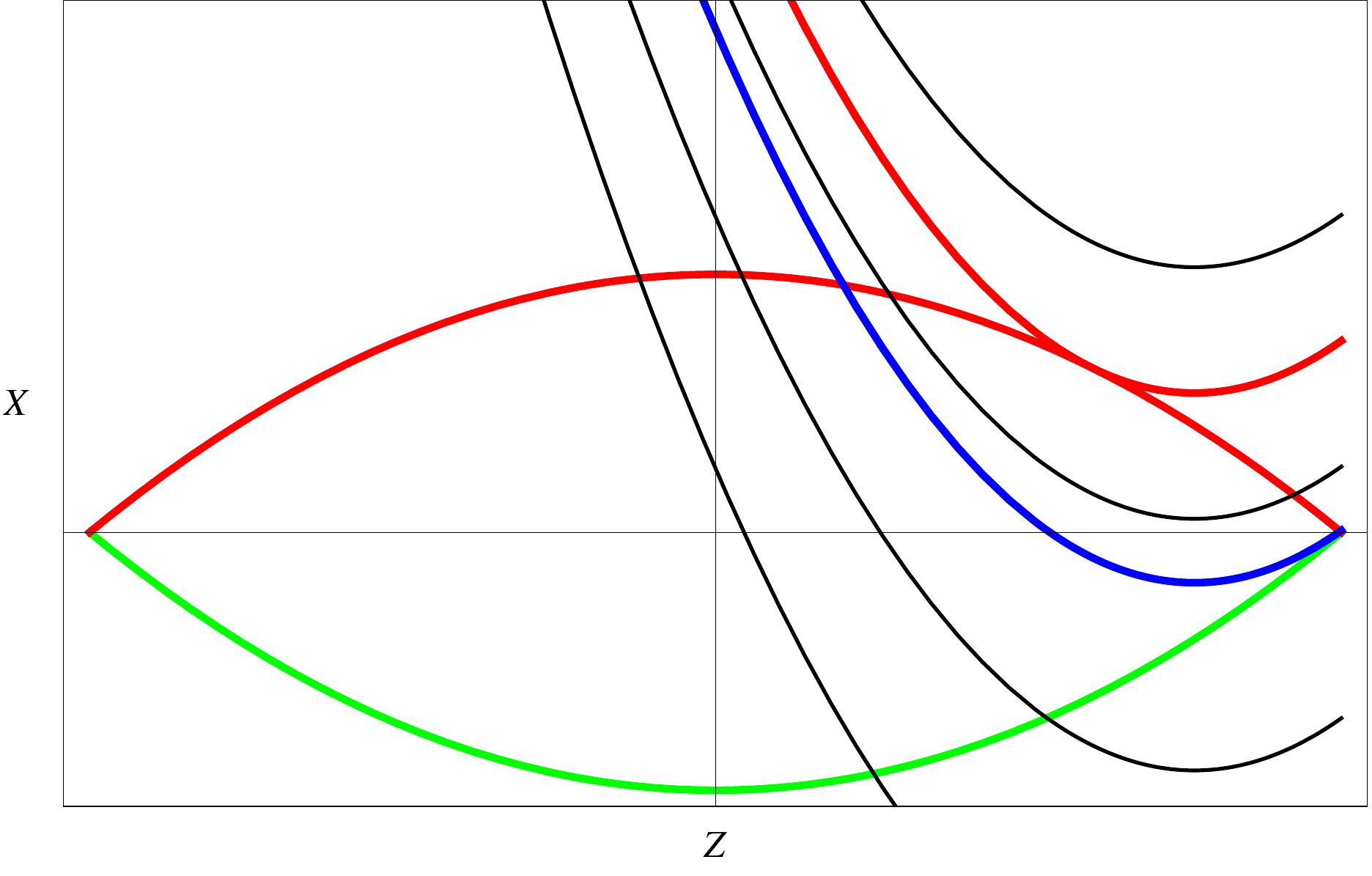}
\includegraphics[width=5.25cm]{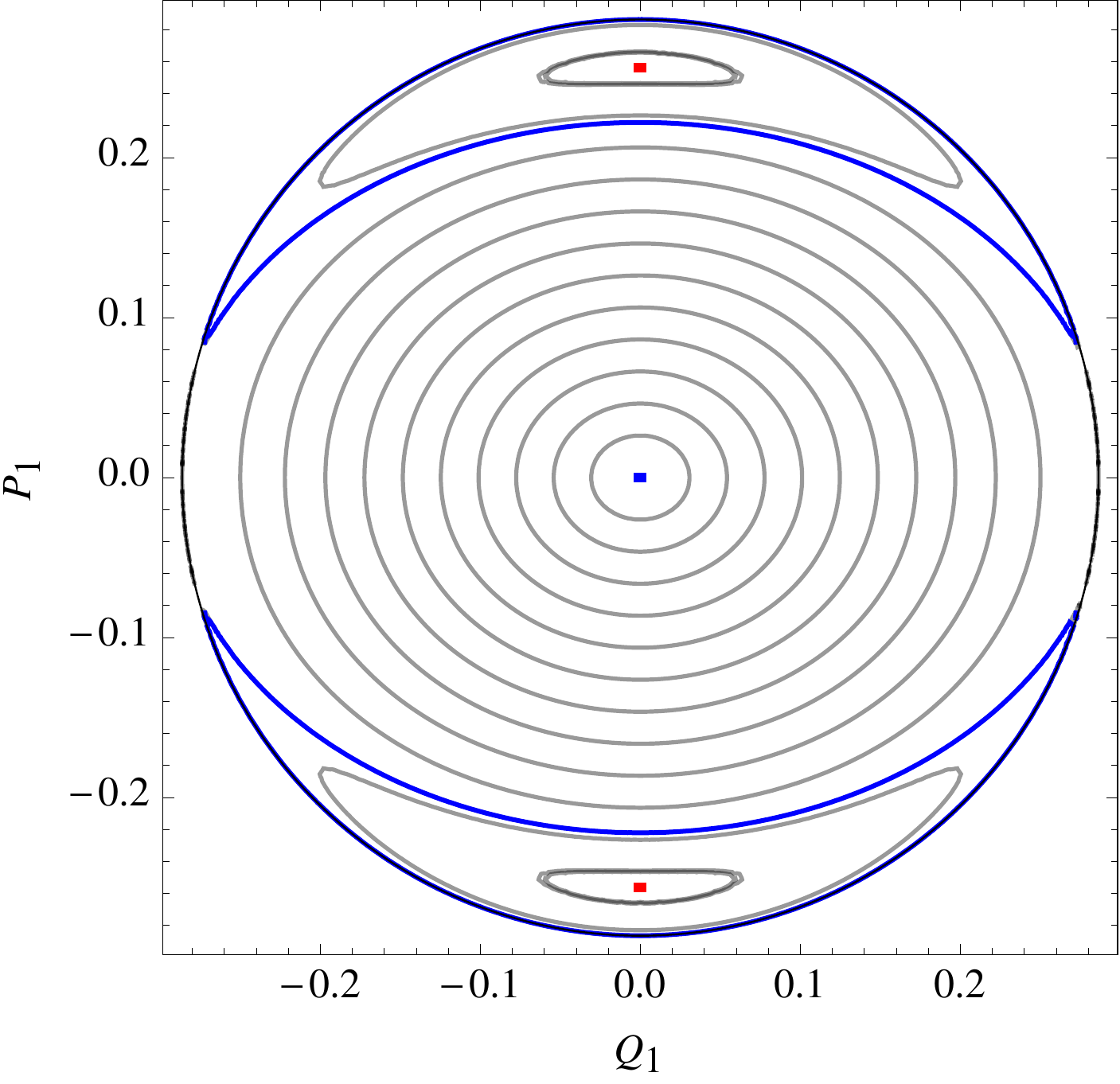}
\end{center}
\caption{Possible intersections between the reduced Hamiltonian and phase-space and corresponding surfaces of section: Reference Case ($A<0,C>0$, sub-case 1), $A=-\frac{11}{15},B=6,C=\frac15,\Delta=-\frac15$.}
\label{RC1A} 
\end{figure}

A stability/instability transition of a normal mode is generally associated with the bifurcation of new periodic orbits. If this is the case, one or more tangency points arise between the reduced phase space section $\C_{\pm}$ and the parabola \eqref{parabola}. Suppose that a tangency point occurs between $\X$ and the  upper arc of the contour $\C$. Let us denote such a point by ${\mathcal Q}_U$. Since  its $X$-coordinate must be positive, if we invert the coordinate transformation \eqref{tr_lem}, it corresponds to \emph{two} points on the section $I_2=0$ of the sphere \eqref{phase_sphere}, the two  \emph{inclined orbits} of \eqref{Ia}-\eqref{Ib}. Similarly, if a tangency point, say $\mathcal Q_L$, occurs on the lower arc of the contour $\C$, it determines the two \emph{loop orbits}  \eqref{La}-\eqref{Lb}. 

In the following subsection we are going to examine all possible cases: to proceed with this general study, we recall that only non-negative value of $\E$ are allowed. For $\E=0$ the dynamics are trivial since the reduced phase space coincides with the origin. Increasing the value of the distinguished parameter the area delimited by the phase-space contour $\C$ increases. Thus, at a certain critical value for $\E$, the contour $\C$ will touch the parabola $\X(Z)$. By varying $h$, $\X$ shifts upward or downward and we have to take into account that the shape and position of the parabola depends not only on $\E$ but on all the parameters of the system. In fact, its concavity is given by $-A/C$ and its vertex is at
\begin{equation}\label{zm}
Z_V=-\frac{B\E+\Delta}{2A},\quad
X_V=\frac1{C}\left(h+\frac{(B\E+\Delta)^2}{4A}\right).
\end{equation}
To study all possible equilibria of the reduced system, the symmetric setting of the problem will help us to simplify the analysis. Restoring the dependence on $Y$ we get the generic intersections between the surfaces $\H$ and $\mathcal L$. These are simple closed curves on the reduced phase space providing all generic trajectories. Their overall structure is uniquely defined by the nature of the critical points we are going to investigate and, for each value of $\E$, by the range of values of $\H$ from the first to the last contact of the two surfaces. 
 
\subsection{Reference case}\label{caso_particolare}

In the case $A<0$, $C>0$ and $\Delta<0$, the parabola \eqref{parabola} is upward concave with a minimum in $Z_V$ which, as shown by eq. \eqref{zm}, does not depend on $h$ and is negative for sufficiently small values of $\E$. The tangency points between $\X$ and $\C$ can be found by imposing that the systems
\begin{equation} \label{su_int}
 \left\{
  \begin{array}{ll}
  X= &\X(Z)\\
  X\in &\C_+
    \end{array}
\right.
\end{equation}
and
\begin{equation}\label{sl_int}
\left\{
  \begin{array}{ll}
  X= &\X(Z)\\
  X\in &\C_-
  \end{array}
\right.
\end{equation}
have a unique solution not coinciding with $\left(\pm\E,0\right)$. This happens if the discriminant of the quadratic equations 
\be\label{ZQZ}
\X(Z)=\pm \left(\E^2-Z^2\right)\ee
vanishes. The left panels in figg. \ref{RC1A} and \ref{RC1B} show possible situations and the right panels display the corresponding surfaces of section.

System \eqref{su_int} admits the unique solution

\be\label{QU} 
{\mathcal Q}_U=\left(Z_U, \E^2-Z_U^2\right), \quad Z_U \doteq \frac{B\E+\Delta}{2(C-A)}\ee
if
\begin{equation}\label{HU}
h=h_U\doteq C\E^2+\frac{(B\E+\Delta)^2}{4(C-A)}=C\E^2+(C-A)Z_U^2.
\end{equation}
This result, which gives the solution for $I_{3U}=Z_U$ mentioned in \eqref{Ia}-\eqref{Ib}, determines a contact point ${\mathcal Q}_U$ on $\C_+$ if it satisfies the constraints
\begin{equation}\label{U_existence}
-\E<Z_U<\E.
\end{equation}
 For $\Delta<0$, these inequalities are verified
\begin{itemize}
   \item if $2(A-C)<B\leq2(C-A)$, for 
   \be \label{en_1u}
   \E>\E_{1U}\doteq \frac{\D}{2(A-C)-B} = \frac{\d}{2(\alpha_2-\alpha_4)-\alpha_3}; \ee
   \item if $B>2(C-A)$, for 
   \be\label{en_2u}
   \E_{1U}<\E<\E_{2U}\doteq \frac{\D}{2(C-A)-B} = \frac{\d}{\alpha_3-2(\alpha_1-\alpha_4)}.\ee
 \end{itemize}
We then obtain the threshold values $\E_{1U}$ and $\E_{2U}$ which corresponds to the two inclined orbits  \eqref{Ia}-\eqref{Ib} bifurcating {\it from} NM1 and annihilating {\it on} NM2. In the left panels in figg. \ref{RC1A} and \ref{RC1B} the contact parabola is displayed in red and provides the red fixed points in the surfaces of section on the right. The ensuing separatrices are displayed in blue in both panels. The nature of the fixed point can be assessed by computing its index \cite{Ku}: the contact point between $\X$ and $\C_+$ has index
\be\label{indU}
{\rm ind} ({\mathcal Q}_U) = {\rm sgn} [C(C-A)].\ee
In the reference case $C>0>A$, therefore ${\rm ind} ({\mathcal Q}_U) > 0$ and the inclined orbits are {\it always stable}. The red fixed points are elliptic in both figg. \ref{RC1A} and \ref{RC1B}.

In a similar fashion, we can proceed with system \eqref{sl_int} and we expect to find now the conditions for the existence of loop orbits. We find the contact point
\be\label{QL} 
\mathcal Q_L=\left(Z_L, \E^2-Z_L^2\right), \quad Z_L \doteq -\frac{B\E+\Delta}{2(A+C)}\ee
if
\begin{equation}\label{HL}
h=h_L\doteq -C\E^2-\frac{(B\E+\Delta)^2}{4(A+C)}=-C\E^2 + (A+C)Z_L^2
\end{equation}
and if it satisfies the constraints
\begin{equation}\label{L_existence}
-\E<Z_L<\E.
\end{equation}
This gives the solution for $I_{3L}=Z_L$ introduced in \eqref{La}-\eqref{Lb} and provides a tangency point $\mathcal Q_L$ on $\C_-$ 
\begin{itemize}
   \item if $2(A+C)<B\leq-2(A+C)$ for 
   \be\label{en_1l} \E>\E_{1L} \doteq \frac{\D}{2(A+C)-B} = \frac{\d}{2(\alpha_2+\alpha_4)-\alpha_3};\ee
   \item if $B>-2(A+C)$ for 
   \be\label{en_2l}
   \E_{1L}<\E<\E_{2L} \doteq \frac{\D}{-2(A+C)-B} = \frac{\d}{\alpha_3-2(\alpha_1+\alpha_4)},\ee
 \end{itemize}
where $\E_{1L}$ and $\E_{2L}$ are the threshold values for the bifurcation of the loop orbits. 

However, we now see that the nature of the contact depends on the ratio $C/A$. Since it is necessary that $A+C\neq0$, in order to proceed we have to distinguish among the three sub-cases:
\begin{enumerate}
  \item $A+C<0 \quad (\e_1=\e_2=-1)$,
  \item $A+C>0 \quad (\e_1=-1,\e_2=1)$,
  \item $A+C=0$ (\rm{degenerate}),
\end{enumerate}
where, for completeness, we mention the correspondence with the coefficients of the germ of the universal deformation \eqref{F_uni} as discussed in the Appendix \cite{VU11}. 

By introducing the {\it curvature} of the parabola $\X$, that is $|C/A|$ and that of the contour \eqref{lemon_arcs}, which is simply equal to one, now we have that the curvature of the parabola is larger than that of the contour in the first sub-case and smaller in the second one. This affects the stability of the bifurcating family associated with the contact with the lower arc. In fact, the contact point between $\X$ and $\C_-$ has index
\be\label{indL1}
{\rm ind} (\mathcal Q_L) = {\rm sgn} [C(A+C)].\ee
In the first sub-case, $A+C<0$, since $C$ is positive, ${\rm ind} (\mathcal Q_L) < 0$ and loop orbits are {\it unstable}: looking at the left panel of fig. \ref{RC1A} we see that the contact with the lower arc is {\it internal} to the lemon and the corresponding fixed points in the surfaces of section are hyperbolic. In the second sub-case, $A+C>0$, ${\rm ind} (\mathcal Q_L) > 0$ and loop orbits are {\it stable}: in the left panel of fig. \ref{RC1B} we see that the contact of the green parabola with the lower arc is {\it external} to the lemon and the corresponding green fixed points in the surfaces of section are elliptic.


\begin{figure}[!h]
\begin{center}
\includegraphics[width=6.5cm]{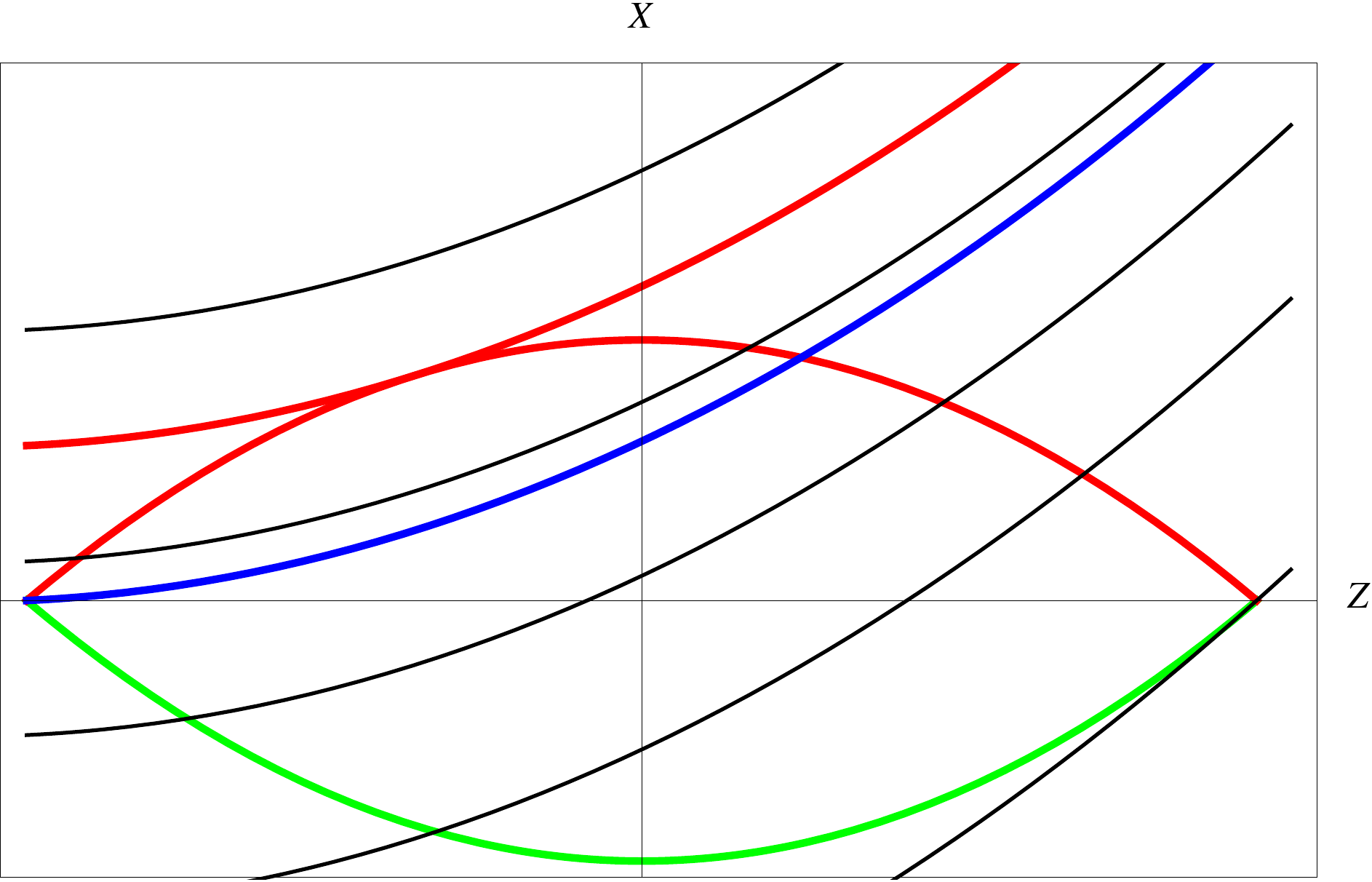}
\includegraphics[width=4.5cm]{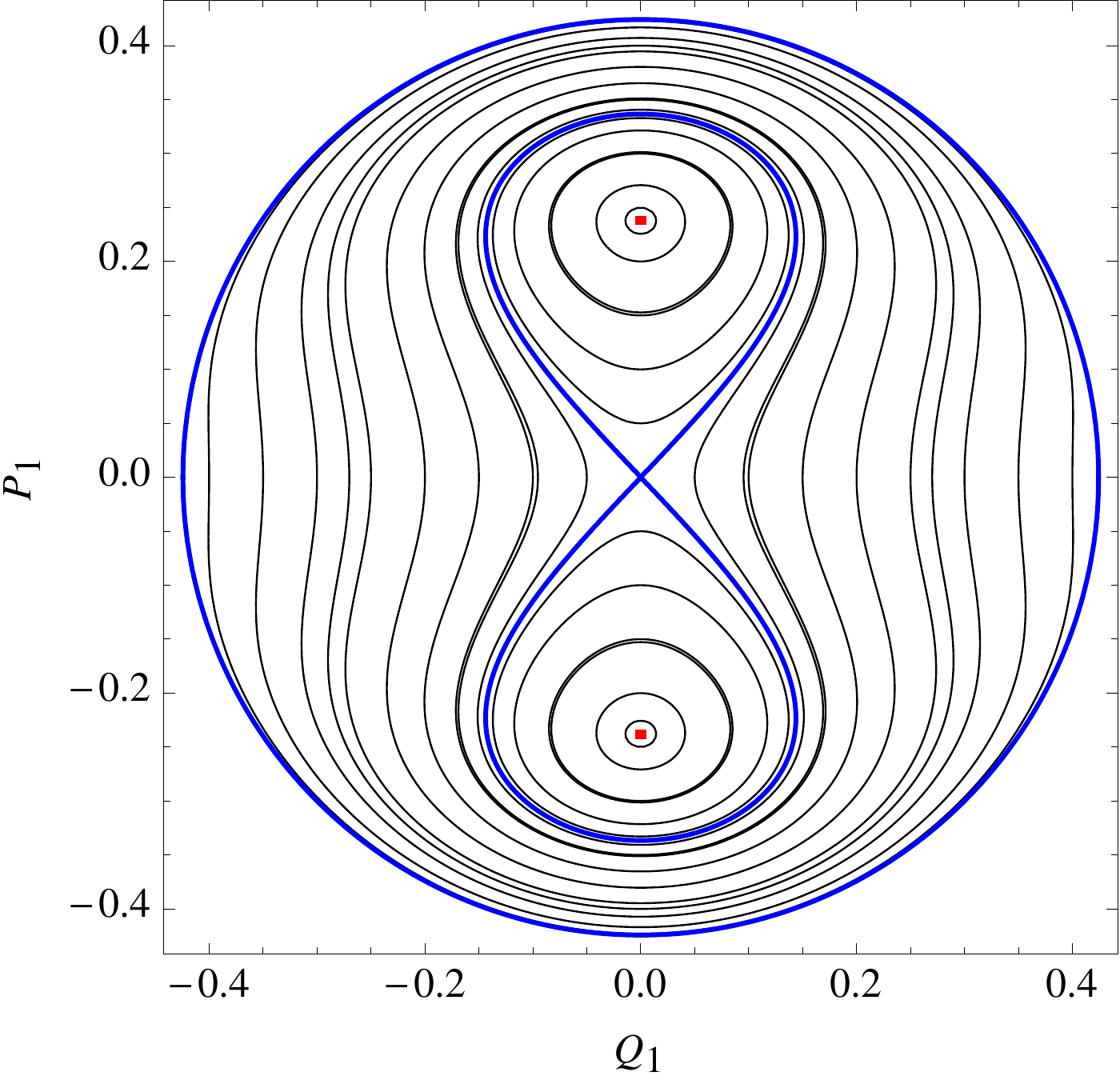}\\
\includegraphics[width=6.5cm]{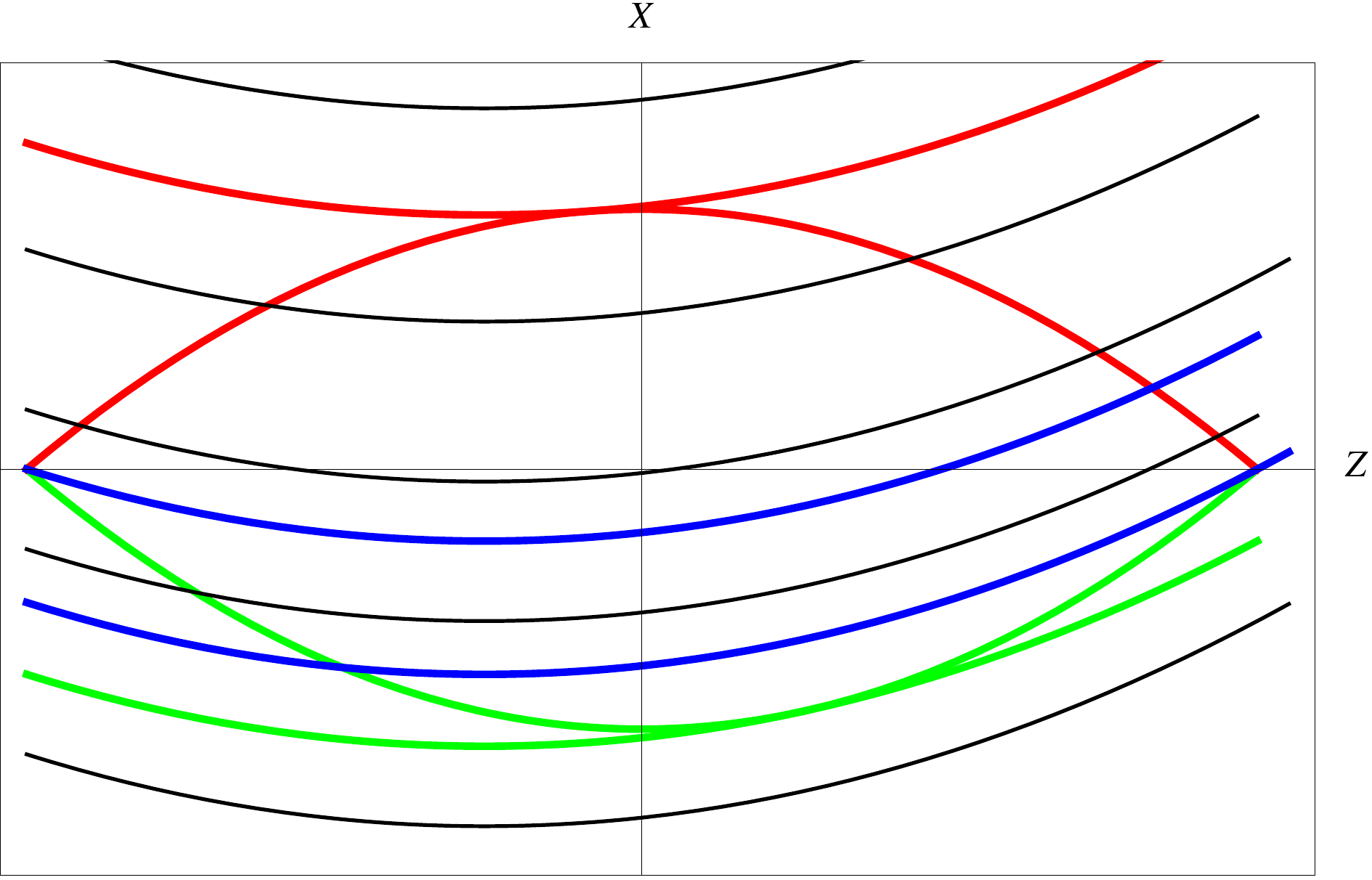}
\includegraphics[width=4.5cm]{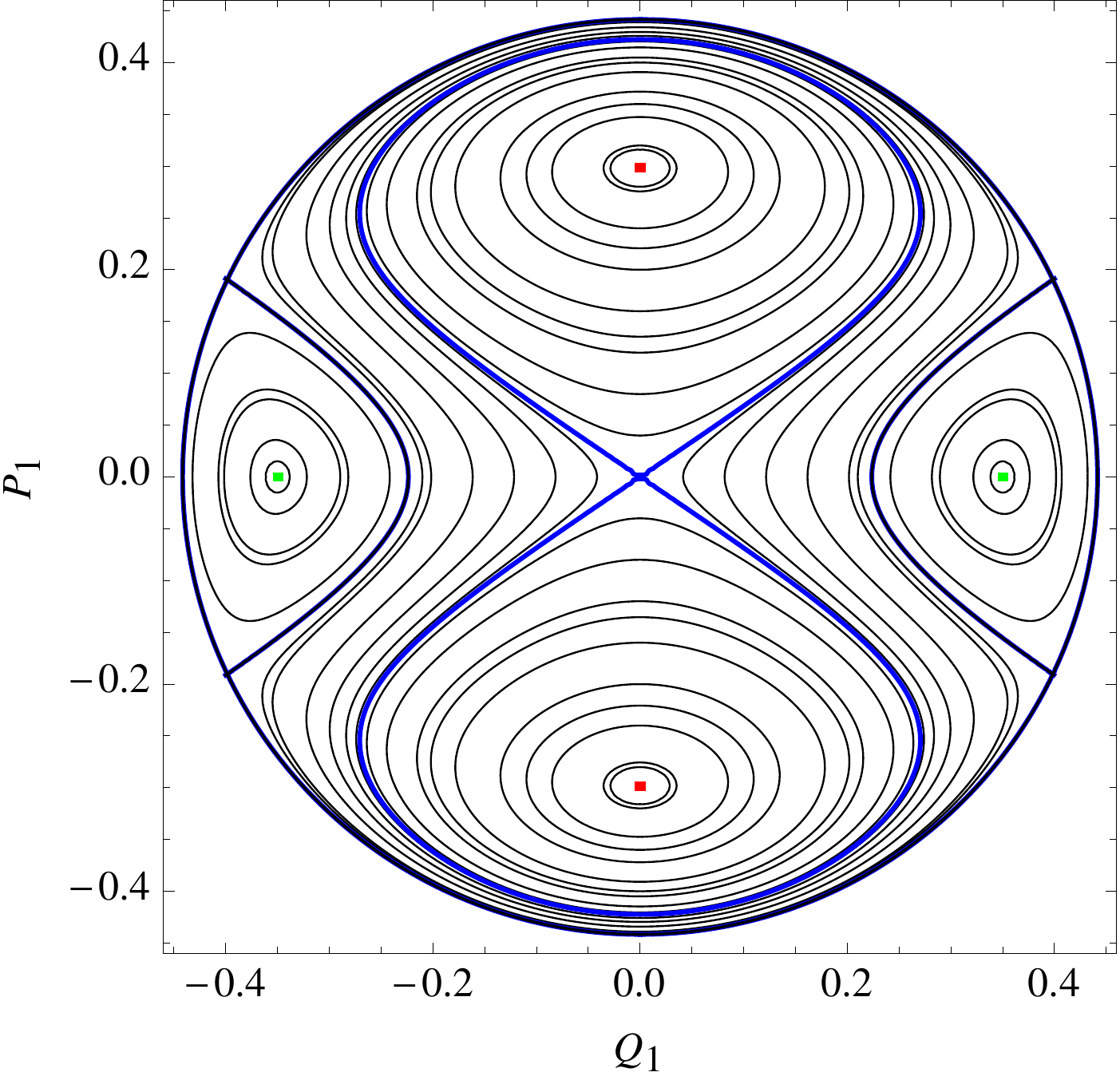}\\
\includegraphics[width=6.5cm]{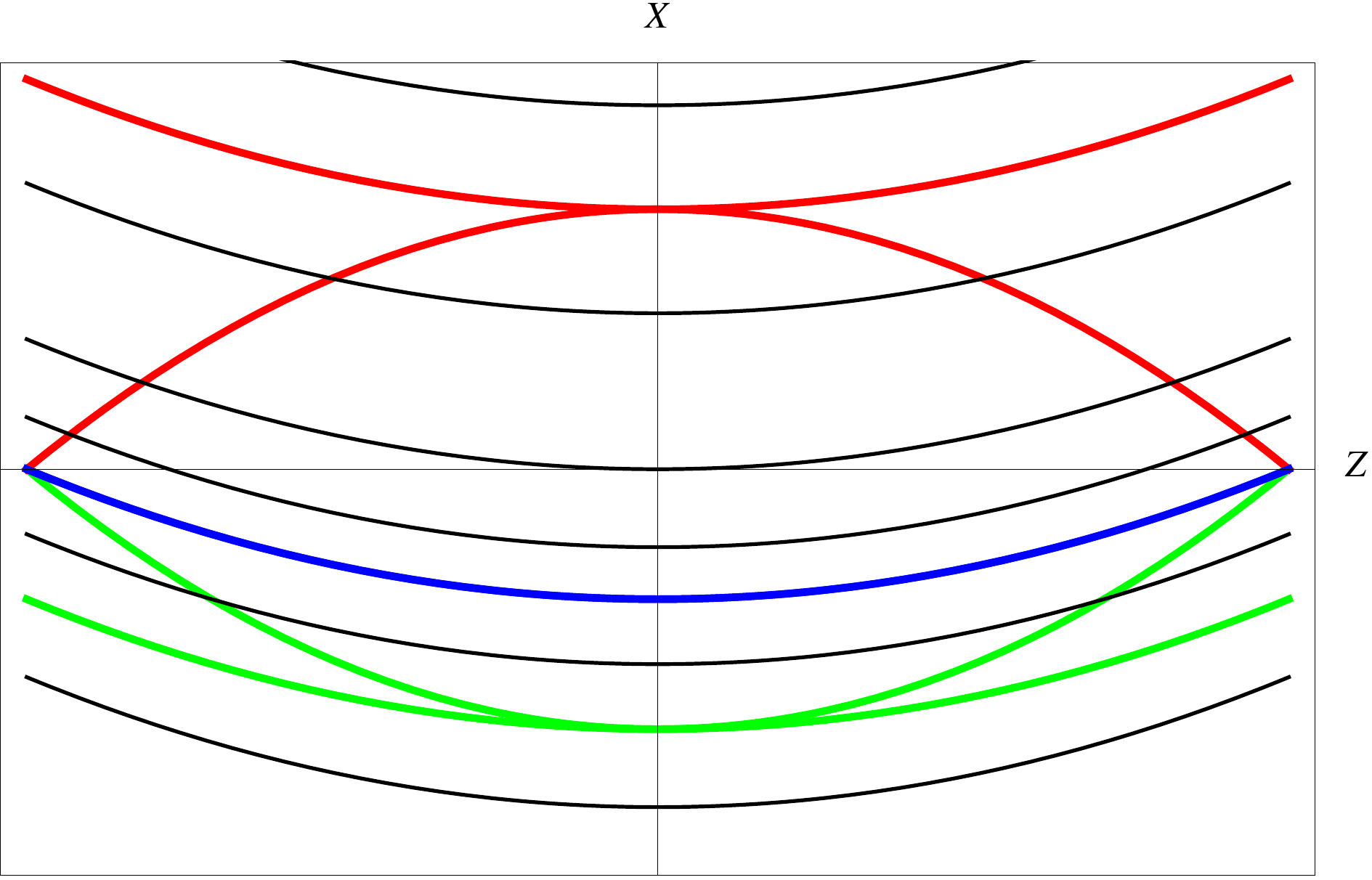}
\includegraphics[width=4.5cm]{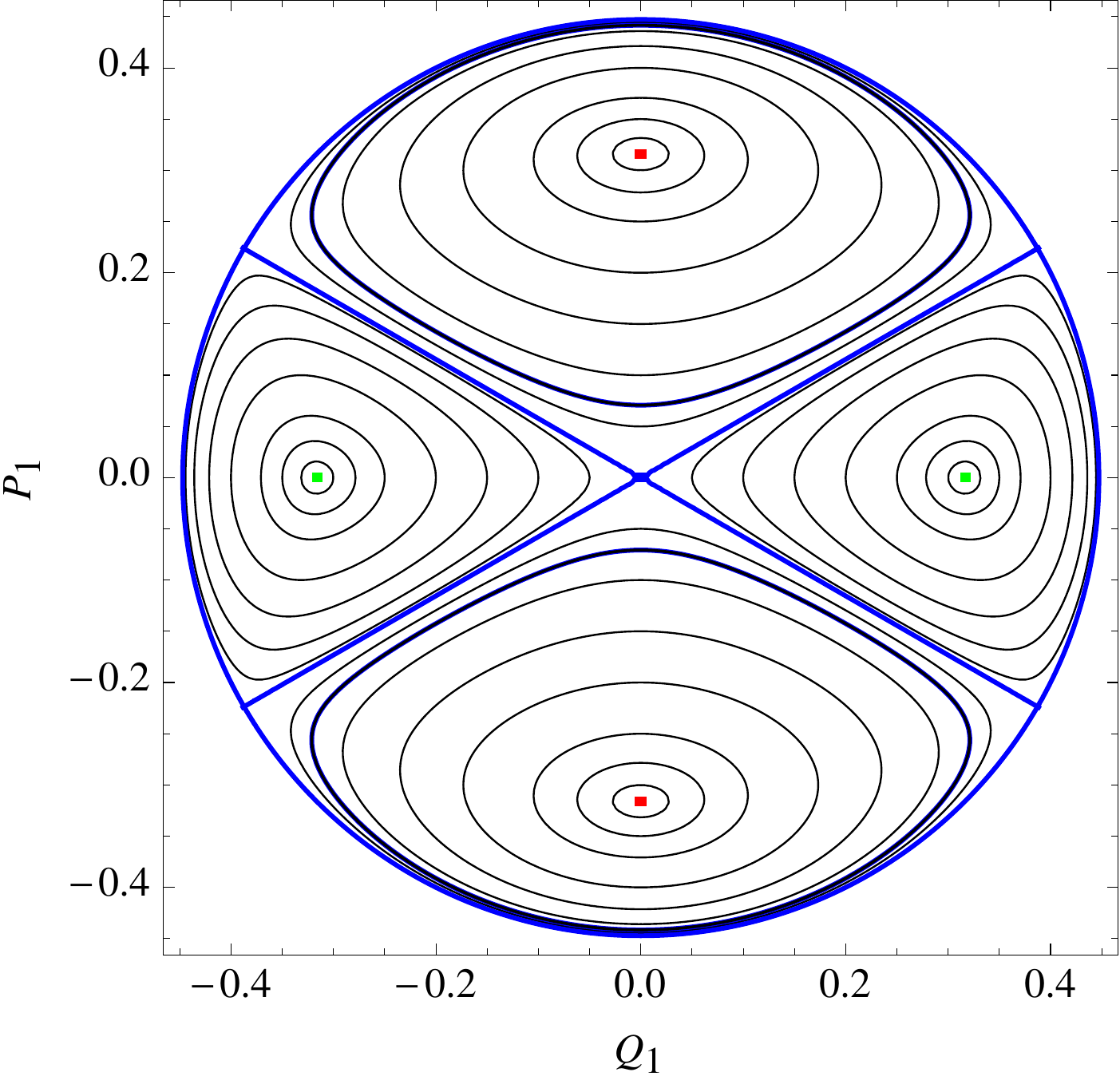}\\
\includegraphics[width=6.5cm]{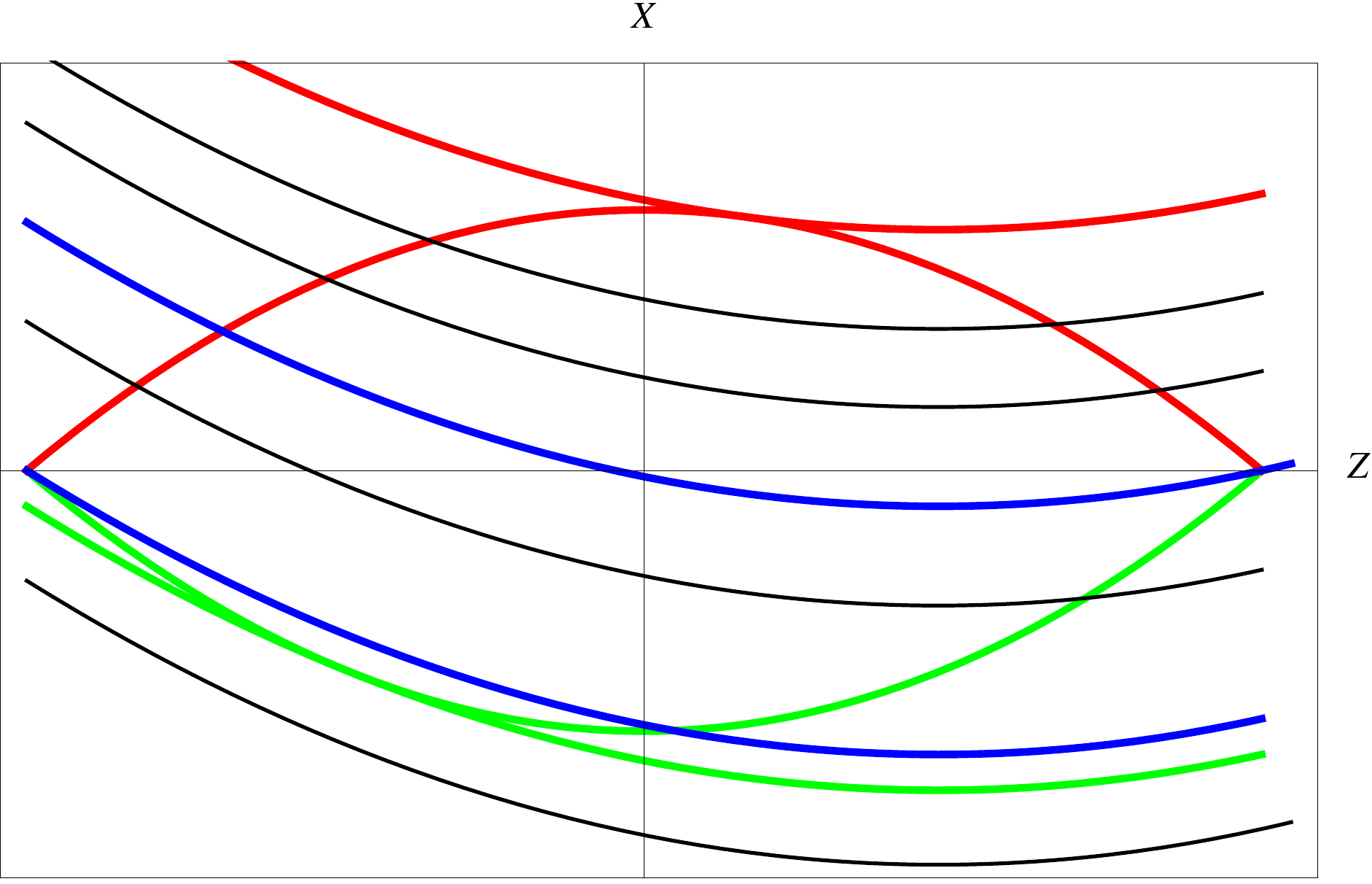}
\includegraphics[width=4.5cm]{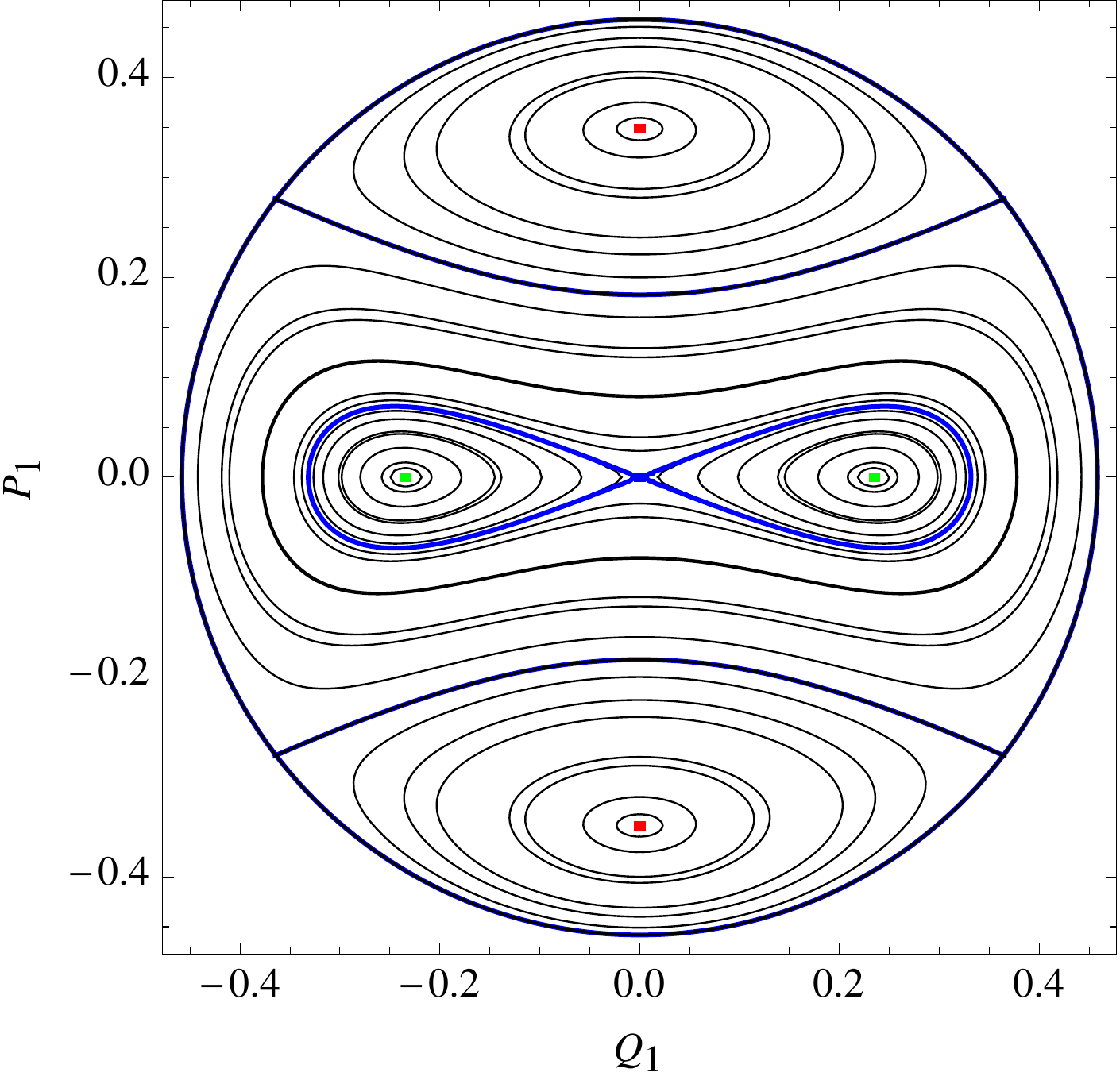}\\
\includegraphics[width=6.5cm]{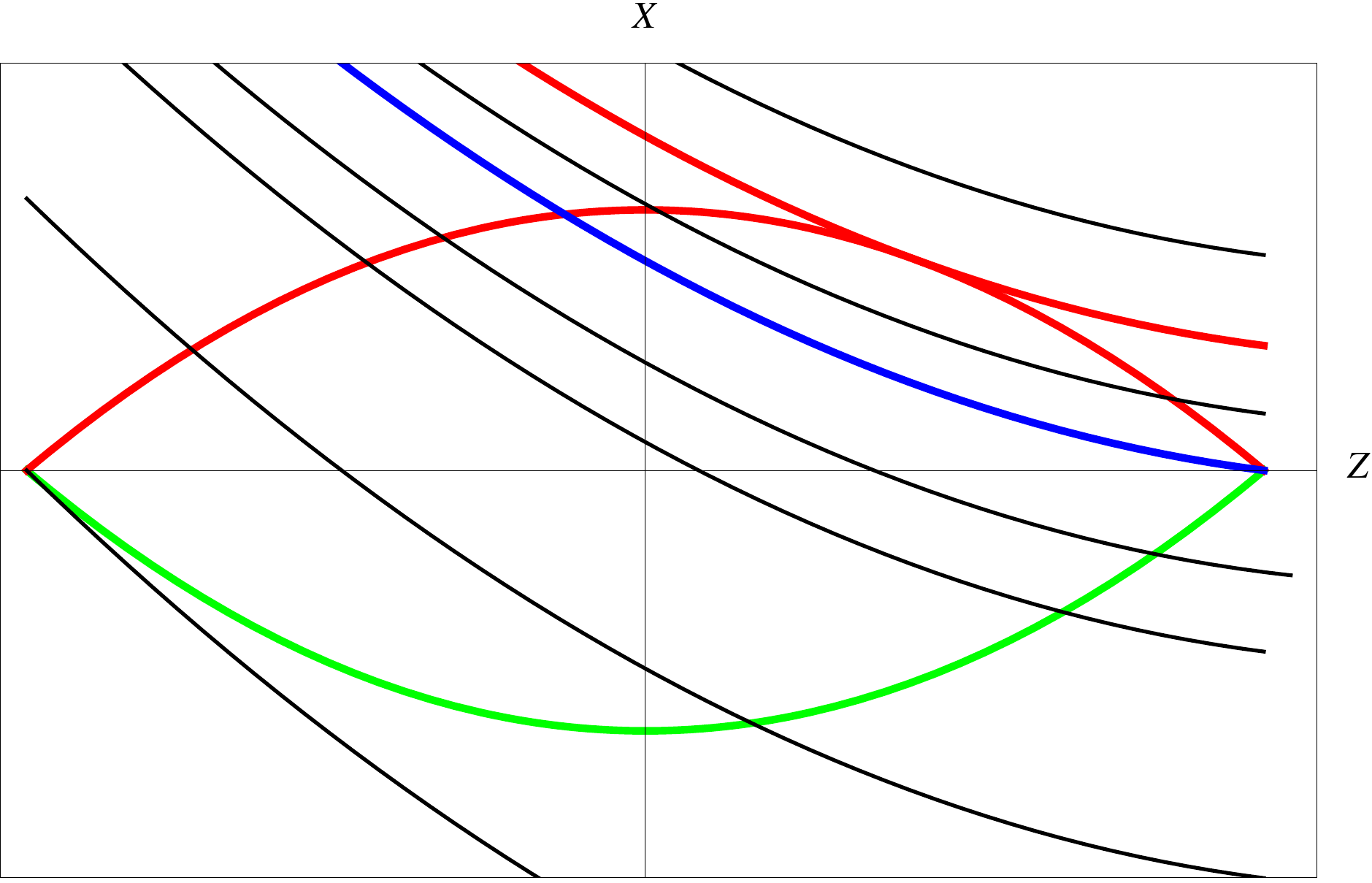}
\includegraphics[width=4.5cm]{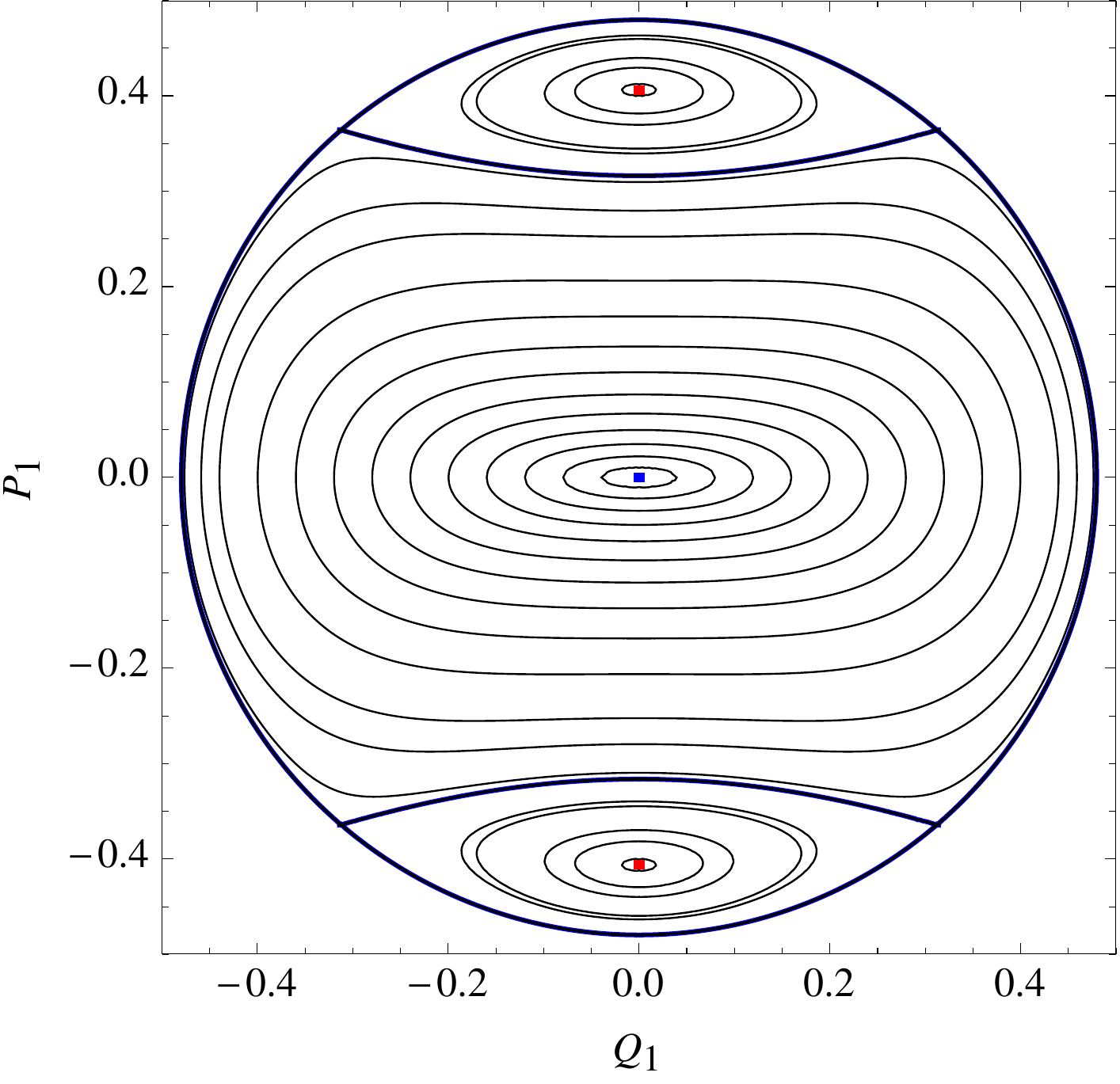}
\end{center}
\caption{\small{Reference case, sub-case 2: $A=-\frac1{10},B=2,C=\frac15,\D=-\frac15$.}}
\label{RC1B} 
\end{figure}


The threshold values of $\E$ are the same as above but, in the second sub-case, the conditions for tangency with the lower arc at $\mathcal Q_L$ are now satisfied 
\begin{itemize}
   \item if $-2(C+A)<B\leq2(C+A)$, for 
   \be\label{en_2L2} \E>\E_{2L};\ee
   \item if $B>2(C+A)$, for 
   \be\label{en_1L2}
   \E_{2L}<\E<\E_{1L},\ee
 \end{itemize}
where $\E_{1L}$ and $\E_{2L}$ are still those defined above in (\ref{en_1l}--\ref{en_2l}). Therefore, now the order of bifurcation is reversed. 

The peculiarity of this sub-case is the `global bifurcation'. Let us consider the critical value of the distinguished parameter 
 \be \E_{GB} \doteq-\frac{\D}{B}.\label{engb}\ee
 Comparing with \eqref{zm}, \eqref{QU}  and \eqref{QL}, we observe that
 \be 
 Z_U (\E_{GB}) = Z_L (\E_{GB}) = Z_V (\E_{GB}) = 0
 \ee
and we have a family of parabolas with axis coinciding with the $X$-axis. From \eqref{yax_energy}--\eqref{xax_energy}, at the value of the Hamiltonian
\be
h_1=h_2=\frac{A \ \D^2}{B^2}\doteq h_{GB},\ee
the parabola passes through \emph{both} points $\mathcal Q_1$ and $\mathcal Q_2$ and a simple computation shows that its minimum is negative but bigger then $\E^2$. The geometrical setting is that of the central panels of fig.\ref{RC1B}. 


In the degenerate case $A=-C$, $\X$ and the lower arc of $\C$ have the same curvature. Hence, by a simple geometrical argument we see that if
$Z_V\neq0$, it is impossible to have isolated intersections point different from $\mathcal Q_1$ between $\X$ and $\C_-$: \emph{all} points of the lower arc of $\C$ are tangency points between $\X$ and $\C$. Thus if $B>0$ $(B<0)$ and $\D<0$ $(\D>0)$, for $ \E=\E_{GB} $
we find infinite (non-isolated) equilibria given by all the points on $\C_-$. They correspond to the circle $I_1=0$ on the sphere \eqref{phase_sphere}. Only inclined orbits may bifurcate as isolated periodic orbits and this happens when a contact between $\X$ and $\C_+$ does occur. 

\begin{oss}
The case $\D>0$ can be treated as $\D<0$ by a transformation which exchange the coordinate axes in the original phase space. On the reduced phase space  it corresponds to the reflection $R_1$. As a consequence, the equilibrium points $\mathcal Q_1$ and $\mathcal Q_2$ are exchanged and the parabola $\X$ is reflected into its symmetric with respect to the $X$-axis. Thus the situation is `specular' with respect to the previous one. The new equilibria ${\mathcal Q}_U$ and $\mathcal Q_L$ bifurcate now from the fixed point $\mathcal Q_2$, instead of $\mathcal Q_1$ and the bifurcation sequences are reversed.
\end{oss}

\subsection{Complementary cases}\label{compl}

In the previous section we considered the `reference' case $A<0$ with $C>0$. Now we are going to study the dynamics of the system in the `complementary' cases:

\begin{description}
  \item[a)] $A<0$, $C<0$;
  \item[b)] $A>0$, $C<0$;
  \item[c)] $A>0$, $C>0$.
\end{description}
As observed above, by applying the transformations \eqref{R_z}, $\eqref{R_x}$ and their compositions, the orbital structure of the system in these cases can be deduced from the analysis of subsection \ref{caso_particolare}.

In case {\bf a)}, the critical point $Z_V$ does not change its sign, but the parabola $\X$ turns out to be downward concave.
However we can reverse its concavity by a simple application of $R_2$. This allows us to deduce the dynamics of the system from the analysis given above. We have to understand how $R_2$ operates on the fixed points of the reduced system. Since $R_2$ is a symmetry with respect to the $Z$-axis, the two degenerate equilibria are invariant under its action. On the other hand, if a tangency point occurs on of $\C_+$ it is reflected into  a tangency point on $C_-$ and vice-versa. This implies that the role of loop and inclined orbits is exchanged (cfr. the right panel in table \ref{T1}). Namely, the first periodic orbits to appear from NM1 are now the loop orbits. The corresponding threshold value for the distinguished parameter is again $\E=\E_{1L}$ given in \eqref{en_1l}.
The bifurcation of inclined orbits is possible from NM1 in the case $A<C$ for $\E>\E_{1U}$ and from NM2 for $\E>\E_{2U}$ in the case $C<A<0$. In fig.\ref{C1} we can see the relation between the Reference Case and the complementary case {\bf a)}: in the upper panel we see the sections after the first and the second bifurcation in the case $A<0, C>0$, in the lower panel the corresponding bifurcations in the case $A<0, C<0$.

\begin{figure}[htbp]
\centering%
{
\includegraphics[width=6cm]{rc1S2.pdf}
\includegraphics[width=6cm]{rc1S3.pdf}
\includegraphics[width=6cm]{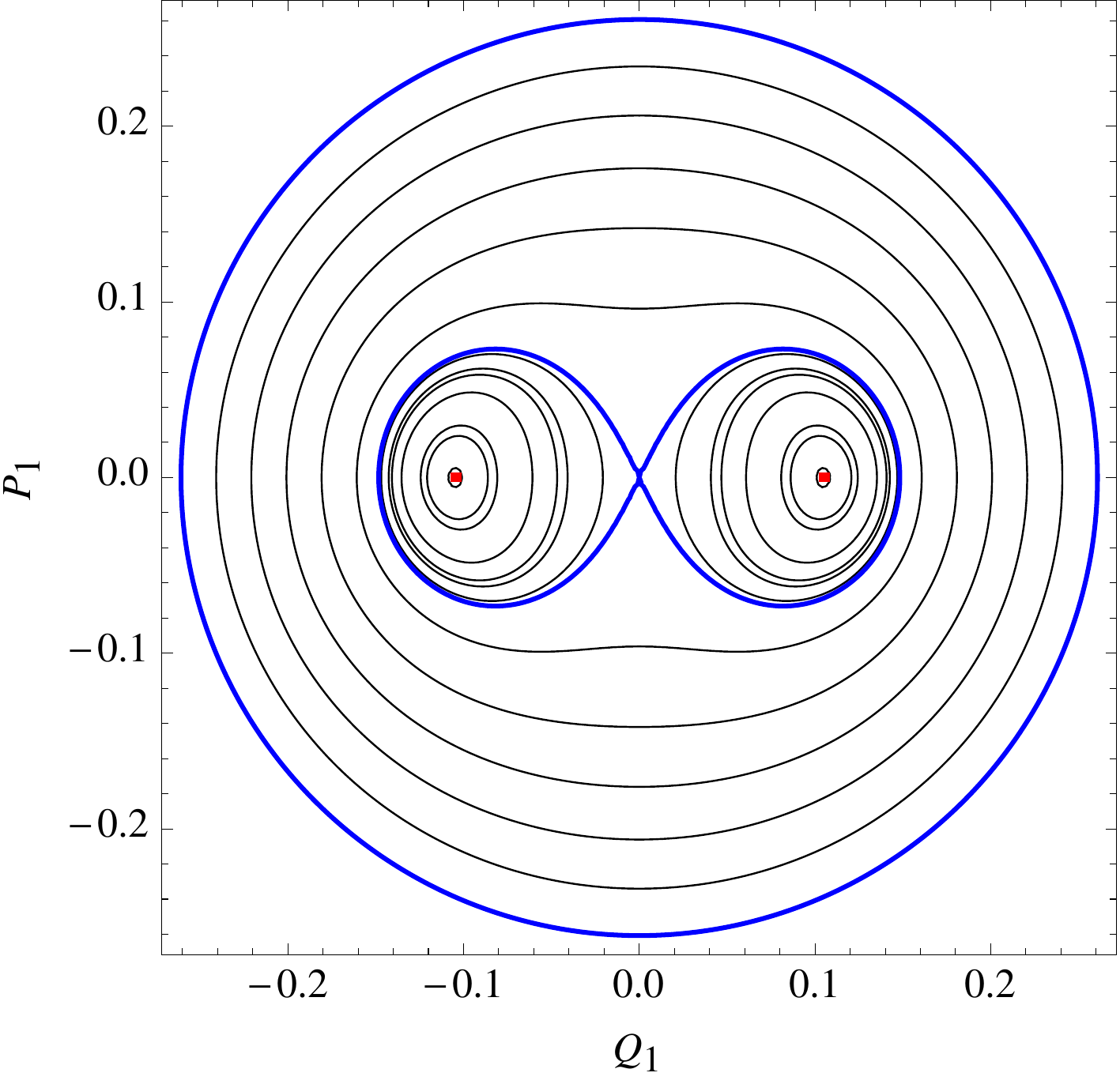}
\includegraphics[width=6cm]{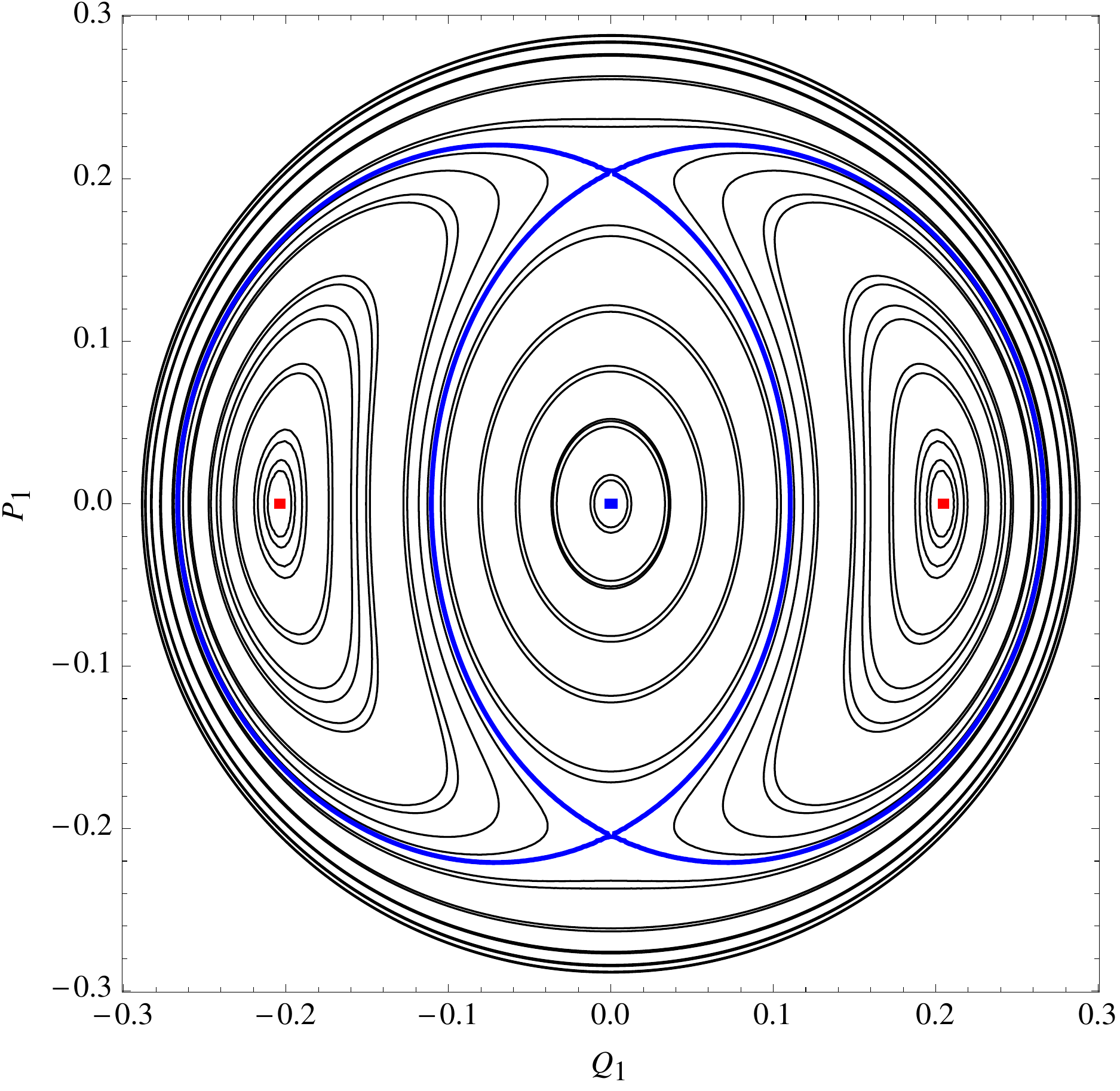}}
\caption{\small{Upper panel, Reference Case  ($A<0, C>0$). Lower panel, Complementary Case  ($A<0, C<0$)}}
\label{C1}
\end{figure}

\begin{figure}[htbp]
\centering%
{
\includegraphics[width=6cm]{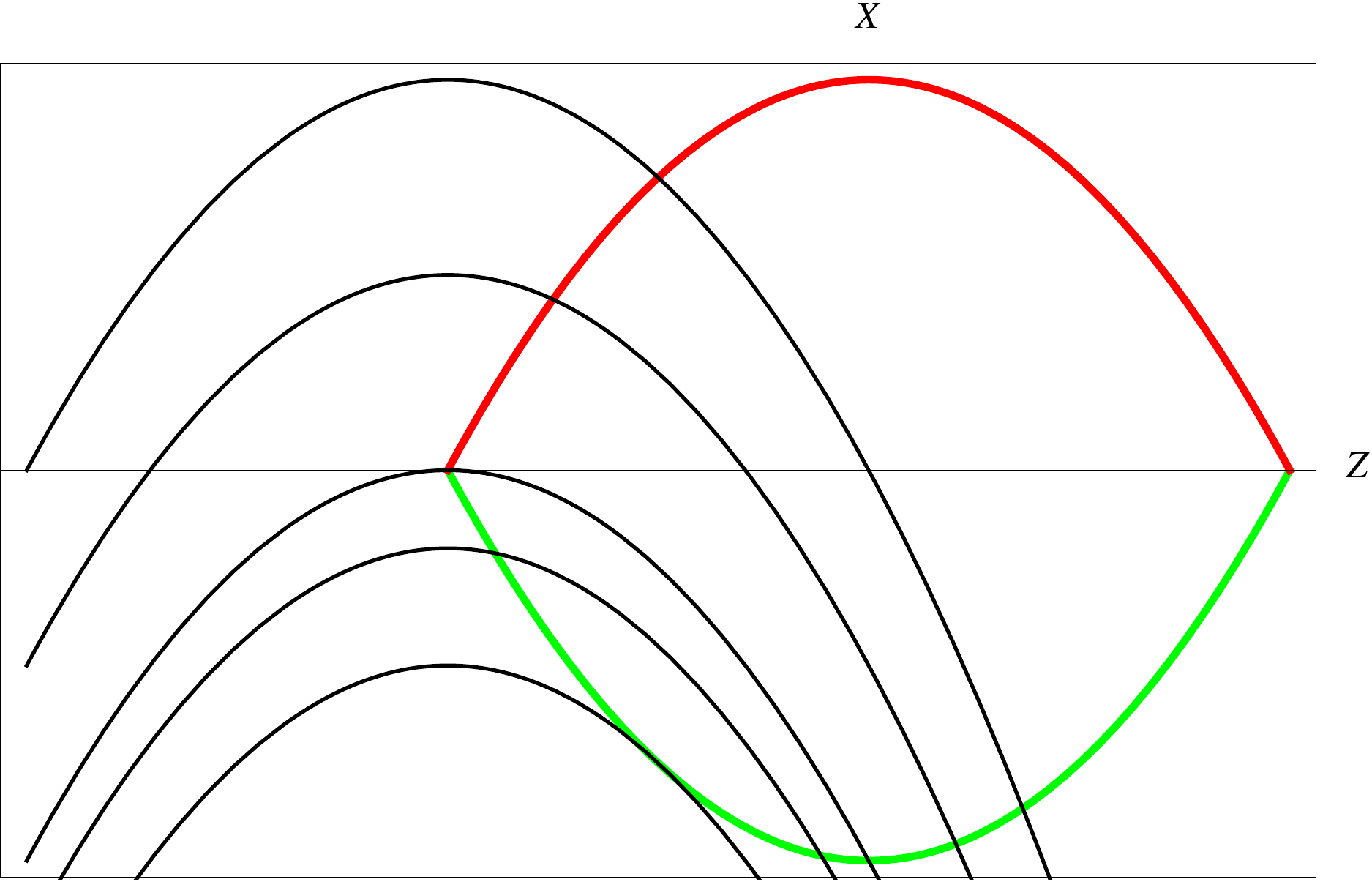}}
\caption{\small{The particular case $A=-\frac1{16},B=0,C=-\frac1{16},\Delta=-\frac14$, $\E>\E_{1L}=1.$}}
\label{C1D1}
\centering%
{
\includegraphics[width=5cm]{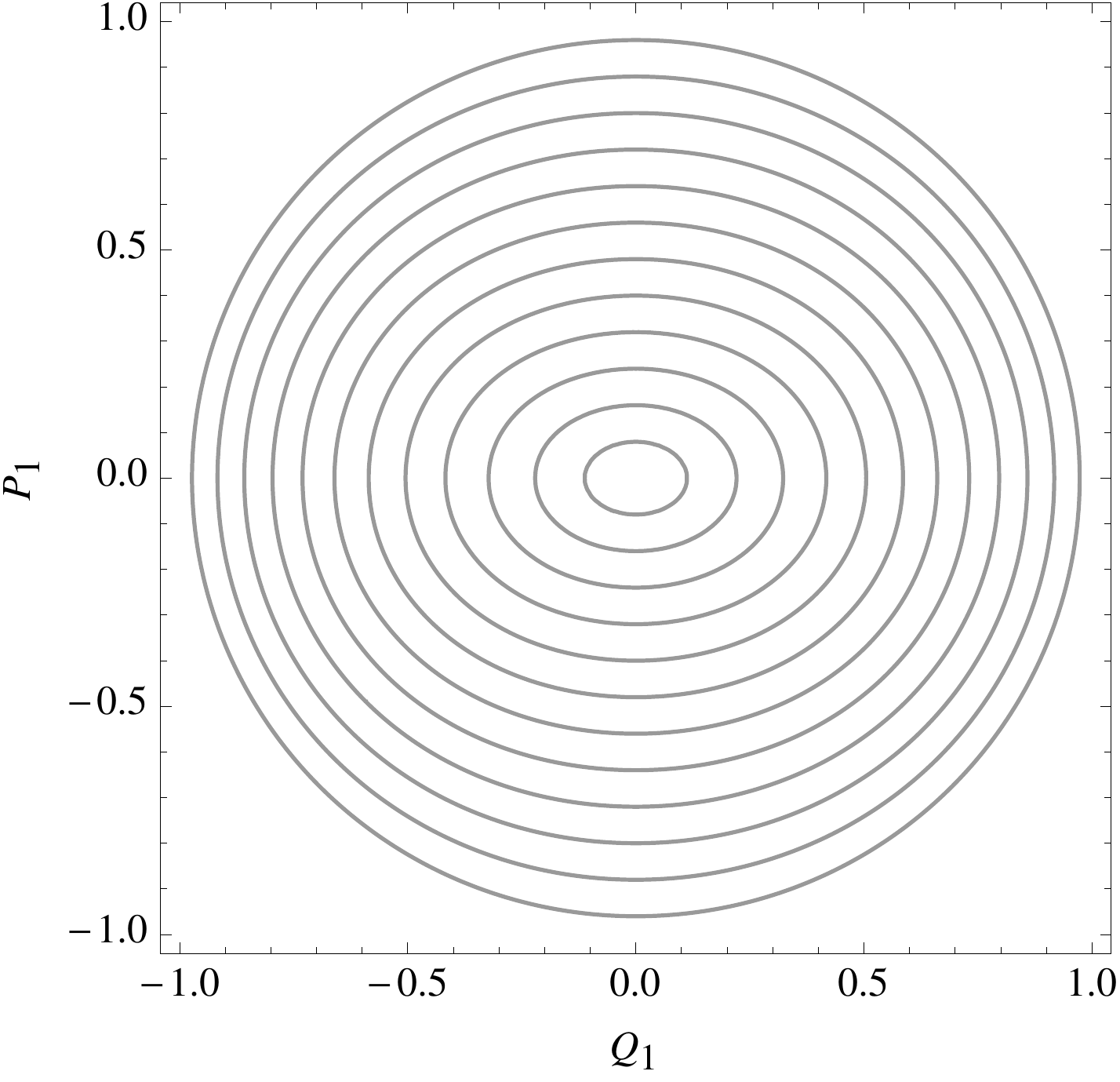}
\includegraphics[width=5cm]{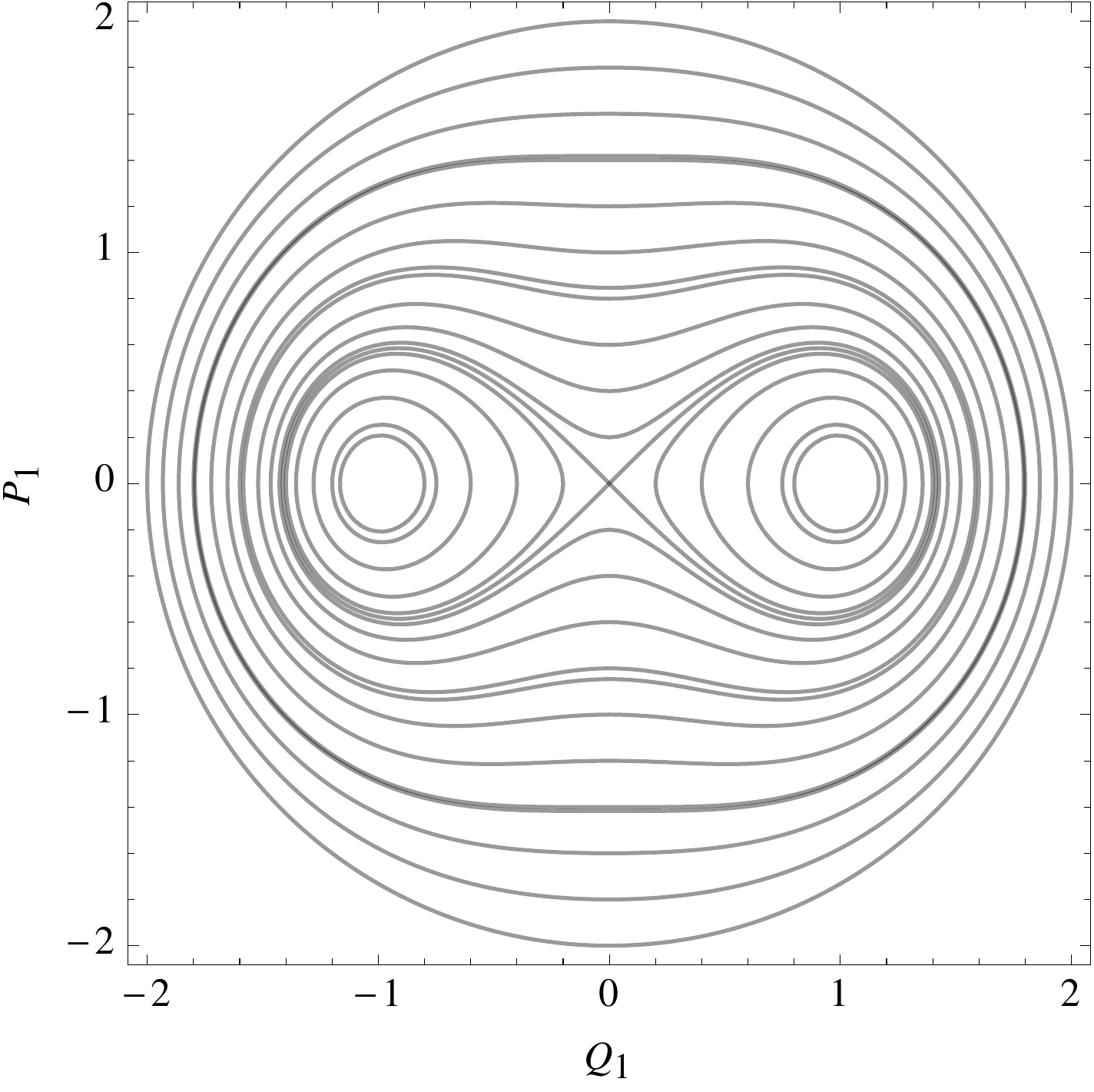}}
\caption{\small{Two PSS corresponding to the case of fig.\ref{C1D1}: $\E=0.5, \; \E=2.$}}
\label{C1D2}
\end{figure}

The degenerate case $A=C$ is specular with respect to the case $A=-C$ with $C>0$. It admits as an interesting example the family of natural systems with elliptical equipotentials \cite{MP13a} often used in galactic dynamics. Inclined are forbidden and only loop orbits may bifurcate as isolated periodic orbits when a contact between $\X$ and $\C_-$ occurs: see figg.\ref{C1D1}--\ref{C1D2}. Even more degenerate is the case in which, in addition to $A=C$ we also assume $\D=B=0$: to this subset belongs the famous H\'enon-Helies system to which we come back in subsection \ref{degenz}.

In case {\bf b)}, $\X$ is upward concave and its maximum point lies on the positive $Z$-axis. Thus, by applying $R_1$ we can deduce the orbital structure of the system from the case \ref{caso_particolare}. Under the action of $R_1$ the degenerate equilibria of the reduced system are exchanged. Furthermore, each tangency point between $\X$ and $\C$ is reflected into its symmetric with respect to the $X$-axis (cfr. the right panel in table \ref{T1}). Namely,
$$\mathcal Q_L\equiv(Z_L,X_L)\rightarrow\widetilde {\mathcal Q}_L\equiv(-Z_L,X_L),$$
$${\mathcal Q}_U\equiv(Z_U,X_U)\rightarrow\widetilde {\mathcal Q}_U\equiv(-Z_U,X_U).$$
Anyway, due to the singularity of the transformation \eqref{tr_lem}, to the points $\mathcal Q_L$ and $\widetilde {\mathcal Q}_L$ correspond the same two points on the section $I_1=0$ of the sphere \eqref{phase_sphere}, that is the same loop orbits for the two degree of freedom system. Thus loop orbits are invariant under the action of $R_1$. By a similar argument it follows the invariance of inclined orbits. However, since the degenerate equilibria on the reduced phase space are exchanged, if in the case \ref{caso_particolare}
a periodic orbit bifurcates form NM2, in case {\bf b)} it bifurcates from NM1 and vice-versa.

Finally, by applying $R_2 \circ R_1$ we obtain the stability analysis in case {\bf c)} from the case  \ref{caso_particolare}. The fixed points of the reduced system change according to the right panel of table \ref{T1}. As a consequence, the normal modes exchange their roles and the bifurcation order of inclined and loop orbits is reversed. 

For sake of completeness, we remind the correspondence of the three complementary cases listed above with the classification of the germ of the universal deformation \eqref{F_uni} according to the following scheme: 

\begin{description}
  \item[a)] $A<C: \; \e_1=\e_2=-1; \quad C<A<0: \; \e_1=-\e_2=1$;
  \item[b)] $A>C: \; \e_1=\e_2=1; \quad -C>A>0: \; \e_1=-\e_2=1$;
  \item[c)] $A>C: \; \e_1=\e_2=1; \quad C>A>0: \; \e_1=-\e_2=-1$.
\end{description}

\subsection{Degenerate cases}\label{dege}

We recall that the  degenerate cases of the system correspond to the parameters values $A=\pm C$, $C=0$ and $A=0$. Let us examine them separately.

\subsubsection{$A=\pm C$}
These cases have already been described in subsections \ref{caso_particolare} and \ref{compl}. Here we just want to remark the insights provided by their investigations. They represent those non-generic (but non-trivial) cases in which the present analysis is unable to give reliable results and therefore calls for a higher-order analysis. Subsection \ref{cata} offers another view of these occurrences by introducing the relation between degeneracy and structural instability.

\subsubsection{ $C=0$}
For $C=0$ the parabola $\X$ degenerates into a pair of straight lines both parallel to the vertical $X$-axis. Thus, for all positive values of $\E$, the reduced system has only two equilibria represented by the singular points $\mathcal Q_1$ and $\mathcal Q_2$: the only periodic orbits allowed by the two degree of freedom Hamiltonian are the normal modes.
This is not surprising since this case corresponds to two uncoupled non-linear oscillators.

\subsubsection{$A=0$}
In the case $A=0$ and $C>0$, the parabola $\X$ degenerates into the straight line

\begin{equation}
X =\frac{h}{C}-\frac{B\E+\Delta}{C}Z.
\end{equation}
Let us denote it  by ${\mathcal Y}(Z)$. Its angular coefficient is given by
\begin{equation}
m\doteq-\frac{B\E+\Delta}{C}.
\end{equation}
For $\D<0$, $m$ is positive if and only if $B\leq0$ or $B>0$ and $\E<\E_{GB}$. Thus, for $\E<\E_{GB}$, if  $\mathcal Y$ passes through the point $\mathcal Q_1$, it may intersect the contour phase space $\C$ only in one further point  on its upper arc. The corresponding value for $h$ is given  by
\begin{equation}
h=\overline h:=-\frac{(B\E+\Delta)\E}{2C}.
\end{equation}
If this is the case, the fixed point $\mathcal Q_1$ results to be an unstable equilibrium. A similar argument shows that, if $m<0$ and $h=\bar h$, $\mathcal Y$ may intersects $\C$ only in one further point on its lower arc. Thus the critical value $\E=\E_{GB}$
does not determine a transition to stability/instability for the fixed point $\mathcal Q_1$. As in the case $-C<A<0$, it corresponds to a global bifurcation for the system. In fact, for $m=0$, the straight line $\mathcal Y$ becomes parallel to the $Z$-axis and for $h=0$ it passes through both the degenerate fixed points. Hence, for $\E=\E_{GB}$ they turn out to be both unstable and their stable and unstable manifolds coincide.

Thus, the stability analysis of the normal mode NM1 for $\D<0$ gives
\begin{itemize}
  \item if $B\leq -2C$ stability for all positive values of $\E$;
  \item if $-2C<B\leq 2C$ instability for $\E>\E_{1U}$;
  \item if $B> 2C$ instability for $\E_{1U}<\E<\E_{2U}$;
\end{itemize}
where the critical threshold are now given by
\begin{eqnarray}
\E_{1U}&=&-\frac{\D}{B+2C} \\
\E_{2U}&=&-\frac{\D}{B-2C}.
\end{eqnarray}
Now, by the symmetry of the reduced phase space,  if  $\mathcal Y$ intersect $\C$ on its upper arc for $h=\overline h$, then by decreasing $h$ enough it will intersect the contour phase space at $\mathcal Q_2$ and on one further point on $\C_-$. Thus, the fixed point $\mathcal Q_2$ turns out to be unstable exactly when the equilibrium point $\mathcal Q_1$  is also unstable! Indeed an easy computation shows that
$$\E_{1U}=\E_{2L}, \;\;\; \E_{1L}=\E_{2U}.$$
 Moreover, by the same argument used above, we see that a tangency point may occur on the upper arc of $\C$ if and only if a tangency point arises between $\mathcal Y$ and $\C_-$.
Hence   the fixed points $\mathcal Q_U$ and $\mathcal Q_L$ (and, as a consequence, loop and inclined orbits) bifurcate at the same time:
\begin{itemize}
  \item if $-2C<B\leq2C$ for $\E>\E_{1U}$;
  \item if $B>2C$ for $\E_{1U}<\E<\E_{2U}$.
\end{itemize}

\begin{figure}
\begin{center}
 \includegraphics[width=12cm]{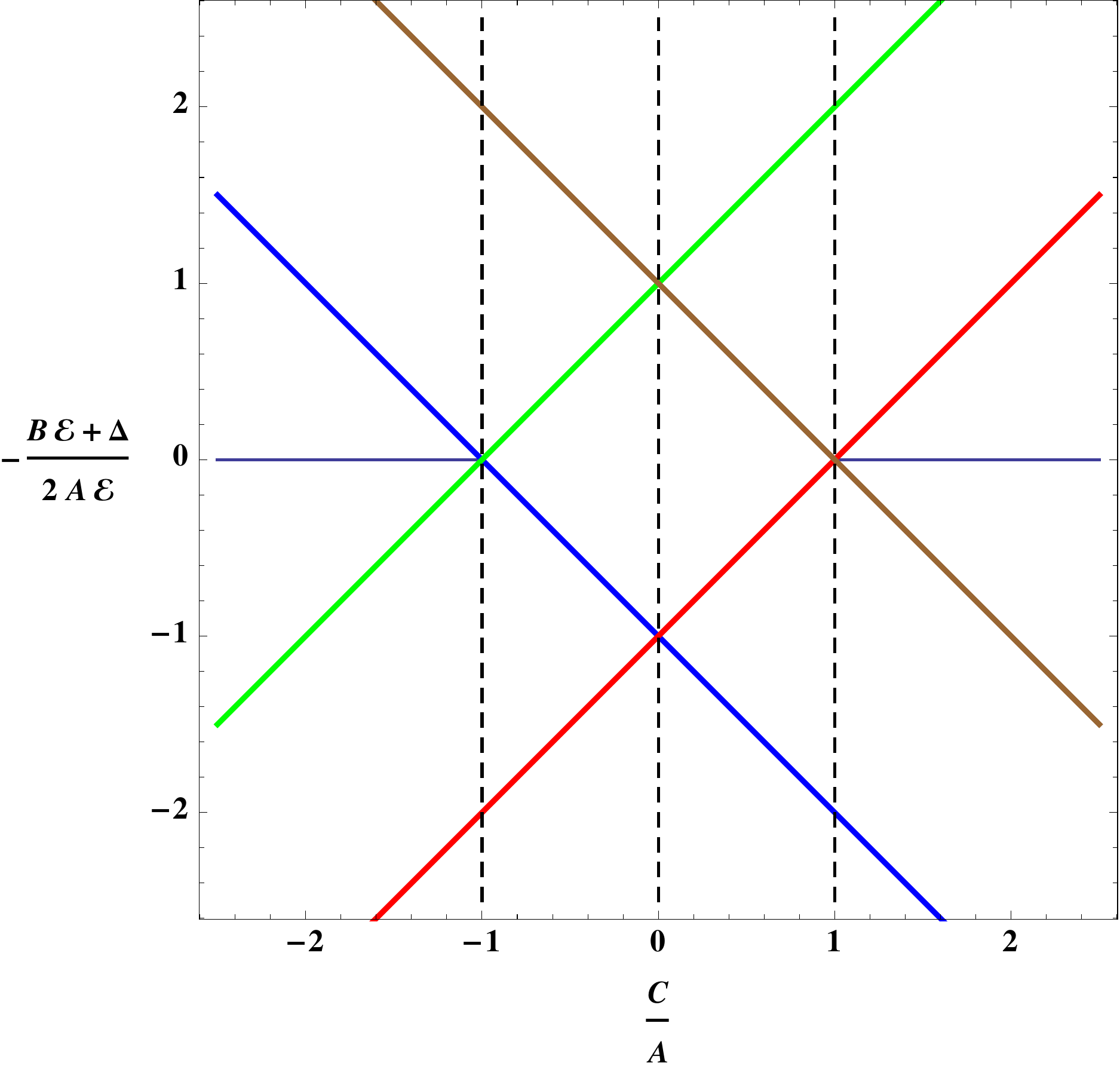}
 \end{center}
\caption{Catastrophe map: the bifurcation lines are associated with $\E_{1U}$ (eq.\eqref{zu1}, red line), $\E_{2U}$ (eq.\eqref{zu2}, brown line), $\E_{1L}$ (eq.\eqref{zl1}, blue line), $\E_{2L}$ (eq.\eqref{zl2}, green line). The vertical dashed lines denote the degenerate cases of subsection \ref{dege}.}
\label{cplot}       
\end{figure}

Since loop and inclined orbits bifurcate at the same time, in the case $C<0$ the orbital structure of the system does not change, even if $\X$ reverses
its concavity.

\subsection{Bifurcation sequences}\label{cata}

A comprehensive way to illustrate the general results described above is obtained by introducing  a pair of  combinations of the internal and external parameters and plot the bifurcation relations on the plane of this pair. This has been derived for the present class of systems in Pucacco \& Marchesiello [2014] and is also referred to as the `catastrophe map' in the physical-chemical literature \cite{STK}. Recalling the determining role of the curvature of the parabola and the four cases generated by the signs of $A$ and $C$ (the `reference' and the complementary cases), we can use 
$C/A$ 
as `coupling' parameter. On the other hand, a parameter which usefully combines the internal parameters $\E,\D$ with the remaining control parameter $B$ is the `asymmetry' parameter
\begin{equation}\label{zmp}
\frac{Z_V(\E)}{\E}=-\frac{B\E+\Delta}{2A\E}.
\end{equation}

On the bifurcation lines we get
\ba
\frac{Z_V(\E_{1U})}{\E_{1U}}&=&\frac{C}{A}-1,\label{zu1}\\
\frac{Z_V(\E_{2U})}{\E_{2U}}&=&1-\frac{C}{A},\label{zu2}\\
\frac{Z_V(\E_{1L})}{\E_{1L}}&=&-1-\frac{C}{A},\label{zl1}\\
\frac{Z_V(\E_{2L})}{\E_{2L}}&=&1+\frac{C}{A},\label{zl2}\ea
whereas the line
\begin{equation}\label{zmgb}
\frac{Z_V(\E_{GB})}{\E_{GB}}=0
\end{equation}
is associated with the global bifurcation. Plotting these lines on the plane of the coupling and asymmetry parameters (see fig.\ref{cplot}),  produces regions with no, one or two  families of periodic orbits in general position. The two triangular regions with bases on the lower/upper sides of the plot are below/above any bifurcation line, therefore they admit only normal modes. The central square is the locus with two bifurcations and therefore admits two families (one stable, the other unstable). The two triangular regions with bases on the lateral sides of the plot have two stable families: the horizontal segments are the loci of global bifurcation. The remaining regions have only one stable family of either type. 

\begin{table}
  \tbl{List of bifurcation sequences as determined by the possible combinations of detuning and control parameters.}
  {\begin{tabular}{@{}ccccccccc@{}}
  
  \hline
  \hline
     
$C/A$   & $\Delta \cdot A > 0 $  & & \multicolumn{4}{c}{Bifurcation sequences} & & $ \Delta \cdot A < 0 $ \\
\hline
  \hline
 $C/A<0$  & $A + C < 0 $ & $\rightarrow$ & 1U    &  1$\ell $ & 2$\ell$ & 2U & $\leftarrow$ & $A + C > 0 $ \\
 $C/A<0$  & $A + C > 0 $ & $\rightarrow$ & 1U    &  2L & 1L & 2U & $\leftarrow$ & $A + C < 0 $ \\
\hline
$C/A>0$  & $A - C < 0 $ & $\rightarrow$ & 1L    &  1$u$ & 2$u$ & 2L& $\leftarrow$ & $A - C > 0 $ \\
 $C/A>0$  & $A - C > 0 $ & $\rightarrow$ & 1U    &  2U & 1U & 2L & $\leftarrow$ & $A - C < 0 $ \\\hline
\end{tabular}}

\label{TT}
\end{table}

The bifurcation sequences as functions of the distinguished parameter $\E$ are summarized in Table 2, where, in the first two lines we have the reference case (the arrow points out that the corresponding sequences are covered from left to right) and the first complementary case (with the arrow pointing out that the sequences are covered from right to left). In the third and fourth lines are represented the remaining complementary cases with analogous rules for the sense along the sequence. Each bifurcation is labeled by the figure 1 or 2 referring to the normal mode from/to which it pertains and by the letters U/$u$ for the (stable/unstable) inclined families and L/$\ell$ for the (stable/unstable) loop families. Overall, sequences of bifurcation of always stable families appear when the condition
\be\label{stabfam}
\bigg\vert{\frac{C}{A}}\bigg\vert > 1\ee
is satisfied.

Table 2 and the plot in fig.\ref{cplot} are clearly equivalent: if $A$ and $\Delta$ have the same sign, the asymmetry parameter is monotonically increasing with $\E$; therefore, a vertical line in the plot (corresponding to a given $C/A$) represents a sequence going from left to right in the Table. If $A$ and $\Delta$ have different signs, the asymmetry parameter is monotonically decreasing with $\E$ and a vertical line in the plot represents now a sequence going from right to left in the Table. Fig.\ref{cplot} can also be used to represent the exact resonance, $\Delta=0$. In this case, on the plane of the control parameters $C/A, -B/2A$, each point is associated with a fixed phase-space structure. 

Any point of the plane is {\it structurally stable}, except the four crossing points of the bifurcation lines. By structural stability we mean that in any points in a small neighborhood (and so for any small perturbation) of the given point, the phase-space structure does not change. The four degenerate points are structurally {\it unstable}: to understand what happens if they are perturbed, the first-order approach used up to now is not enough and we have to go to a higher order. Subsection \ref{degenz} illustrates this occurrence with a nice simple example and the Appendix adds clues on its analytic meaning.

\section{The energy-momentum map}\label{em}

So far we have not made direct use of the most relevant feature of a family of integrable systems, namely the fact that we have (indeed by construction) a complete commuting set of integrals of motion. Indeed, the surfaces of section we have used above are a standard tool with which anyone is familiar with and are directly associated to the integrability of the 2-DOF system. In this and the subsequent sections we explore in more details these integrability properties: we introduce the `energy-momentum map', which provides a global view of the phase-space of each member of the family and later on we write explicit formulas to compute action-angle variables.

The integrable dynamical system defined by the normal-form Hamiltonian 
$$K (p_1, p_2,q_1, q_2) = K_0 + K_2 + ...,$$
where $K_{2j},\;j\ge1,$ have vanishing Poisson bracket with $K_0$, determines the two-component map \cite{CB}
\ba
\E\M: && T^* \mathbb R^2 \longrightarrow  \mathbb R^2,\\
&& w \longmapsto \left(K_0 (w),K (w)\right),\\
&& (p_1, p_2,q_1, q_2) \longmapsto  \left(K_0 (p_1, p_2,q_1, q_2),K (p_1, p_2,q_1, q_2)\right).\ea
Any phase-space point $w \in T^* \mathbb R^2$ is called {\it regular} if the rank of the matrix of  the differential of the energy-momentum map
\be
d\E\M=\begin{pmatrix}
\partial_w K_0 \\
\partial_w K \\
\end{pmatrix}
\ee
is two, otherways it is called {\it critical}. The theorem of Liouville-Arnold \cite{Arn1} implies that, chosen a regular value $w$ of $\E\M$, there is a neighborhood $W(w)$ such that $\E\M^{-1} (W)$ is isomorphic to $W \times {\mathbb T}^2$. ${\mathbb T}^2$ denotes a {\it two dimensional torus} immersed in the phase-space and this confirms that the phase-space of our system is a torus-bundle with (possible) singularities. By explicitly constructing the $\E\M$ map we can assess the nature of these singularities and how they are related with the critical values of the map. At critical values where the Hessian matrix has rank less than two, it is immediate to deduce that the curves of critical values on the image of the map are connected with the bifurcation lines found above and that the pre-image of the critical values coincide with the 1-tori of the periodic orbits in general position \cite{CDHS,SD}.
%
%
%
%

For our purposes it is enough to consider the truncated form $K=K_0+K_2$ and use $K_0=\E$ as first component and 
$$
\H = K - \left(1+\Delta\right)\E - A_1 \E^2 = h
$$ 
itself as the second component. Since the rank of $d\E\M$ is zero at equilibrium and is one where the differential of the two components are linearly dependent and not both zero, these conditions for the singular values of the map correspond exactly to those obtained in the geometric analysis performed in the previous section. The $\H$ component takes its extrema just on the normal modes and on the bifurcating families. Therefore, the curves defined by \eqref{yax_energy}--\eqref{xax_energy} give the boundary branches of the image of the energy-momentum map {\it up to} the first bifurcation. In addition to these, the values of $\H$ at the contact points between the reduced phase-space and the second integral given by the functions \eqref{HU} and \eqref{HL} provide new branches starting and/or ending at bifurcating points. {\it External branches are produced by bifurcations of stable families}, while the {\it internal branches appear when bifurcations of unstable families are accompanied by the return to stability of a normal mode}. 
\begin{figure}
 \includegraphics[width=14cm]{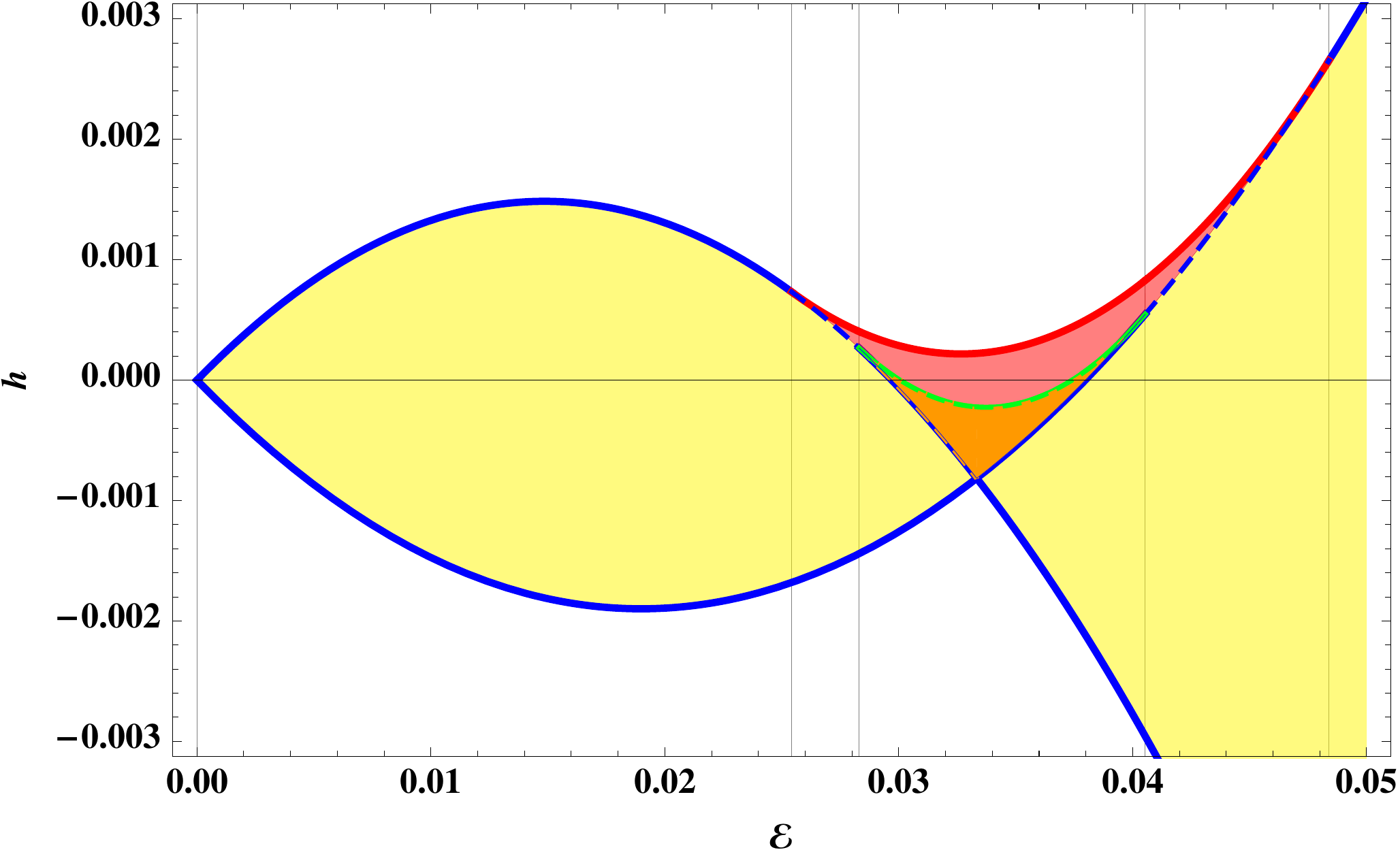}\\
 \centering{
\includegraphics[width=3.25cm]{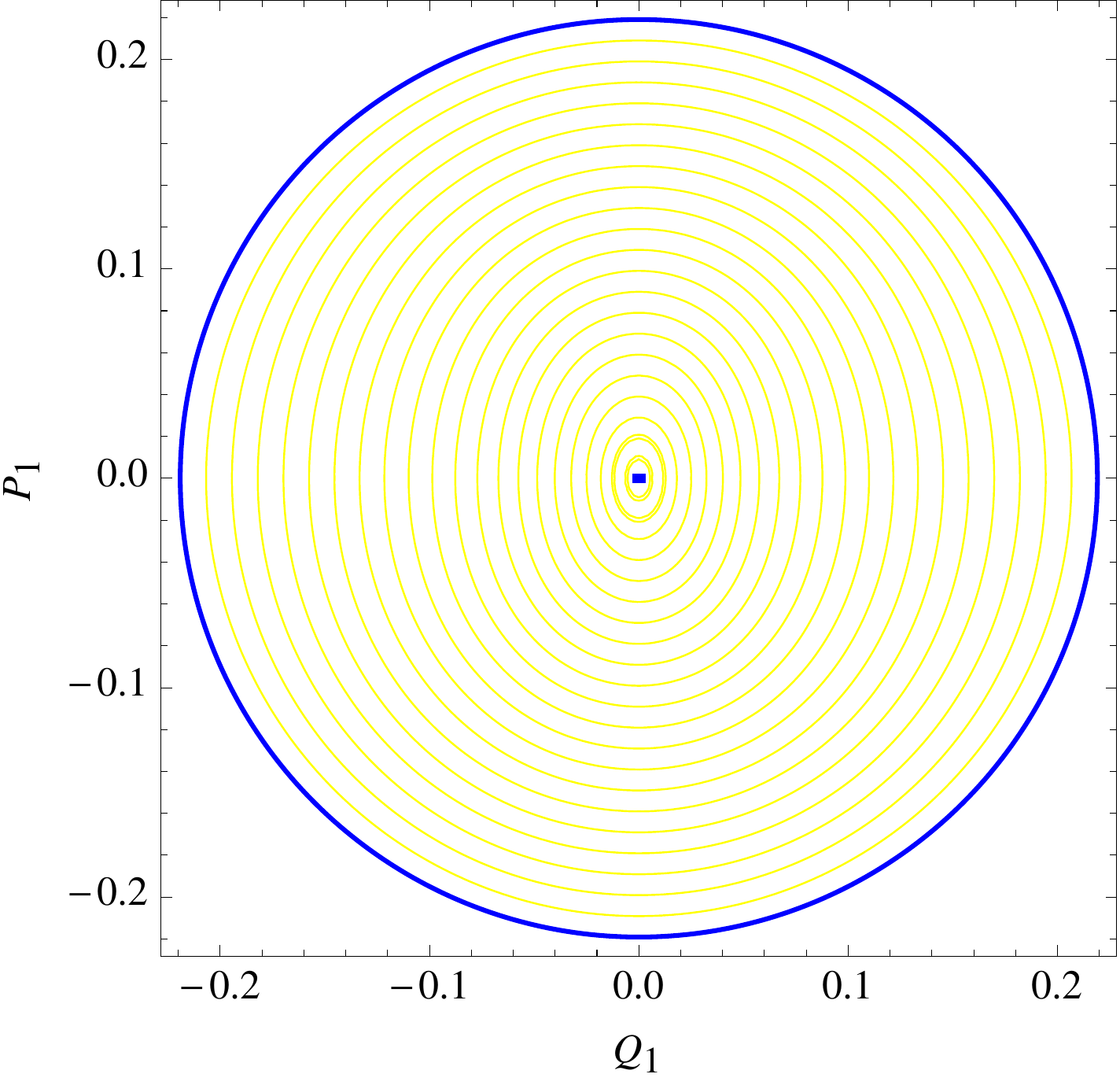}
\includegraphics[width=3.25cm]{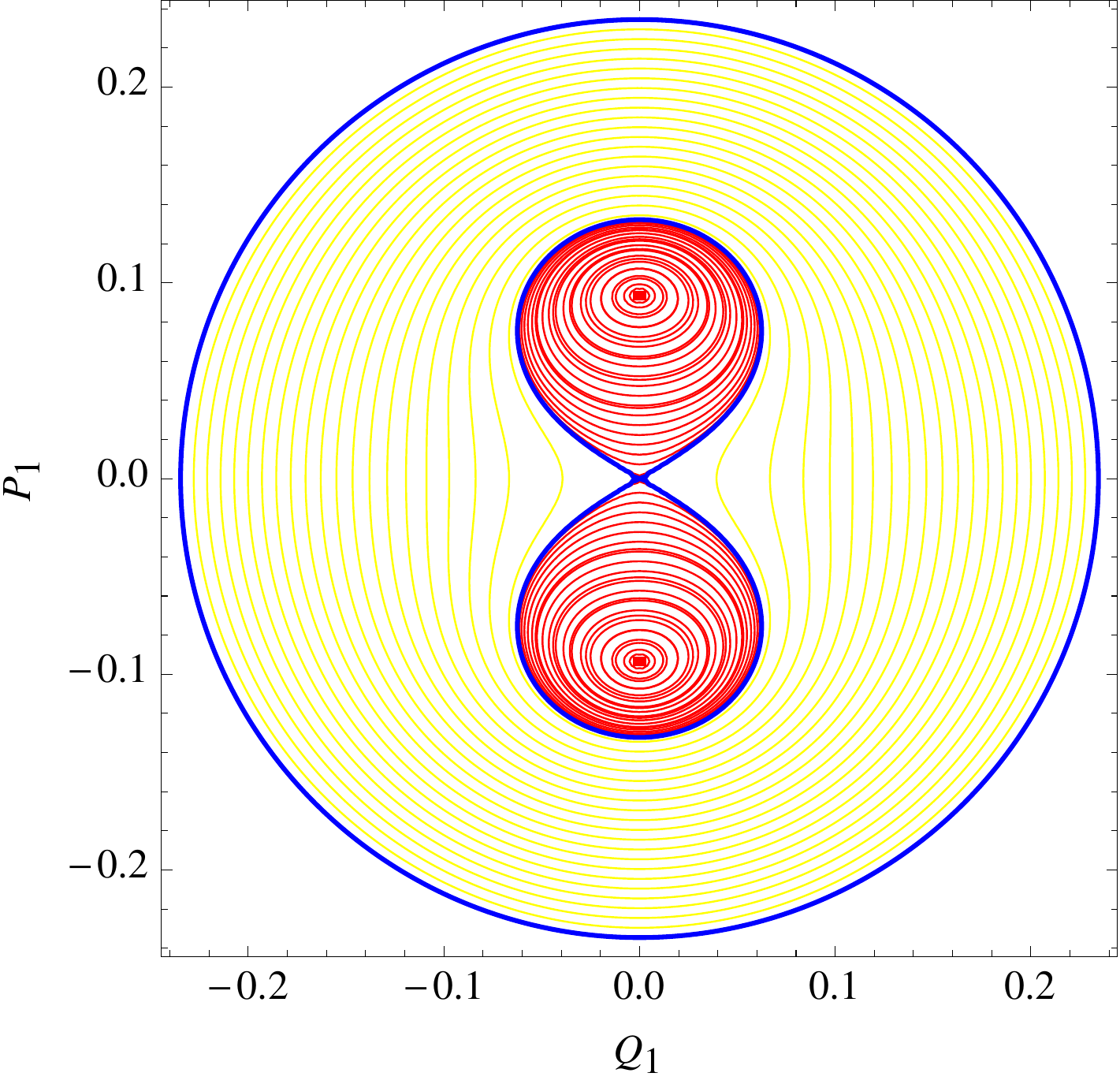}
\includegraphics[width=3.25cm]{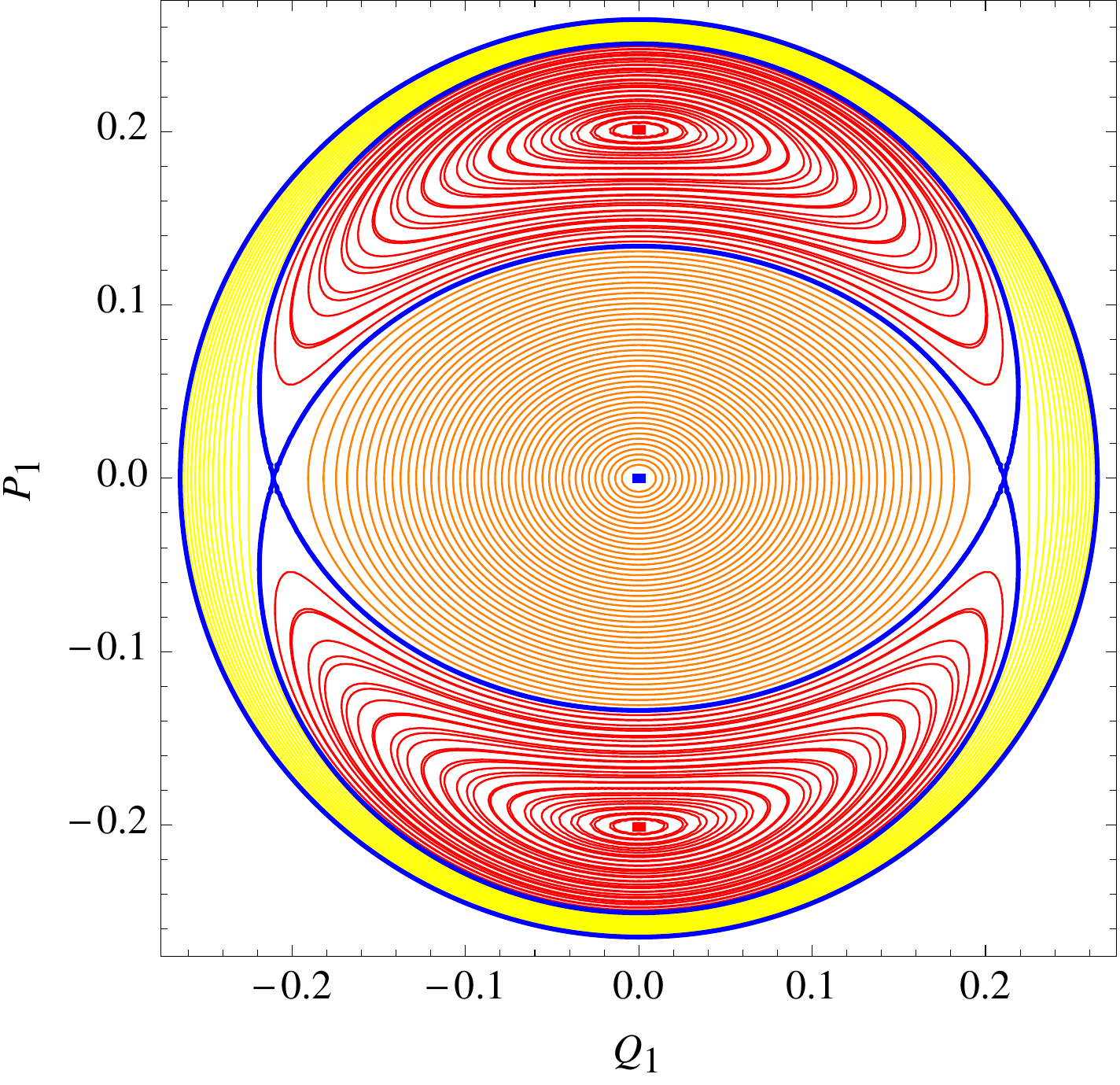}
\includegraphics[width=3.25cm]{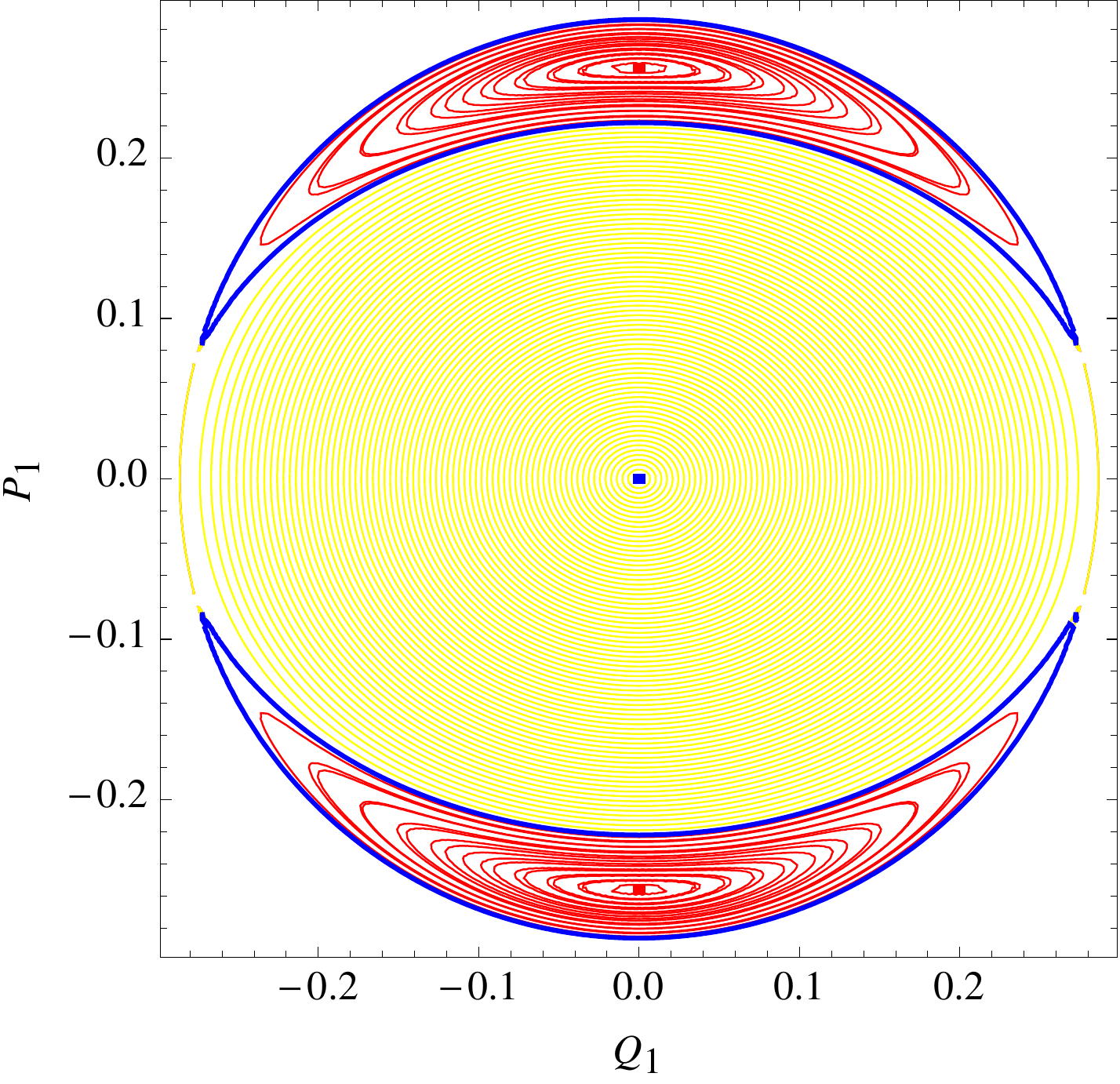}
}\\
\includegraphics[width=13.3cm]{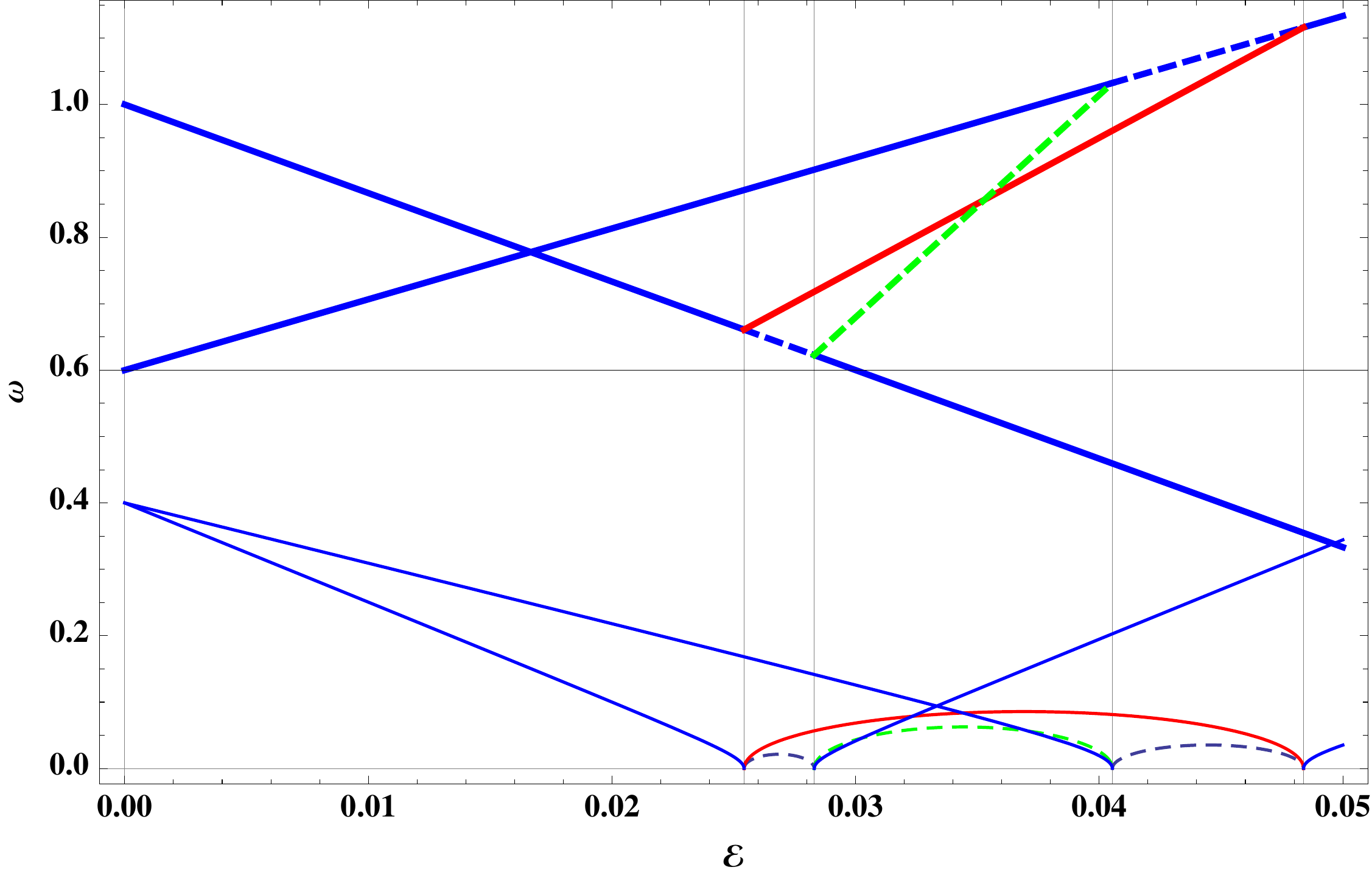}
\caption{Upper panel: Image of the energy-momentum map in the reference case, sub-case 1: $A=-\frac{11}{15},B=6,C=\frac15,\D=-\frac15$. Middle panel: corresponding surfaces of section at levels $\E=0.024,\; \E=0.028 ,\;  \E=0.035,  \; \E=0.041$.
Lower panel: Frequencies on the periodic orbits (thick lines); reduced (variational) frequencies (thin lines).}
\label{rc1plot}       
\end{figure}

All these features are nicely displayed in the bifurcation plots of the image of the map in the $(\E,h)$-plane. Let us consider for definiteness the reference case of subsection \ref{caso_particolare}. In the upper panel of fig.\ref{rc1plot} we see the plot corresponding to the first sub-case, that with $A+C<0$  (first line of Table 2, left): the vertical lines represent the sequence $\E_{1U},\E_{1L},\E_{2L},\E_{2U}$ and the range of the map is the union of the 3 domains 

\begin{itemize}
  \item $\{0\le\E\le\E_{1U},h_2 \le h\le h_1 \}$, 
  \item $\{\E_{1U}\le\E\le\E_{2U},{\rm min}(h_2,h_1)\le h \le h_U \}$,
  \item $\{\E\ge\E_{2U},h_1 \le h\le h_2 \}.$ 
\end{itemize}  
  The blue curves correspond to the two normal modes and are the only boundaries of the first and third domains. In the ranges in which the normal modes are unstable the blue curves are dashed. Invariant tori around the normal modes are shaded in yellow. The upper boundary of the second domain is the red curve which is associated with the bifurcation of the stable family of the inclined orbits. The green dashed curve is completely inside the domain and is associated with the bifurcation of the unstable family of the loop orbits. The `chamber' below it, shaded in orange, is occupied by invariant-tori around NM1 (again stable after $\E_{1L}$) turning again in yellow when NM2 becomes unstable at $\E_{2L}$. The chamber between the red and the green curves is occupied by tori, shaded in red, around the stable inclined. The plots in the middle panel resumes the situation with a suitable sequence of sections and the same color code to denote the families of tori.

\begin{figure}
 \includegraphics[width=15cm]{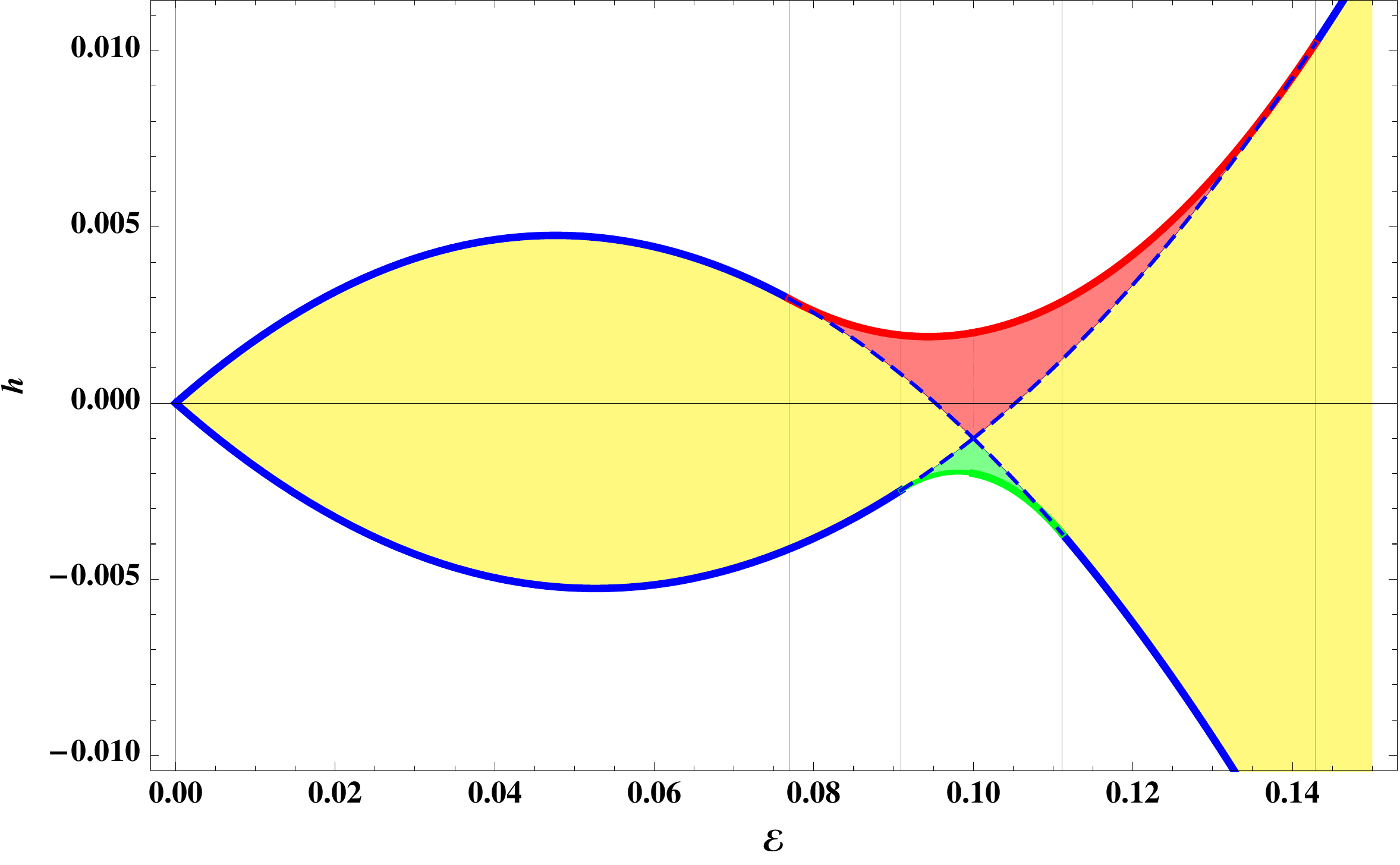}
 \centering{
\includegraphics[width=3cm]{rc1S2C.pdf}
\includegraphics[width=3cm]{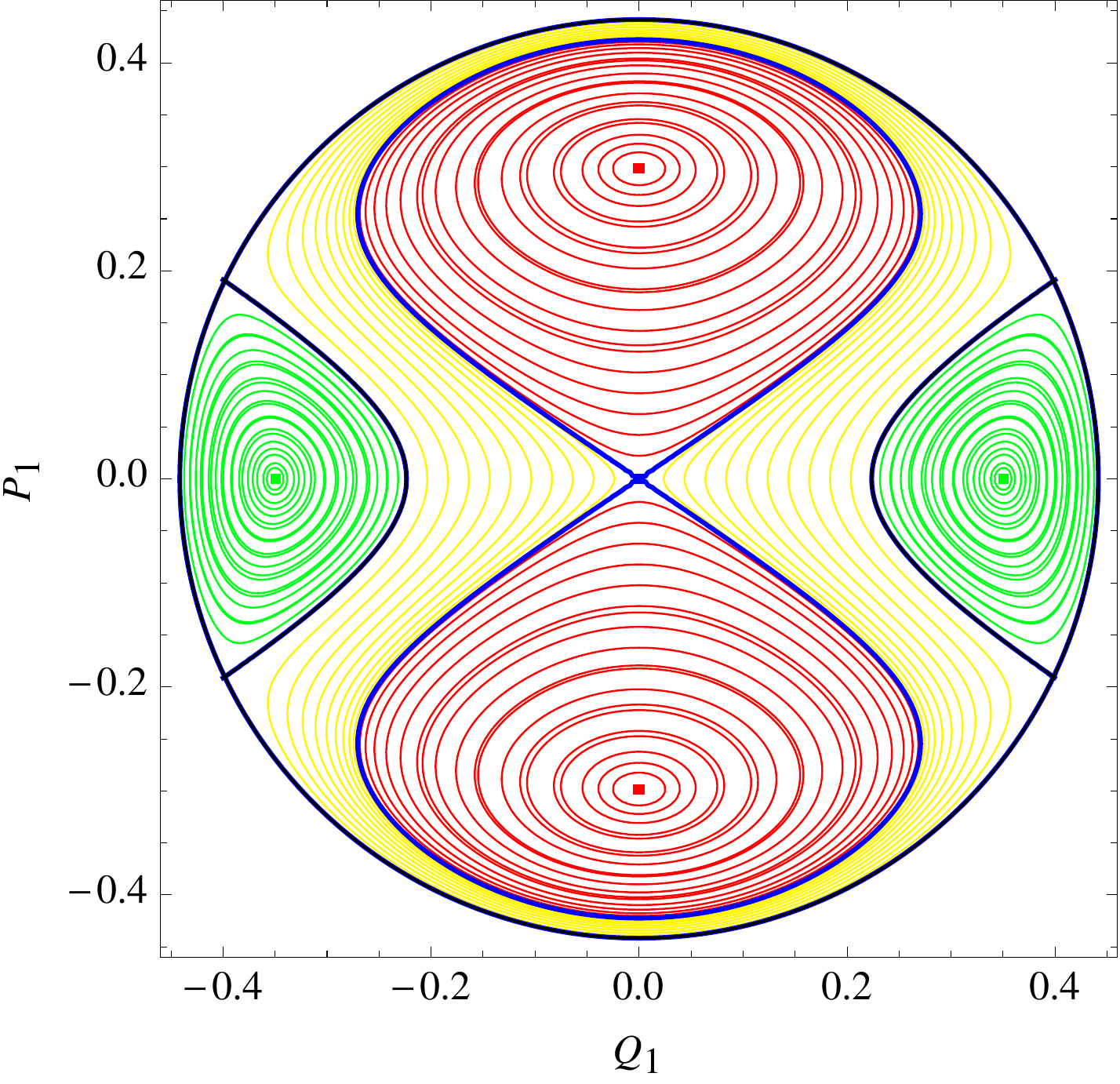}
\includegraphics[width=3cm]{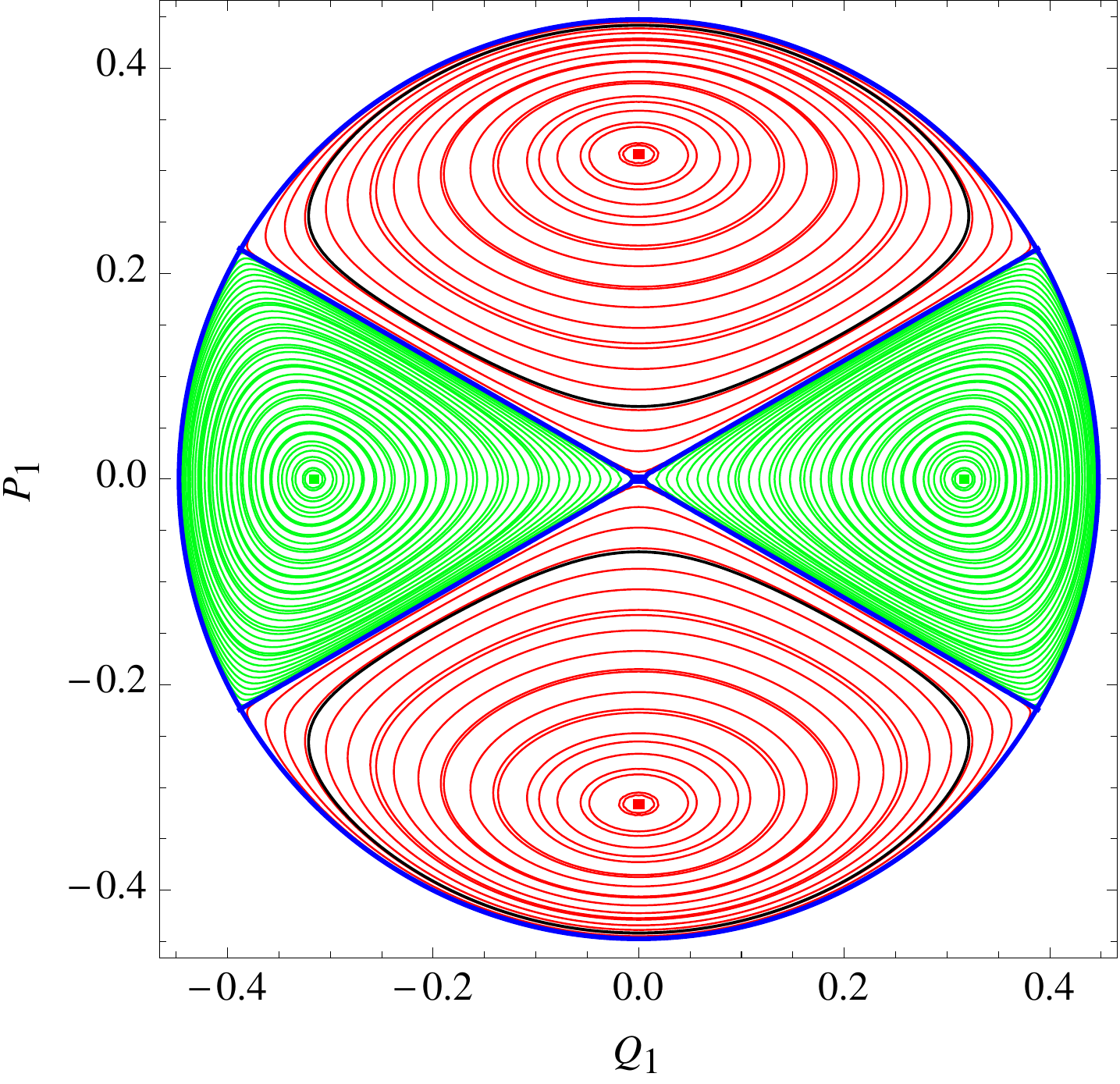}
\includegraphics[width=3cm]{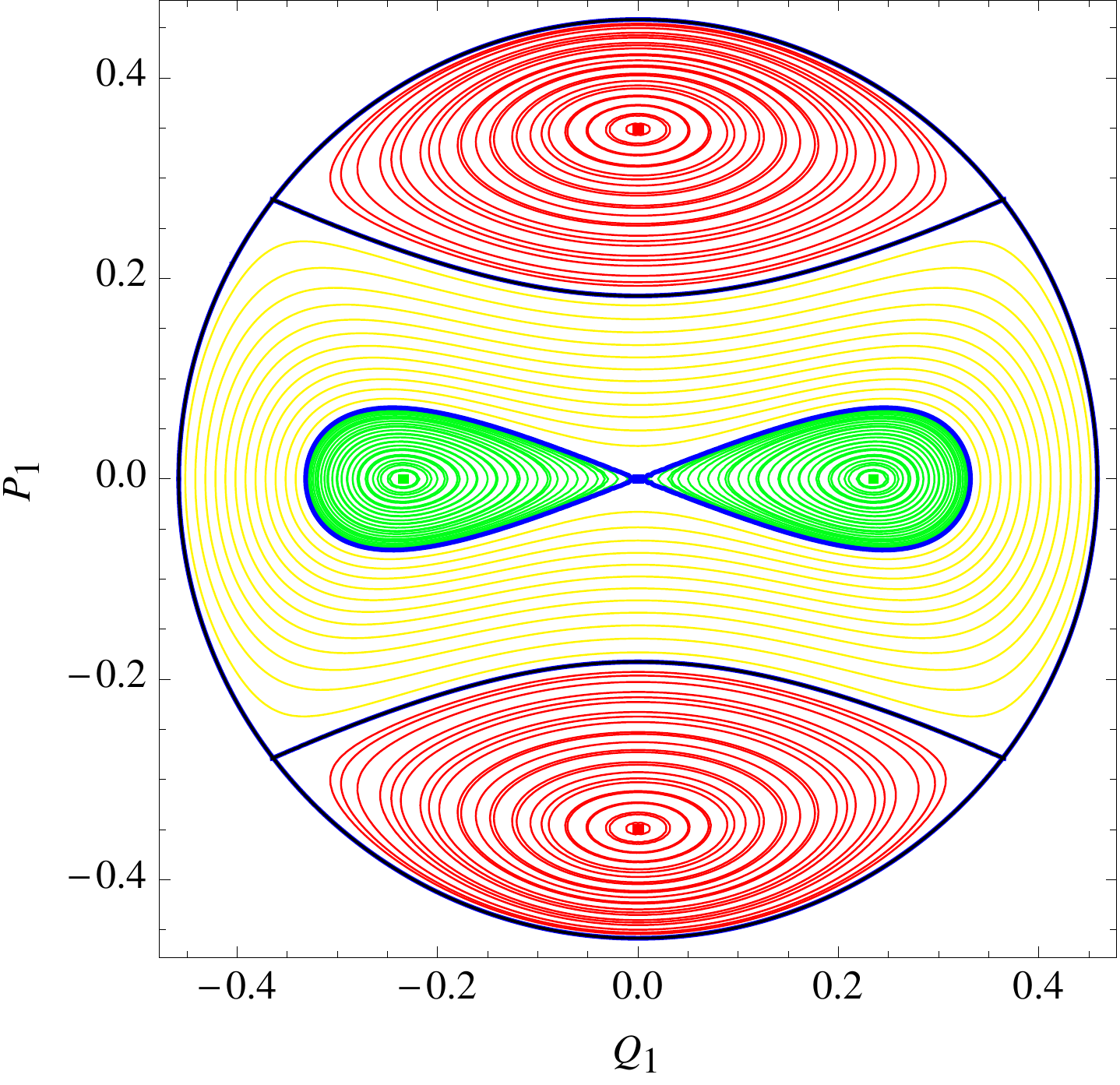}
\includegraphics[width=3cm]{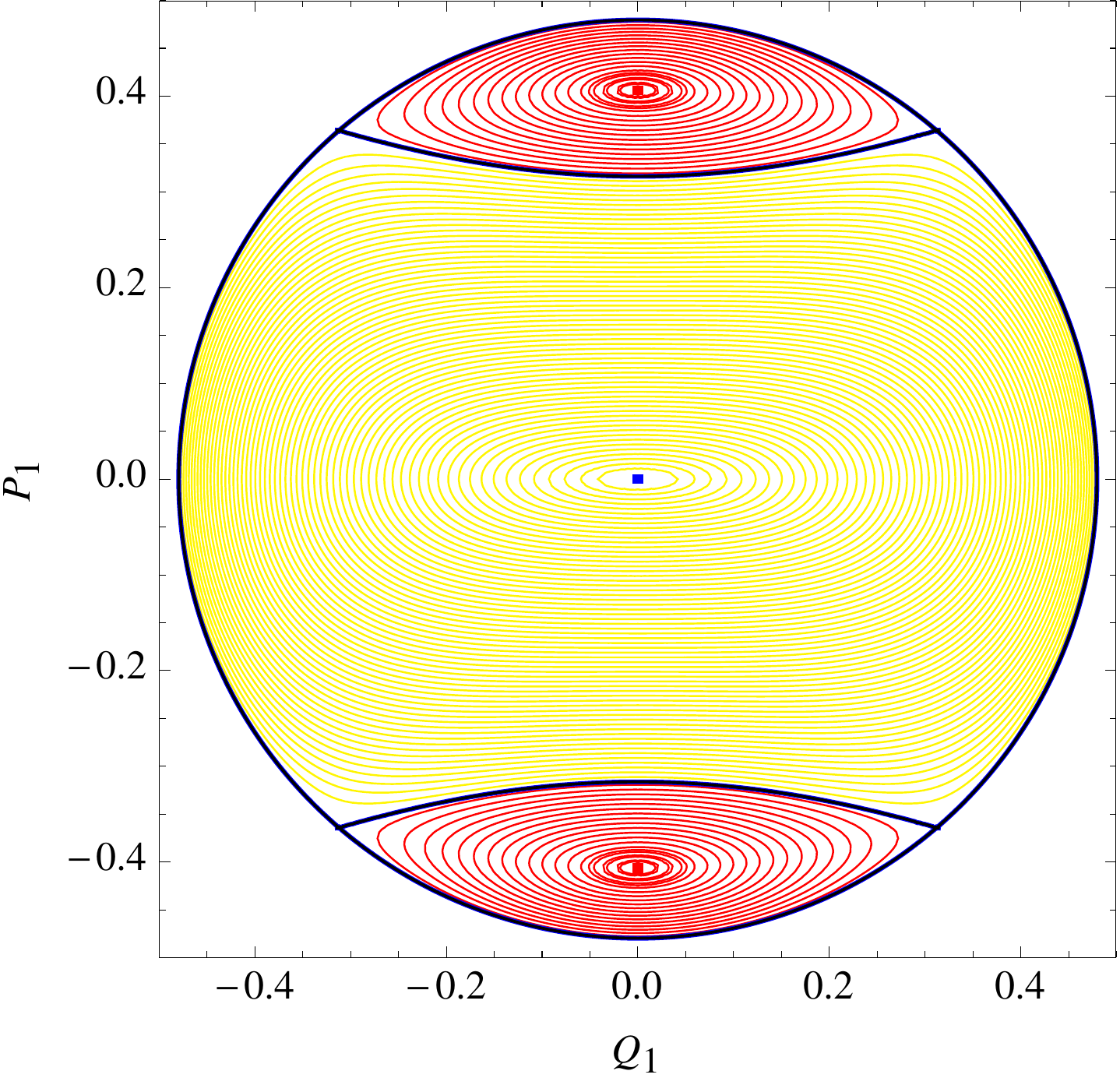}
}
\includegraphics[width=13.3cm]{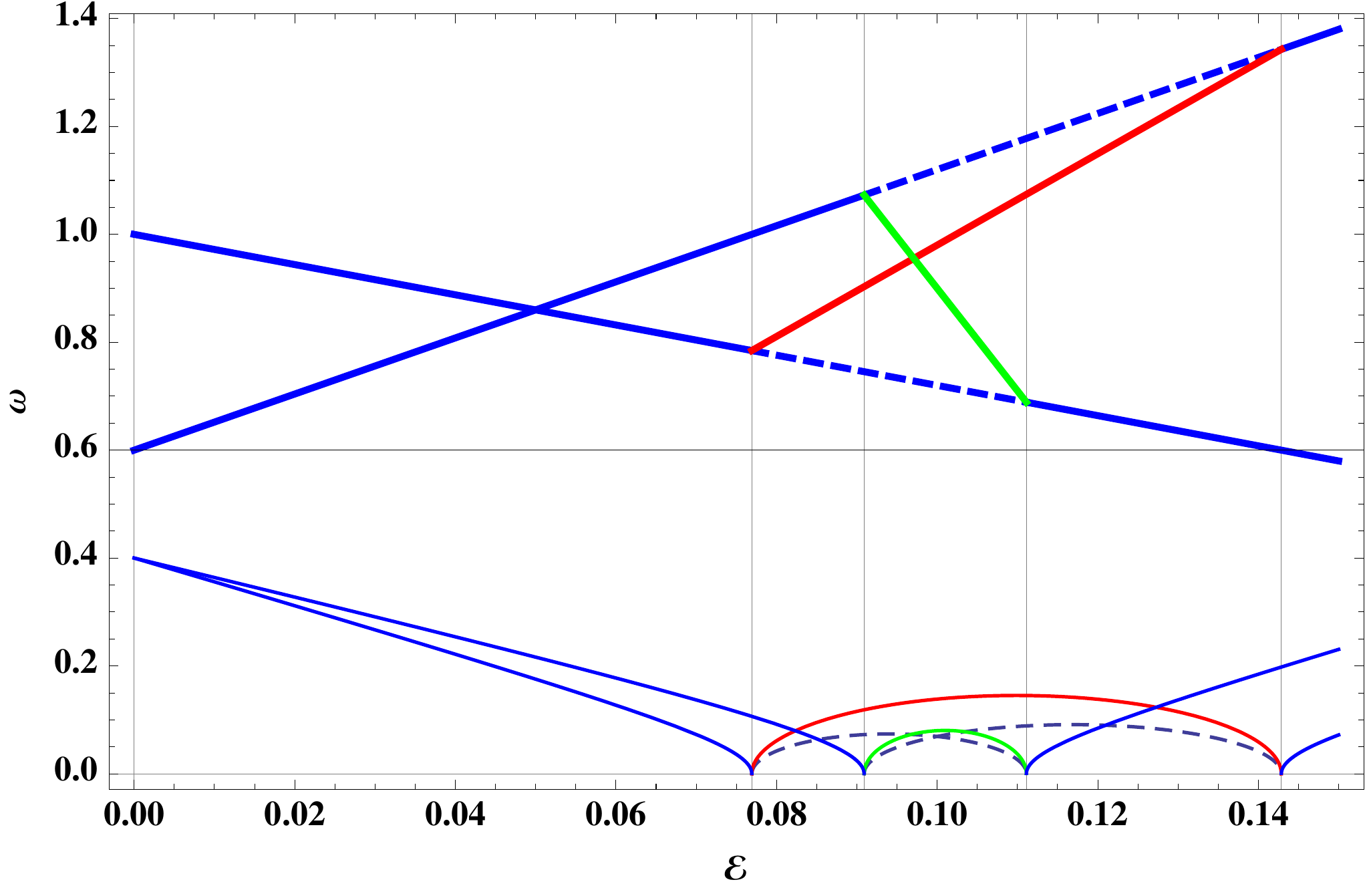}
\caption{Upper panel: Image of the energy-momentum map in the reference case, sub-case 2: $A=-\frac1{10},B=2,C=\frac15,\D=-\frac15$. Middle panel: corresponding surfaces of section at levels $\E=0.09, \; \E=0.098,  \;  \E=\E_{GB}=0.10,  \; \E=0.105,  \; \E=0.115$. Lower panel: Frequencies on the periodic orbits (thick lines); reduced or variational frequencies (thin lines).}
\label{rc2plot}       
\end{figure}

In the upper panel of fig.\ref{rc2plot} we see the plot corresponding to the second sub-case, that with $A+C>0$   (second line of Table 2, left): the bifurcation sequence now is $\E_{1U},\E_{2L},\E_{GB},\E_{1L},\E_{2U}$ and the range of the map is the union of the 5 domains 

\begin{itemize}
  \item $\{0\le\E\le\E_{1U},h_2 \le h\le h_1 \}$, 
  \item $\{\E_{1U}\le\E\le\E_{2L},h_2 \le h \le h_U \}$, 
  \item $\{\E_{2L}\le\E\le\E_{1L},h_L \le h \le h_U \}$, 
  \item $\{\E_{1L}\le\E\le\E_{2U},h_1 \le h \le h_U \}$, 
  \item $\{\E\ge\E_{2U},h_1 \le h\le h_2 \}$. 
 \end{itemize}  
  The red curve is again associated with the stable family of the inclined orbits whereas the continuous green curve is associated with the bifurcation of the now {\it stable} family of the loop orbits and the chamber {\it above} it, shaded in green, is occupied by invariant-tori generated by them. At the value $\E_{GB}$ of the global bifurcation the phase-space fraction of tori around the normal modes vanishes as can be seen in the third central section of the sequence of the central panel of the figure.

Apart for symmetry transformations, these two cases exhaust the inequivalent types of energy-momentum maps for this family of systems. In both instances the control parameters have been chosen in order to get positive values for {\it all} the thresholds: otherwise, one or more branching points are lacking and the ensuing chambers are unbounded. All other complementary cases of subsection \ref{compl} can be obtained by applying the transformation rules of table 1.


\section{Action-angle coordinates}
\label{AARN}

According to a corollary of the Liouville-Arnold theorem \cite{gior} there exists a set of action-angle variables $({\mathcal J}_{\ell},{\vartheta}_{\ell}), \ell=1,2,$ such that the Hamiltonian can be written in the form
\be\label{TAA}
{\mathcal K} = {\mathcal K} ({\mathcal J}_1,{\mathcal J}_2)\ee
and the angles evolve according to
\be\label{TAAng}
{\vartheta}_{\ell}=\omega_{\ell} (t-t_0),\ee
where the `frequencies' are found by means of the derivatives
\be\label{TFR}
\omega_{\ell} ({\mathcal J}_1,{\mathcal J}_2)=\frac{\partial {\mathcal K}}{\partial {\mathcal J}_{\ell}}.\ee
The problem of computing the actions ${\mathcal J}_{\ell}, \ell=1,2,$ is simplified by the fact that $K_0=\E$ is already one of them, ${\mathcal J}_1\doteq\E$. To find expressions for the other action, the reduced dynamics suggests to look for integrals in $Z$. In fact, the canonical variables adapted to the resonance can be slightly modified with a linear transformation such that the canonical set is given by $\E, Z, \eta_+, \eta_-$, with $\eta_{\pm} = (\phi_1\pm\phi_2)/2$. The second `non-trivial' action can therefore be computed by means of 
\be
{\mathcal J}_2 (\E,h)= - \frac1{2 \pi} \oint \eta_- dZ,\ee
where the contour of integration is the cross-section of the invariant torus fixed by $\E$ and $h$ on the $(Z,\eta_-)$-plane. By applying the linear transformation to (\ref{inv1}--\ref{inv2}) and by using the first of \eqref{tr_lem} we find
\be
\eta_-=\frac14 \arccos \left( \frac{X}{\E^2-Z^2} \right).\ee
The reduced dynamics is encapsulated in the relation \eqref{parabola} which defines the parabola $X=\X(Z;\E,h)$. Henceforth, the following integral is the formal expression for the non-trivial action
\be\label{NTA}
{\mathcal J}_2 (\E,h)=- \frac1{8 \pi} \oint \arccos \left( \frac{\X(Z;\E,h)}{\E^2-Z^2} \right) dZ. \ee

From the Liouville-Arnold theorem \cite{gior}, in view of \eqref{TAA} and \eqref{TFR}, we can compute the Jacobian matrix
 \be
\frac{\partial (\E,{\mathcal K})}{\partial \left({\mathcal J}_1,{\mathcal J}_2 \right)} = 
\left( \begin{array}{cc}
                        1 & 0 \\
                       \omega_1 & 
                       \omega_2 \\
                      \end{array}
                    \right).\ee
By computing the inverse matrix, it follows that
\be
\frac{\partial ({\mathcal J}_1,{\mathcal J}_2)}{\partial \left(\E,{\mathcal K} \right)} = 
\left( \begin{array}{cc}
                        0 & 1 \\
                       \frac{\partial {\mathcal J}_2}{\partial \E} & 
                       \frac{\partial {\mathcal J}_2}{\partial \mathcal K} \\
                      \end{array}
                    \right) = \frac1{\omega_2} \left( \begin{array}{cc}
                        \omega_2 & 0 \\
                       -\omega_1 & 
                       1 \\
                      \end{array}
                    \right).\ee
 Then, we can define the `reduced period', that is the time required to complete a cycle of the reduced Hamiltonian,
 \be\label{PR}
 T_2 =\frac{2 \pi}{\omega_2}=2 \pi \frac{\partial {\mathcal J}_2}{\partial {\mathcal K}}\ee
 and the {\it rotation number} $W$, namely the ratio between the two frequencies or, as we will see later on, ($1/ 2 \pi \, \times$) the advance of the angle conjugate to the first action in a period $T_2$:
 \be\label{RN}
 W=\frac{\omega_1}{\omega_2}=-\frac{\partial {\mathcal J}_2}{\partial \E}.\ee
By using \eqref{NTA}, the reduced period \eqref{PR} is then explicitly given by the integral
\be\label{QPR}
 T_2(\E,h) =\frac1{4C} \oint \frac{dZ}{\sqrt{Q(Z)}},\ee
 where we introduce the bi-quadratic 
 \be\label{Quartic}
 Q(Z;\E,h)=(\E^2-Z^2)^2-\left(\X(Z;\E,h)\right)^2.\ee
The rotation number is given by the partial derivative \eqref{RN}, where, in order to perform the derivation, one has to be careful by recalling the dependence encoded in \eqref{eneconv} of the reduced energy on $\E$. The result is:
\be\label{QRN}
 W(\E,h) =\frac1{8 \pi C} \oint \frac{(\E^2-Z^2)\left(1+\Delta+2 A_1 \E + BZ\right) + 2 C \E \X(Z;\E,h)}{(\E^2-Z^2)\sqrt{Q(Z)}} dZ.\ee
These integrals can be computed by extending to the complex plane and choosing a suitable contour determined by the roots of the polynomial $Q(Z)$. On periodic orbits we have double roots due to the tangency between the Hamiltonian and the reduced phase-space surfaces, therefore we obtain
\be\label{PROP}
 T_2(\E,h) = \oint_{\gamma} \frac{dZ}{(Z-Z_C)\sqrt{a_0(Z-Z_1)(Z-Z_2)}},\ee
 where $Z_C$ is the contact point, $Z_{1,2}$ the other two roots of $Q(Z)=0$ and $\gamma$ is a cycle in the complex plane around the point $Z_C$. In the reference case, the constant $a_0$ is
 \be\label{azero}
 a_0=16(C-A)(C+A);\ee
in the complementary cases a different choice of the sign can be necessary. Integrals of the form \eqref{PROP} can be computed with the method of residues. On the family of inclined, the double root is given by $Z_U$ in \eqref{QU}, so that
\be\label{PROPU}
 T_U \left(\E,h_U(\E) \right) = \frac{2 \pi i}{4C} {\rm Res} \left\{ \frac1{\sqrt{Q(Z_U)}} \right\}=
\frac{2\pi}{\sqrt{a_0(Z_{1L}-Z_U)(Z_U-Z_{2L})}},\ee
 where $Z_{1L,2}$ are the two distinct solutions of \eqref{ZQZ} evaluated at the reduced energy $h_U$ of \eqref{HU}. By explicitly computing the solutions and passing to the frequency we get
 \be\label{FNOPU}
 \omega_{2U}(\E) \doteq \frac{2 \pi}{T_U} =
2\sqrt{\frac{2C}{C-A}}{\sqrt{\left( (2(A -C) - B)\E - \D\right)\left(\D + (2(A - C) + B)\E \right)}}.\ee
 Recalling the threshold values (\ref{en_1u}--\ref{en_2u}) we see that, in the reference case, the reduced frequency of inclined periodic orbits is real in their existence range $\E_{1U}\le\E\le\E_{2U}$, coherently with its interpretation as their normal frequency \cite{CCP}. 
 
 Proceeding in an analogous manner, with $Z_L$ double root of \eqref{ZQZ}, we get
 \be\label{PROPL}
 T_L \left(\E,h_L(\E) \right) = \frac{2 \pi i}{4C} {\rm Res} \left\{ \frac1{\sqrt{Q(Z_L)}} \right\}=
\frac{2\pi}{\sqrt{a_0(Z_{1U}-Z_L)(Z_L-Z_{2U})}},\ee
 where $Z_{1U,2}$ are the two distinct solutions of \eqref{ZQZ} evaluated at $h_L$ of \eqref{HL}. Accordingly
 \be\label{FNOPL}
 \omega_{2L}(\E) \doteq \frac{2 \pi}{T_L} =
2\sqrt{\frac{2C}{C+A}}{\sqrt{\left((2(A + C) - B)\E - \D \right)\left(\D + (2(A + C) + B)\E \right)}}.\ee
From the threshold values (\ref{en_1l}--\ref{en_2l}) we again find that, in the reference case, we have to distinguish the two sub-cases $C+A<0$ and $C+A>0$: in the former, in the existence range $\E_{2L}\le\E\le\E_{1L}$, the argument of the square root is negative confirming the fact that the family of loops is unstable; in the latter, their reduced (normal) frequency is real and the family is stable. The two normal frequencies of the periodic families are plotted in the lower panels of figg.\ref{rc1plot},\ref{rc2plot} showing the connection with the stability of each family. 
 
 We can use the integral for the rotation number to compute very easily the frequency $\omega_1$ of the periodic orbit itself. For the normal modes it is straightforward to get:
 \ba
 \omega_{11}(\E) &=& 1  + 2 \Delta + 2 (A_1 + A + B) \E,\\
 \omega_{12}(\E) &=& 1  + 2 (A_1 + A - B) \E.\ea
 Let us denote for brevity with ${\cal A}(Z)$ the argument of the integral in the expression \eqref{QRN}. On the family of inclined we obtain
 \be\label{RNOPU}
 W_U \left(\E,h_U(\E) \right) = \frac{i}{4C} {\rm Res} \left\{ {\cal A}(Z_U) \right\}=
 \frac{i \left(1+\Delta+2(A_1+C)\E + BZ_U \right)} {4C} {\rm Res} \left\{ \frac1{\sqrt{Q(Z_U)}} \right\},\ee
 from which, comparing with \eqref{PROPU}, we get
 \be\label{FROPU}
 \omega_{1U}(\E) = 1+\Delta+2(A_1+C)\E + B \frac{B\E+\Delta}{2(C-A)}.\ee
 Analogously, on the family of loops we have
 \be\label{RNOPL}
 W_L \left(\E,h_L(\E) \right) = \frac{i}{4C} {\rm Res} \left\{ {\cal A}(Z_L) \right\}=
 \frac{i \left(1+\Delta+2(A_1-C)\E + BZ_L \right)} {4C} {\rm Res} \left\{ \frac1{\sqrt{Q(Z_L)}} \right\},\ee
from which, comparing with \eqref{PROPL}, we get
 \be\label{FROPL}
 \omega_{1L}(\E) = 1+\Delta+2(A_1-C)\E - B \frac{B\E+\Delta}{2(C+A)}.\ee
In the lower panels of figg.\ref{rc1plot},\ref{rc2plot} we can see the frequencies of the periodic orbits in the two sub-cases of the Reference Case.

For a generic quasi-periodic solution on a torus specified by the actions ${\mathcal J}_{\ell}, \ell=1,2,$ the angles ${\vartheta}_{\ell}$ evolve linearly according to \eqref{TAAng}. ${\vartheta}_{2}$ can be computed through \eqref{PR} when the reduced period \eqref{QPR} is known. The `fast' angle 
$$ \eta_+ = (\phi_1+\phi_2)/2$$ 
is the variable conjugate to $\E$. From the canonical equations we have
$$\dot\eta_+ = \frac{\partial K}{\partial \E} = 1 + \Delta + 2 A_1 \E + BZ(t) + 2 C \E \cos 2 \psi(t).$$
To find the first frequency we could take 
\be\label{appf} \omega_1^{(0)} = 1 + \Delta + 2 A_1 \E + BZ_0, \ee
which simply assumes that the oscillating terms have zero average. However it is clear that in general, even a tiny error on the fast frequency, rapidly produces an appreciable error. Therefore to get the exact value is necessary
to integrate to obtain
\be\label{exaf}  
\omega_1 = 1 + \Delta + 2 A_1 \E + \frac1{T_2} \int_0^{T_2} \left[ BZ(t) + 2C\E  \cos 2 \psi(t) \right] dt. \ee
To prove that this is actually the exact expression of the fast frequency we have that 
$$ \omega_1 = \omega_2 W $$
so that, by using \eqref{PR} and \eqref{QRN}, 
we get
\ba \omega_1 &=& \frac{2 \pi}{T_2} W \label{rotae} \\
&=& \frac1{4 C T_2} \oint \frac{(\E^2-Z^2)\left(1+\Delta+2 A_1 \E + BZ\right) + 2 C \E \X(Z;\E,h)}{(\E^2-Z^2)\sqrt{Q(Z)}} dZ \\
&=& \frac1{4 C T_2} \left[\left(1+\Delta+2 A_1 \E \right) \oint \frac{dZ}{\sqrt{Q(Z)}} + B \oint \frac{Z dZ}{\sqrt{Q(Z)}} + 2 C \E \oint \frac{\X dZ}{(\E^2-Z^2)\sqrt{Q(Z)}} \right] \\
&=& 1+\Delta+2 A_1 \E+ \frac1{T_2} \oint \left(B Z (t) +2 C \E \frac{\X(t)}{\E^2-Z^2(t)} \right) dt, 
\ea
which coincides with \eqref{exaf}. Eq. \eqref{rotae} confirms the statement that $W$ can be regarded as ($1/ 2 \pi \, \times$) the advance of $\eta_+$ in a period $T_2$.


\section{Applications}
\label{Applications}

We examine a selection of problems (mainly chosen in the field of astrodynamics) to illustrate the role of bifurcation theory in reconstructing the dynamics of a resonant Hamiltonian system.

\subsection{A `physical' application}\label{appl}

An interpretation in terms of a simple physical model of the classification obtained above concerns the relation between the phase-space structure and the strength of the nonlinear interaction between the two degrees of freedom. 

Condition \eqref{stabfam} provided by the catastrophe map, recalling the definitions \eqref{ABC}, is equivalent to state that, if
\be\label{interval}  
\bigg\vert \frac {2\alpha_4}{\alpha_1+\alpha_2 - \alpha_3}\bigg\vert >1\ee
   the system admits only stable bifurcating families in general position, whereas if this condition is not satisfied, the sequences contain also unstable families. In the natural case with cubic and quartic terms  (see \ref{11acq}--\ref{11dcq}) this conditions translate into
   \be\label{stabnat}
   \bigg\vert \frac {2 (b_{22} -2 b_{12}^2 + b_{12} b_{30})}
   {6 (b_{40} +b_{04}) - 4 b_{22} + b_{12}^2  + 12 b_{12} b_{30} - 15 b_{30}^2}\bigg\vert >1.\ee
   In the light of the discussion on structural stability and by applying singularity theory \cite{VU11,VU14}, the inclusion of small higher-order terms does not change these statements.

\subsection{Bifurcation of loop orbits in natural systems with elliptical equipotentials}

Let us consider a two degree of freedom system with a smooth potential with an absolute minimum in the origin, symmetric under reflection with respect to both  coordinate axes. By expanding around this equilibrium, the Hamiltonian is given by
\ba \label{Hamiltonian}
	H(\bm{v},\bm{x})=\frac12(v_1^2+v_2^2) + \sum_{k=0}^N \sum_{j=0}^{k+1} C_{2j,2(k+1-j)}x_1^{2j} x_2 ^{2(k+1-j)}
\ea 
where 
the truncation order $N$ and the coefficients $C_{ij}$ are determined by the problem under study.

In particular, we are interested in a fairly general class of potentials with self-similar elliptical equipotentials of the form
\ba
\label{pota}
{\mathcal V}_{\alpha} (x_1,x_2;\alpha,b) =\left\{
                    \begin{array}{ll}
                      \frac{1}{\alpha}\left(1 + x_1^2 + \frac{x_2^2}{b^2}\right)^{\alpha/2},\;\; & 0<\alpha<2 \\
                      \frac12 \log\left(1 + x_1^2 + \frac{x_2^2}{b^2}\right),\;\; & \alpha = 0.
                    \end{array}
                  \right.
\ea
The ellipticity of the equipotentials is determined by the parameter $b$: usually, we will speak of an `oblate' figure when $b < 1$ and a `prolate' figure when $b > 1$. The profile parameter $ \alpha $ with
\be\label{alfa}
-1 < \alpha \le 2, \ee
determines the radial profile. The non-vanishing coefficients of the expansion up to order $N=4$ are the following:
\ba
C_{20} & = & \frac12 \omega_1^2 = \frac12, \label{q20} \\
C_{02} & = & \frac12 \omega_2^2 = \frac1{2 b^{2}}, \label{q02} \\ 
C_{40} & = & \frac{\alpha-2}{8}, \label{q40} \\
C_{22} & = & \frac{\alpha-2}{4 b^{2}}, \label{q22} \\
C_{04} & = & \frac{\alpha-2}{8 b^{4}} \label{q04} \ea
Eqs.(\ref{q20},\ref{q02}) define the unperturbed frequencies so that, after the scaling transformation
\be v_1 \longrightarrow  \sqrt{\omega_1} \, p_{1} = p_{1}, \quad  v_2 \longrightarrow \sqrt{\omega_2} p_{2} = p_2 / \sqrt{b}, \quad  
x_1 \longrightarrow q_1/\sqrt{\omega_1} = q_1,  \quad
x_2 \longrightarrow q_2 /\sqrt{\omega_2} = \sqrt{b} q_2,\ee
the quadratic part of the original Hamiltonian system is cast in the `standard' form \eqref{HZO} and the normal form truncated to the first non zero resonant term is given by expression (\ref{GNF}) where
\be\label{deltalfa}
\delta = b - 1\ee 
and the coefficients (see again \ref{11acq}--\ref{11dcq}) are given by
\ba
\alpha_1 & = & \frac34 \frac{C_{40}}{C_{20} \sqrt{2 C_{02}}} = \frac3{16} b (\alpha-2),\\
\alpha_2 & = & \frac34 \frac{C_{04}}{C_{02} \sqrt{2 C_{02}}} = \frac3{16 b} (\alpha-2),\\ 
\alpha_3 & = & 4 \alpha_4 = \frac12 \frac{C_{22}}{C_{02} \sqrt{2 C_{20} }} = \frac14 (\alpha-2).\label{11dsc}\ea
By using (\ref{en_1u},\ref{en_2u},\ref{en_1l},\ref{en_2l}) and the expression \eqref{deltalfa} for the detuning, we find the thresholds to first order in the detuning
\ba \label{EPOL1}
{\cal E}_{1L} = - {\cal E}_{2L}  = \frac{4  }{2 - \alpha}(1-b), \quad
{\cal E}_{1U} = {\cal E}_{2U}  = \frac{8}{6 - 3 \alpha}. 
\ea
In view of the rescaling and of the expansion of the energy as a truncated series in the parameter ${\cal E}$, we have that $E = {\omega_2} {\cal E} = {\cal E} / q$ so to first order these are also the estimate of the `true' energy of the orbital motion. We can use the above critical values to establish the instability threshold for the model problem given by the potentials \eqref{pota}. In the range
\ba
0.7 < b < 1.3,\ea
which can be considered as `realistic' for elliptical galaxies, the thresholds (\ref{EPOL1}) give estimates correct within a 10\% if compared with numerical computations \cite{bbp06,bbp}. When the first of them is satisfied, loop orbits bifurcate from the `short-axis' (the $y$-axial normal mode in the oblate case, the $x$-axial normal mode in the prolate case \cite{MP11}). For the second bifurcation we encounter the degeneracy $A=C$ mentioned in subsection \ref{compl} which implies the absence of the family of inclined periodic orbit. To get a better precision, higher-order terms in the series expansions have been included \cite{MP13a}: if we expand the potential up to order six and truncate the normal form at $N=6$, the critical energy up to order two in the detuning parameter is given by
\begin{eqnarray}
E_{1\ell}&=&\frac{4  }{2 - \alpha}(1-q)+\frac{2 (2+3 \alpha )}{(\alpha-2 )^2} (1-q)^2 , \\
E_{2\ell}&=&-\frac{4 }{2 - \alpha} (1-q)+\frac{2 (5 \alpha -2)}{(\alpha-2 )^2} (1-q)^2 
\end{eqnarray}
respectively in the oblate and prolate case. These provide estimates correct within 3\% in the whole range of realistic values of ellipticity. It is important to remark that this analysis can be extended to the planes of symmetry of a 3-dimensional potential with reflection symmetries (as can be assumed in first approximation for elliptical galaxies): in particular, these arguments account for the important fact that, if both short and long-axis `tube' orbits exist, the intermediate-axis is always unstable. Relative to galactic dynamics, other important extensions are those of finding explicit solutions for periodic orbits \cite{contos,scu,pbbb} and investigating more general $m/n$--resonances  $(m,n \in \mathbb N)$ \cite{mes,MP13b} in which other families of periodic orbits relevant in galactic structures appear.

\subsection{`Levitation' of stars in disk galaxies.}

Consider now the case of motion in an axisymmetric 3-dimensional system which, by exploiting the conservation of axial angular momentum $L$, is effectively the Hamiltonian of a two degrees of freedom system in which $L$ is now a dynamical parameter \cite{CHR}. By using for example the family of potentials \eqref{pota} we have, using cylindrical coordinates $R,\phi,z$,
\be\label{pote}
\Phi_{\alpha} (R,z; L,\alpha,b) = \frac{L^{2}}{2{R}^{2}} + {\mathcal V}_{\alpha} (R,z;\alpha,b).\ee
These potentials have a unique absolute minimum in 
\be\label{amin}
R=R_{\rm c}(\alpha)=L^{\frac2{2+\alpha}},\;\;z=0.\ee
This is a stable equilibrium on the {\it meridional} plane $\phi = {\rm const}$ corresponding to a circular orbit of radius $R_{\rm c}(\alpha)$ of the full three-dimensional problem. Since the dynamics are scale-free we may fix the energy 
\ba
E_{\alpha} &=& \left(\frac12 + \frac{1}{\alpha}\right){\rm e}^{-\frac{\alpha}{2+\alpha}}, \;\; \alpha \ne 0, \label{ena} \\
E_{0} &=& 0, \;\;  \alpha = 0 \label{enlog}.
\ea
This implies that the radius of the circular orbit at this energy is
\be\label{circ}
R_{\rm c}(\alpha)={\rm e}^{-\frac{1}{2+\alpha}}, \;\; -1 < \alpha \le 2\ee
and we can investigate the dynamics at 
\be
E=E_{\alpha}, \;\; \forall \alpha \in (-1, 2],\ee 
by varying $L$ in the range
\be\label{amr}
0 < L \le L_{\rm max} \equiv R_{\rm c}^{\frac{2+\alpha}2} = \frac1{\sqrt{\rm e}}\ee
without any loss of generality.

In order to implement the perturbation method we expand the effective potential around the minimum (\ref{amin}). 
We introduce rescaled coordinates according to
\be\label{xy}
x_1 \doteq \frac{R-R_{\rm c}}{R_{\rm c}}, \quad x_2 \doteq \frac{z}{R_{\rm c}}\ee
with origin in the equilibrium point (\ref{amin}). The potential  (\ref{pote}) is then expanded as a truncated series (in the coordinates $x_1, x_2$) of the form
\be\label{potes}
\Phi_{\alpha}^{(N)} (x_1,x_2;L,\alpha,b) = 
\sum_{k=0}^{N} \sum_{j=0}^{k+2} C_{j,k+2-j}(L,\alpha,b)x_1^{j}x_2^{k+2-j},
\ee
where the truncation order $N$ is determined by the resonance under study and is clearly 4 in our present application to the 1:1 case. From (\ref{amin}) and the rescaling (\ref{xy}), the coefficients  of the expansion have the form \cite{P09}
\be\label{CJK}
C_{j,k+2-j} (L,\alpha,b) = L^{\frac{2\alpha}{2+\alpha}} c_{j,k+2-j} (\alpha,b).
\ee
In order to simplify formulas, we introduce the new parameter
\be\label{beta}
\beta = -\frac{2\alpha}{2+\alpha}, \;\; -1 < \beta \le 2, \ee
with the same range of $\alpha$ in view of (\ref{alfa}). The orbit structure of the original family of potentials (\ref{pote}) at the energy level fixed by (\ref{ena}) is then approximated by the orbit structure of the rescaled Hamiltonian
\be\label{HaS}
\widetilde H = \frac12 \big(v_1^{2} + v_2^{2} \big) + \widetilde \Phi_{\alpha}^{(N)} (x_1,x_2;b),\ee
where
\be\label{pots}
\widetilde \Phi_{\alpha}^{(N)} (x_1,x_2;b) = L^{\beta} \Phi_{\alpha}^{(N)} (x_1,x_2;L,b).\ee
The dynamics given by Hamiltonian (\ref{HaS}) take place in the rescaled time
\be\label{TS}
\tau = t/L^{\beta+1}\ee
at the new `energy'
\be\label{EL}
\widetilde E = L^{\beta} \left(E_{\alpha} - C_{0,0} (L,\alpha) \right) = 
  \frac1{\beta} \left(1-(L/L_{\rm max})^{\beta} \right).\ee
According to (\ref{amr}) the singular value $L=0$ is excluded from the analysis implying that the fictitious energy (\ref{EL}) is always finite and that the expansion around the equilibrium point (\ref{amin}) make sense.

The non-vanishing coefficients of the expansion of $\widetilde \Phi_{\alpha}^{(N)}$ up to order $N=4$ are the following:
\ba
c_{2,0} = & \frac{2+\alpha}2, &\label{c20} \\
c_{0,2} = & \frac1{2 b^{2}}, &\label{c02} \\ 
c_{3,0} = & -\frac{10+3 \alpha-\alpha^{2}}6, &\label{c30} \\
c_{1,2} = & -\frac{2-\alpha}{2 b^{2}}, &\label{c12} \\
c_{4,0} = & -\frac{54+11 \alpha-6\alpha^{2}+\alpha^{3}}{24}, &\label{c40} \\
c_{2,2} = & -\frac{6-5 \alpha-\alpha^{2}}{4 b^{2}}, &\label{c22} \\
c_{0,4} = & -\frac{2-\alpha}{8 b^{4}}. &\label{c04} \ea
The first two of them provide the frequencies of the {\it epicyclic motions} \cite{newast}. Recalling the time rescaling in (\ref{TS}), the radial and vertical harmonic frequencies are respectively
\be\label{rf}
\kappa = \frac{\sqrt{2+\alpha} }{ L^{\beta+1}} \ee
and
\be\label{vf}
\nu = \frac{1}{b L^{\beta+1}}. \ee


The bifurcation equations (\ref{en_1l},\ref{en_1u}) determine {\it critical} values of ${\widetilde E }$ in terms of the parameters $\delta,\alpha$ where now, in view of the values of the epicyclic frequencies, the detuning is
\be\label{detaxis}
\d=\frac{\kappa}{\nu}-1 = b \sqrt{2+\alpha}  - 1.\ee 
To make a quantitative prediction, we look for an expression for the corresponding {\it critical angular momentum} and, since the approach we have followed so far is altogether a perturbation approach truncated to the first non-trivial order, it is natural to look for expansions truncated to the first order in the detuning. Taking into account the rescaling in $\widetilde H $ and the explicit expression of the control parameters in terms of the coefficients in the expansion (see Appendix 1), the first order expansions of the critical values of the fictitious energy (\ref{EL}) are:
\ba\label{c11}
\widetilde E &=& \frac{12 (2+\alpha)}{5(-2-\alpha+\alpha^{2})}  \delta, \quad  \delta < 0,\\
\widetilde E &=& \frac{6 (2+\alpha)}{2+\alpha-\alpha^{2}}  \delta, \quad  \delta > 0.\ea
                                    The first solution corresponds to the bifurcation from the thin tube (the `vertical' normal mode), the second one corresponds to the bifurcation from the disk (the `horizontal' normal mode). By using the relation between $\widetilde E$ and $L$ established by (\ref{EL}), we get the following expressions for the critical values of the angular momentum below which inclined orbits exist:
\ba
L_{\rm crit} &=& \frac1{\sqrt{\rm e}} \left(1- \frac{24 \alpha  (b \sqrt{2+\alpha} - 1)}{5(2+\alpha-\alpha^{2})} \right)^{-\frac{2+\alpha}{2\alpha}}, 
\quad  b < \frac1{\sqrt{2+\alpha}}, \label{Lc11a} \\
L_{\rm crit} &=& \frac1{\sqrt{\rm e}}\left(1+ \frac{12 \alpha (b \sqrt{2+\alpha} - 1)}{2+\alpha-\alpha^{2}} \right)^{-\frac{2+\alpha}{2\alpha}}, 
\quad  b > \frac1{\sqrt{2+\alpha}}. \label{Lc11b} \ea
It is also useful to write the limiting case of the logarithmic potential ($\alpha=0$):
\ba
L_{\rm crit} &=&  {\rm e}^{-\frac{29}{10} + \frac{12}{5} \sqrt{2} b}, \quad  b < \frac1{\sqrt{2}}, \label{LcL11a} \\
L_{\rm crit} &=& {\rm e}^{\frac{11}2 - 6 \sqrt{2} b}, \quad  b > \frac1{\sqrt{2}}.\label{LcL11b}\ea 
A comparison with the outcome of numerical determinations of the bifurcation threshold allows us to evaluate the accuracy of these analytical predictions. 
%
The accuracy is particularly good when the model is close to the exact resonance. Overall, the discrepancy linearly grows with detuning, as can be expected in this first order approach and in all cases available in the literature is of the order of 10\% or lower \cite{P09}. 
An interesting application of these results to the dynamics of disk galaxies consists in the fact that the threshold \eqref{Lc11b} gives an {\it upper limit} of the angular momentum that ensures the stability of the planar orbits in the symmetry plane. When stable, disk orbits are accompanied by small amplitude quasi-periodic oscillations. When the angular momentum is lower than the critical value or, equivalently, if the eccentricity of the disk orbit is great enough, a double family of inclined orbits bifurcates and, in the 3-dimensional space, give rise to two stable axisymmetric `ribbons' (accompanied by their quasi-periodic companions) that may go quite far from the symmetry plane: stars on them can therefore `levitate' quite high with respect to the galactic disk.

\subsection{Bifurcation of `halo' orbits around the collinear Lagrangian points.}

The Lagrangian points are intriguing solutions of the restricted three body problem of celestial mechanics \cite{bp1,Alebook}. If seen in the reference frame rotating with the primaries, they are unstable equilibrium points: the first one $L_1$ is located between the primary, the other two $L_{2,3}$ on their side, all three remaining on a straight line connecting them and the primaries in the rotating frame. All of them are unstable. In fact, linearizing the Hamiltonian around the equilibrium, $H_0$ turns out to be of the form
\be\label{K0LP}
 H_0=\lambda q p + \omega_1 J_1+\omega_2 J_2
 \ee
 in a suitable set of coordinates, where
 \ba
 \lambda&=&\sqrt{ \frac12\left(-2 + c_2 + \sqrt{9 c_2^2-8 c_2} \right)},\\
\omega_1&=&\sqrt{ \frac12\left(2 - c_2 + \sqrt{9 c_2^2-8 c_2} \right)}, \\ 
\omega_2 &=& \sqrt{c_2} 
 \ea
 and $c_2(\mu)$ is the first of a set of real coefficients (uniquely determined by the mass ratio $\mu$ of the primaries) in the expansion of the effective potential around the given Lagrangian point \cite{JM}.
The instability is however quite `mild', in the sense that, the real number $\lambda$ which gives the inverse of the time-scale of the unstable mode, is quite small for all values $0<\mu\leq1/2$. This means that the instability occurs `slowly' and can be easily controlled with space manouvers, so that a spacecraft can be placed nearby one of the collinear points and kept there for a long time exploiting them as privileged locations for space missions. Moreover, the two real frequencies $\omega_{1,2}$, accounting for oscillations in the `planar' and `vertical' directions, are quite close one to each other. We can immediately deduce a dynamics in which, if we neglect the instability, a 1:1 resonance appears with the detuning defined as
$$
\d=\frac{\omega_1-\omega_2}{\omega_2}.$$
A normal form associated to $H_0$ can be constructed along the line sketched above and has the general form (see \cite{CPS}):
\be\label{zero}
K(P,Q,J_1,J_2,\phi_1,\phi_2)=\lambda QP + J_1 + J_2 + \d J_1 +
\sum_{n=1}^N K_{2n}(QP, J_1, J_2,\phi_1-\phi_2) .
\ee
It has the nice property of depending on the unstable degree of freedom only through monomials of the form $(PQ)^n$, so that, from the canonical equations, it is easy to see that, solutions with initial conditions $P=0,Q=0$, preserve this condition for any time. In the theory of dynamical system this procedure is referred to as {\it central manifold reduction} \cite{AbMa}.
This restriction to the central manifold provides us with a 2-DOF Hamiltonian that, truncated to the lowest order, has the standard form \eqref{GNF} to which we can apply our general theory of bifurcations. When the horizontal normal mode (the so called {\it Lyapunov planar orbit}) exceeds a critical amplitude, it becomes unstable and a pair of inclined orbits bifurcates. These orbits, known as {\it haloes} \cite{Richardson,GM}, bring a spacecraft to `levitate' quite far from the symmetry plane making them ideal for observational purposes. The agreement between the theory and numerical computations is very good, especially in the cases of $L_1$ and $L_2$ \cite{CCP} providing a way to predict in a first approximation energy and initial conditions of the family. A useful generalization is that in which one includes the Solar radiation pressure on a `Solar sail' \cite{Sara} which has the effect of lowering the bifurcation threshold easing the procedure of insertion in a suitable halo orbit.

\subsection{`H\'enon-Heiles' case}\label{degenz}

\begin{figure}[htbp]
\centering
{\includegraphics[width=6.5cm]{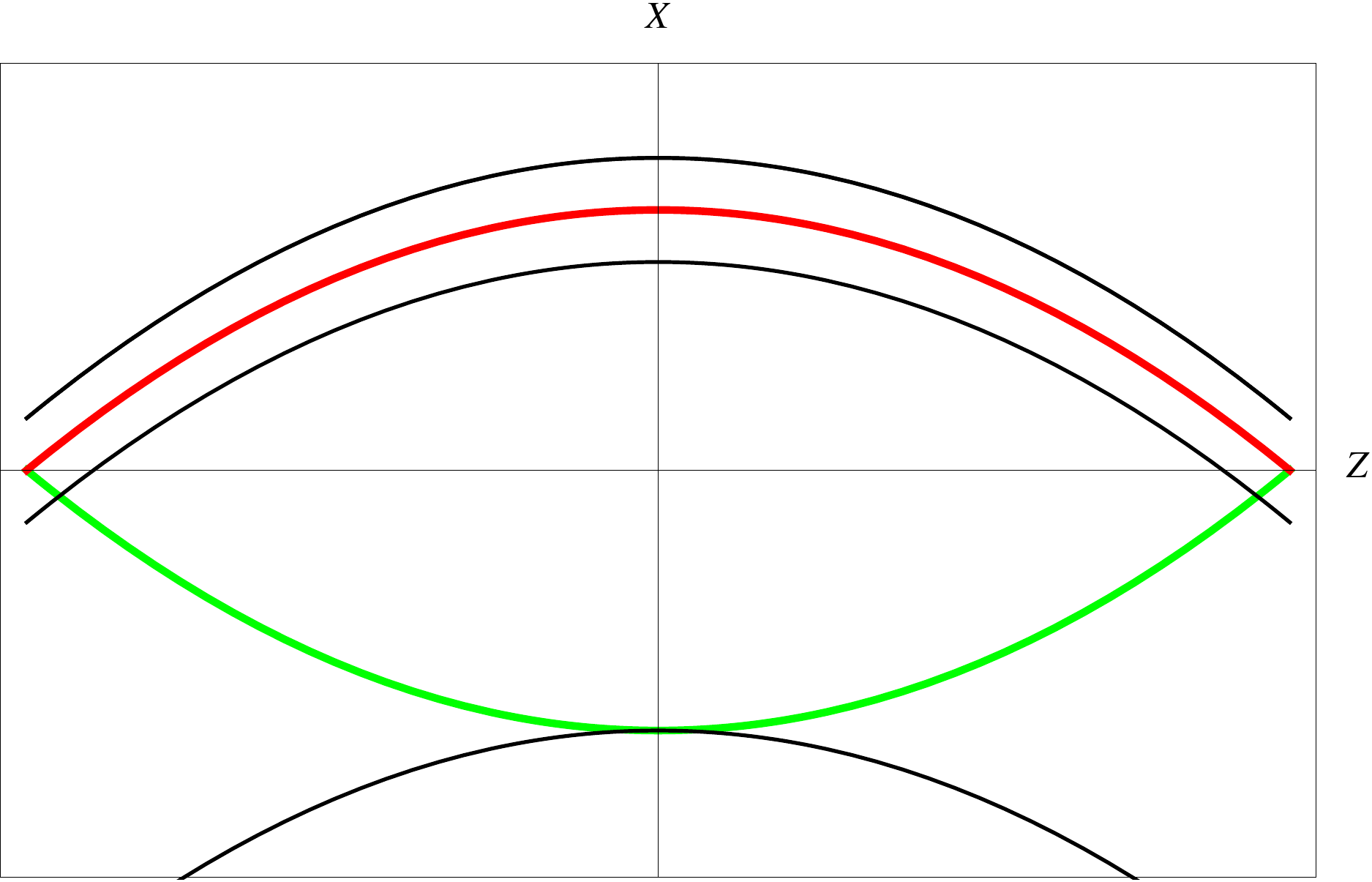}}
{\includegraphics[width=6.5cm]{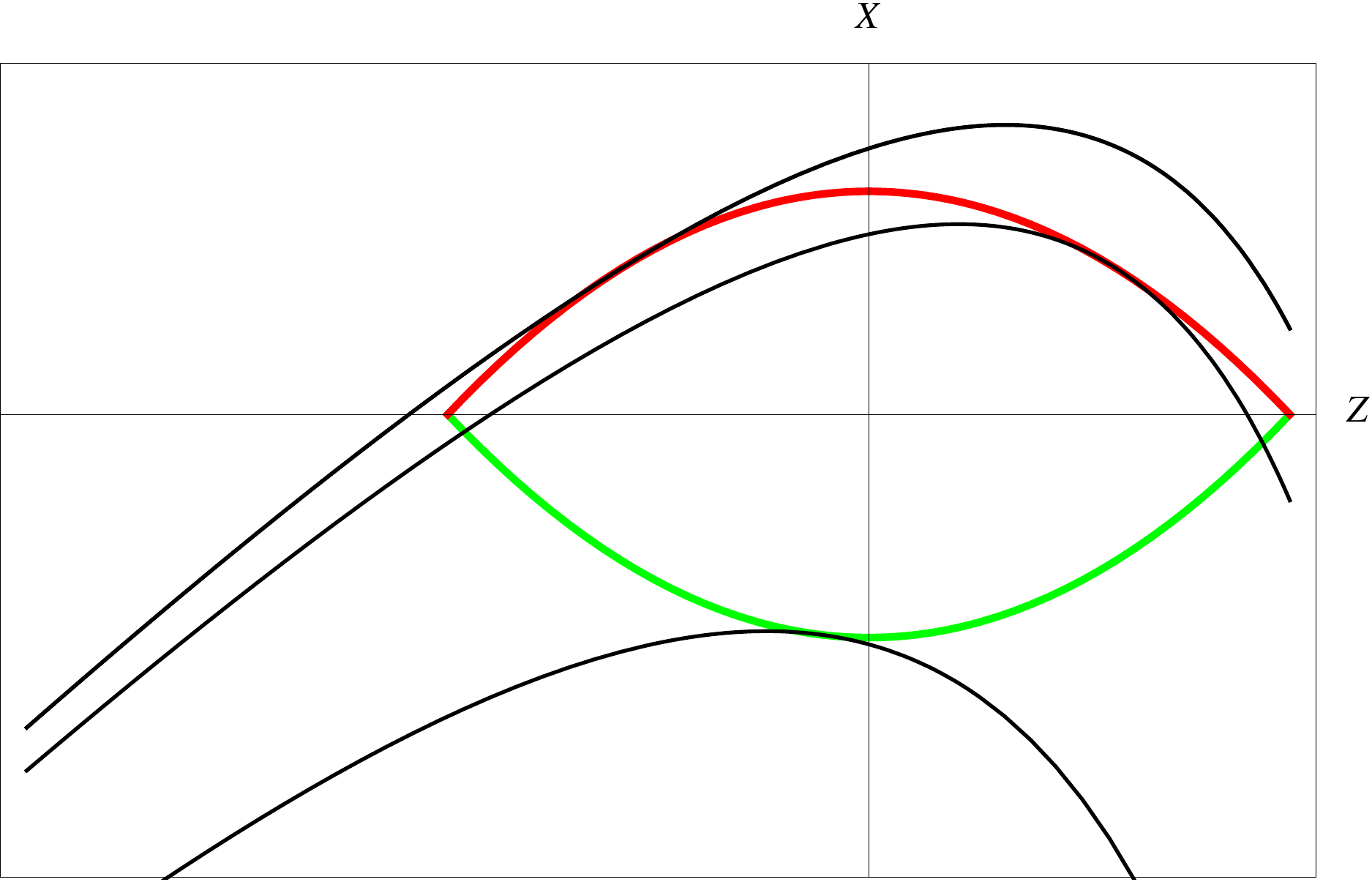}}
\centering{
{\includegraphics[width=6cm]{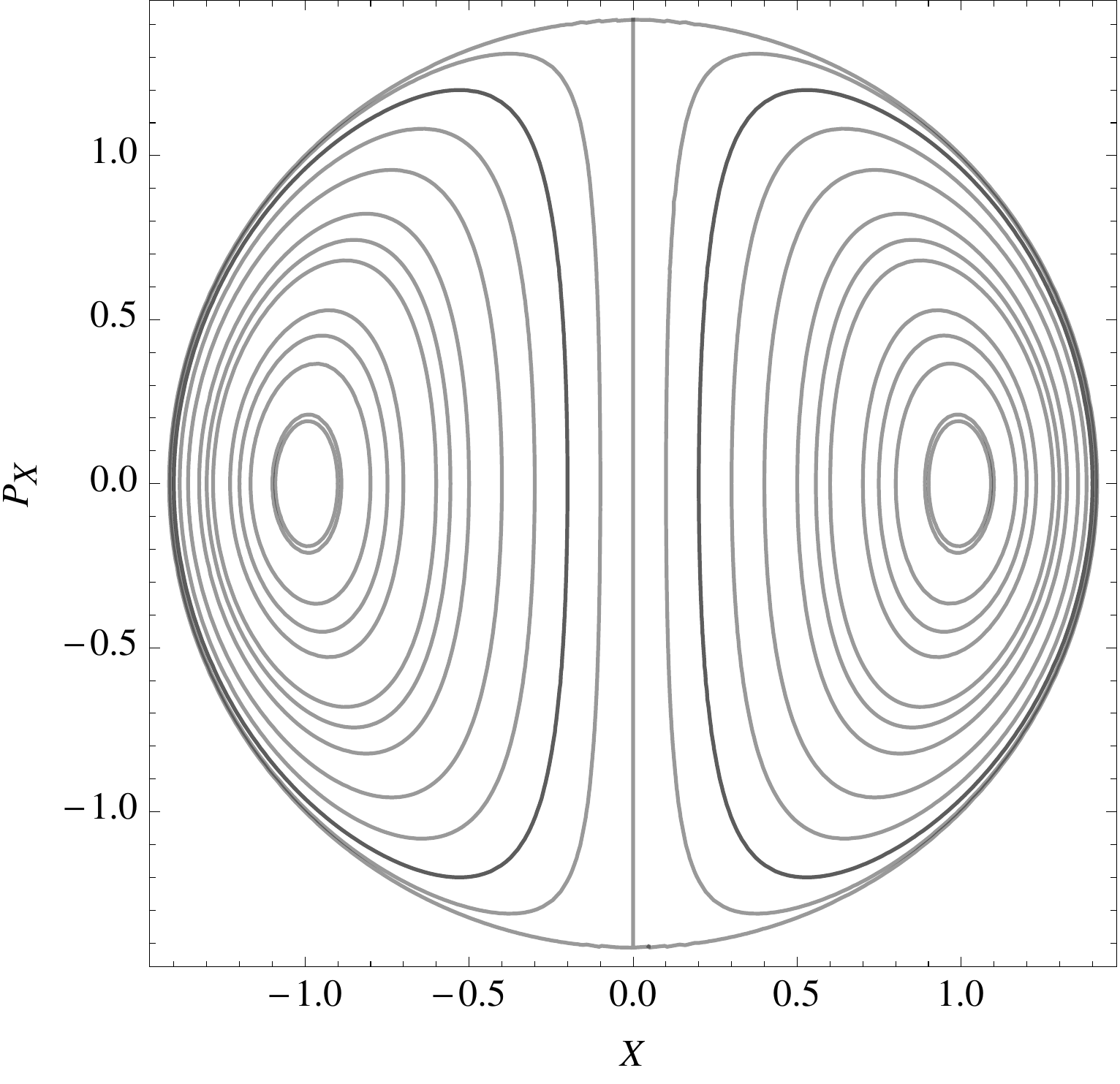}}
{\includegraphics[width=6cm]{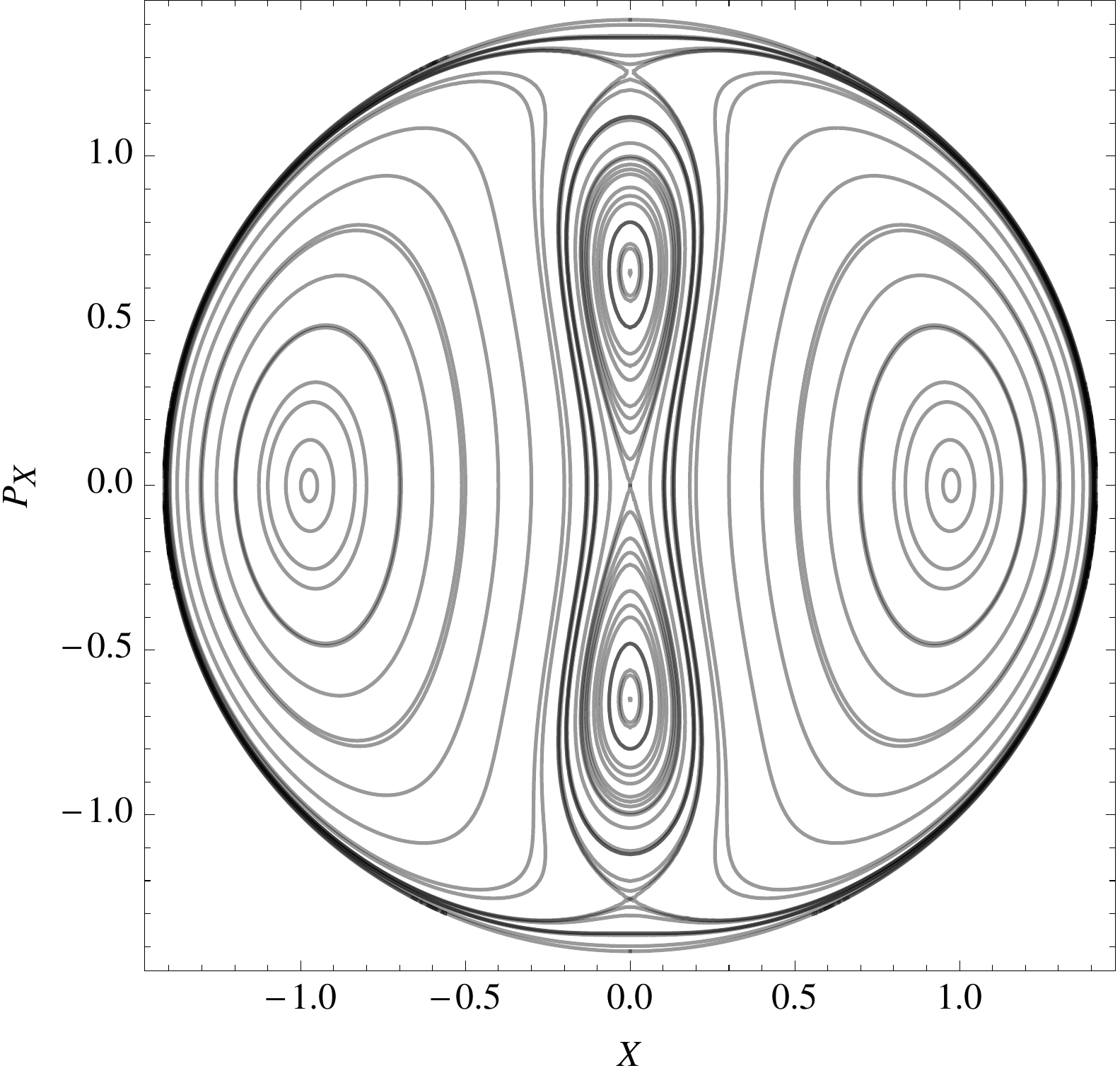}}
}
\caption{\small{Left: $\lambda=0$; right: $\lambda=0.4$}}
\label{dplot}
\end{figure}

\begin{figure}[htbp]
\centering{
{\includegraphics[width=8cm]{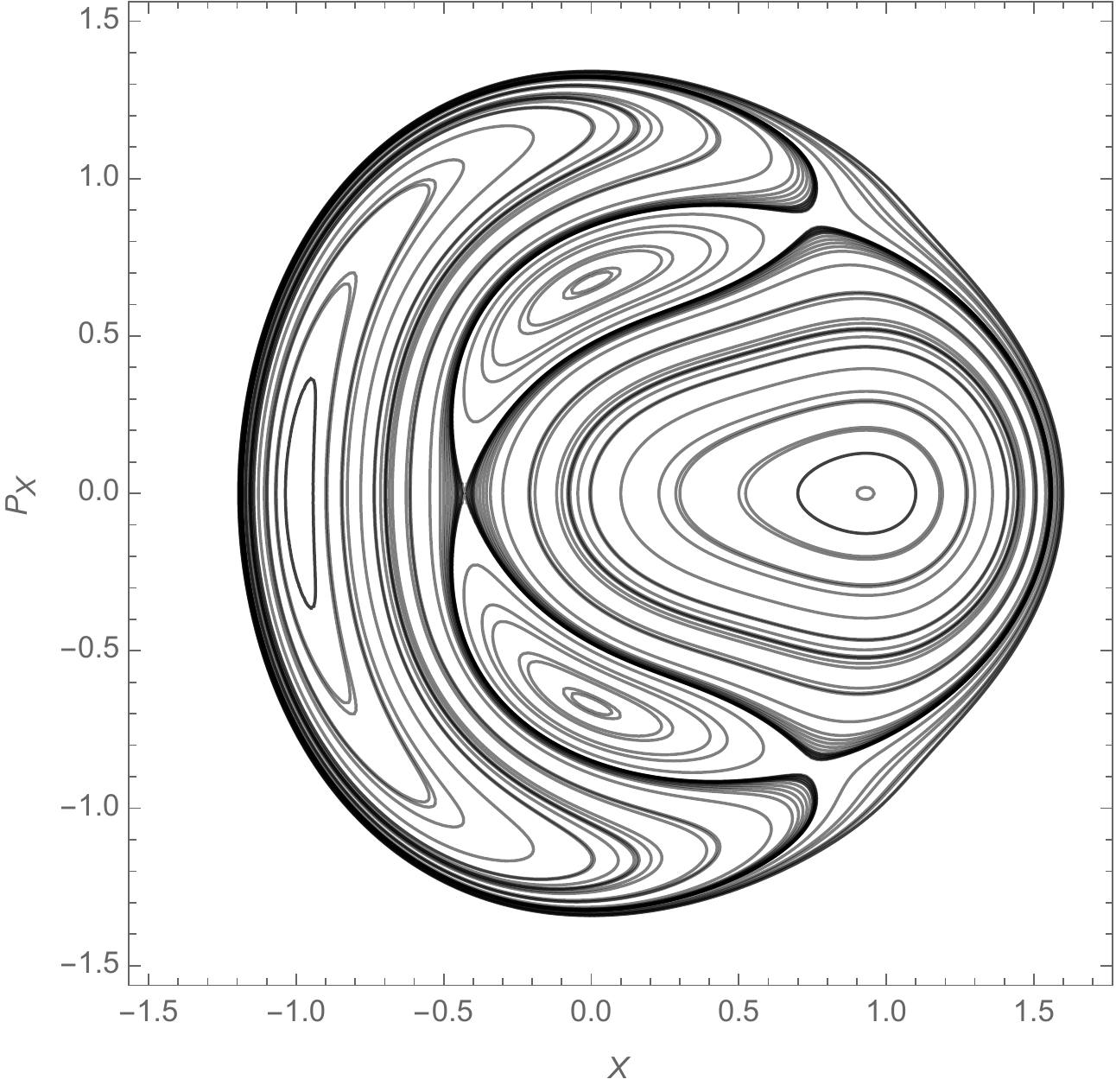}}
}
\caption{\small{The `standard H\'enon-Heiles' case.}}
\label{HHplot}
\end{figure}
Let us turn back to the degenerate point $(1,0)$ in the diagram of Fig.\eqref{cplot}. It corresponds to a `H\'enon-Heiles-type' singularity: the phase-space structure and phase portrait looks like in the upper panels of Fig.\ref{dplot}. With the choice of variables like in \eqref{tr_lem}, that is
\ba
Z & = & J_1 - J_2, \\
X & = & 4 J_1 J_2 \cos 2 (\phi_1-\phi_2),
\ea
 the Hamiltonian at order two has the very simple form 
\begin{equation}\label{hamiltonian_dege}
{\mathcal K}_D(X,Z)=C (X + Z^2)
\end{equation}
where we have also discarded terms depending only on $\E$. The unstable character of this degeneracy is manifest in the line of degenerate critical points appearing on the surface of section \cite{CHH}. A simple choice to stabilize the unfolding is the following
\begin{equation}\label{6deg}
{\mathcal K}_4(X,Z)=\lambda  Z X,
\end{equation}
with $\lambda$ constant. It is in fact clear that, with the assumed symmetries (no dependence on $Y$) and exploiting relations between the invariants, other cubic (in the action variables) terms should be of the form already present in the quadratic terms. Calling 
$$ {\mathcal K}_D + {\mathcal K}_4 = h, $$
this produces the set of `distorted' parabolas
\begin{equation}\label{distorted}
 X (Z; C,\lambda,h) = \frac{h - C Z^2}{C + \lambda Z}.
 \end{equation}
With the choice $\lambda=0.4$ (actually this $\lambda$ is not so small, but is for a clearer picture) we get the situation depicted in the lower left panel of Fig.\ref{dplot}: the distortion (no matter how small it is) produces two contact points with the upper arc, the new one giving two new hyperbolic points. The phase portrait looks like in the right panel. To get the `standard' very well known H\'enon-Heiles phase portrait (see Fig.\ref{HHplot}) it is enough to apply the recipe implicit in \eqref{HK}, namely
\begin{equation}\label{KHHH}
     H(\bm{p},\bm{q}) =  {\rm e}^{-{\cal L}_G} {\mathcal K}(\bm{P},\bm{Q}).\ee 
     One recognizes the familiar surface of section having the same topology of that given by our model problem `distorted' by the normalizing transformation.

 \appendix{Versal deformation of the 2:2 resonance}
\label{VD}

An important result in the framework of singularity theory is that of inducing a generic function, defined around a critical point and depending on several parameters, from a simple germ and deformation depending on a small set of derived parameters \cite{Br1:1,GB1,GB2,HDS,Hlibro}. In the present case, starting from the general setting introduced in \cite{Br1:1}, a versal deformation of the family of systems \eqref{GNF} is obtained in  \cite{VU11}. The easiest way to perform this further normalization is by exploiting the planar reduction and use the stereographic coordinates \eqref{ccoord}. Let us consider the resulting normal form
\ba
{\mathcal K}_b(x,y; \E,\d,\alpha_i^{(j)}) = \E + K_2 (x,y; \E,\d,\alpha_i^{(1)}) + ... + K_{2N} (x,y; \E,\d,\alpha_i^{(N)}).
\ea
It can be shown \cite{Br1:1,VU11,VU14} that there exists a $\mathbb{Z}_2\times\mathbb{Z}_2$-equivariant transformation   which `induces' ${\mathcal K}_b$ from the function 
\be
 F(x,y,u_k)=\e_1x^4+(\mu+u_3)x^2y^2+\e_2y^4+u_1x^2+u_2y^2,\label{F_uni}
 \ee
 namely, there exists a diffeomorphism
 \be
\Phi: \mathbb R^2 \times \mathbb R^{m+2} \longrightarrow  \mathbb R^2\times \mathbb R^3,\quad
 (x,y,\E,\d,\alpha_i^{(j)}) \longmapsto  \left(x,y,u_k \right),\ee
where $m$ is the dimensionality of the external-parameter space, such that ${\mathcal K}_b = F \circ \Phi.$

The coefficients $u_k, k=1,2,3,$ depend on the internal $\E,\d$ and external $\alpha_i^{(j)}$ parameters and are constructed in an algorithmic way with an iterative process carried out up to order $N$. Explicit expressions for $N=2$ are computed in \cite{VU11}. The coefficients $\e_1,\mu,\e_2$ are otherwise determined by the leading-order terms `at the singularity' $\E=\d=0$ and are expressed as the discrete set of constants
\be\label{germ}
\mu=\frac{2 A}{\sqrt{|A^2-C^2|}}, \;\;
\e_1 = \frac{A-C}{|A-C|}, \;\; \e_2 = \frac{A+C}{|A+C|},
\ee
where the coefficients $A,B,C$ have been introduced in \eqref{ABC}.
The function $F(x,y)$ provides the phase portraits on either surfaces of section of the normal form as they are determined by varying the parameters. Quantitative predictions for bifurcations around the resonance are given by the series expansion of the $u$ coefficients in terms of the internal parameters. If we content ourselves with qualitative aspects, these predictions are already determined by their first order expressions
\begin{eqnarray}
u_1&=& \frac{\Delta+(B-2(A-C))\E}{\sqrt{|A-C|}},
\label{u1_11}\\
u_2&=& \frac{\Delta+(B-2(A+C))\E}{\sqrt{|A+C|}},
\label{u2_11}\\
u_3&=&0,
\ea
where $
\Delta=\d/2.$
 
The quartic terms of the function $F(x, y)$ (with $u_3 = 0$ and coefficients like in (B.4)) compose the \emph{central singularity} of this resonance and the remaining terms give its \emph{universal deformation}. For any values of the control parameters such that $\mu^2 \neq 4 \epsilon_1\epsilon_2$, the central singularity has \emph{finite codimension}, given by the number of the unfolding parameters $u_k$ in the function $F$. This means that $F(x,y,u_k)$ is the simplest and most general function that can capture and describe the essential features of the system. In fact, the bifurcation diagram of $F$ in the $u_k$-parameters space, namely the change in the number of critical points, provides, via (B.5)-(B.7), the bifurcation curves (\ref{en_1u},\ref{en_2u},\ref{en_1l},\ref{en_2l}) and completely describe the qualitative features of the system. 

These qualitative predictions cannot change anymore by the addition of higher-order contributions, however they can gradually become quantitatively more accurate by considering higher-order terms up to some optimal order [Efthymiopoulos et al., 2004; Pucacco et al., 2008a].
Notice that $\mu^2 = 4\epsilon_1\epsilon_2$ just at the degenerate case $C=0$ corresponding to the uncoupled system. Moreover, $\epsilon_1,\epsilon_2$ are not well defined if $A=\pm C$. Therefore, these degenerate points give the value of the control parameters at the crossing points in the catastrophe map (Fig. \ref{cplot}) that are structurally unstable. In these degenerate cases, a similar reduction may still be possible; however, in analogy with the example of subsection \ref{degenz}, we may have to retain sixth degree (or even higher) order terms in the central singularity and its deformation. 


\nonumsection{Acknowledgments}

We acknowledge very fruitful discussions with A. Celletti, C. Efthymiopoulos, G. Gaeta, H. Han{\ss}mann and F. Verhulst. A.M. is supported by the European social fund  ``Support of inter-sectoral mobility and quality enhancement of research teams at Czech Technical University in Prague'' (CZ.1.07/2.3.00/30.0034). G.P. is partially supported by INFN, Sezione di Roma Tor Vergata, by the European MC-ITN grant ``Stardust'' and by GNFM-INdAM. We thank the anonymous reviewers for comments and suggestions which helped to improve the presentation.

\end{document}